\documentclass[a4paper,11pt]{article}
\usepackage{stmaryrd}
\usepackage{bbm}
\usepackage{float}
\usepackage{amscd}
\usepackage{mathrsfs}
\usepackage{latexsym,amsxtra}
\usepackage[dvips]{graphicx}
\usepackage[utf8]{inputenc}
\usepackage[T1]{fontenc}
\usepackage{lmodern}
\usepackage[all]{xy}
\usepackage{nicefrac,mathtools,enumitem}
\usepackage{microtype}
\usepackage{lipsum}
\usepackage{mathtools}
\usepackage{xfrac}

\usepackage{geometry}
\usepackage{hyperref}

\usepackage{times}
\usepackage{amsthm}
\usepackage{amsmath,scalerel}
\usepackage{amsfonts}
\usepackage{amssymb}
\usepackage{mathdots}
\usepackage{color}
\usepackage{mathrsfs}
\usepackage{stmaryrd}
\SetSymbolFont{stmry}{bold}{U}{stmry}{m}{n}

\usepackage{bm}

\usepackage{tikz-cd}
\usepackage{tikz}
\usetikzlibrary{matrix,arrows,decorations.pathmorphing}
\usepackage[all]{xy}
\usepackage{graphicx}
\usepackage{subfigure}
\usepackage[toc,page]{appendix}
\usepackage[usestackEOL]{stackengine}
\usepackage{dutchcal}

\usetikzlibrary{positioning}

\newcommand\keywords[1]{\textbf{Keywords}: #1}

\usepackage{theoremref}
\newtheorem{thm}{Theorem}[section]
\newtheorem{lem}[thm]{Lemma}
\newtheorem{cor}[thm]{Corollary}
\newtheorem{pro}[thm]{Proposition}

\theoremstyle{definition}
\newtheorem{ex}[thm]{Example}
\newtheorem{rmk}[thm]{Remark}
\newtheorem{defi}[thm]{Definition}
\newtheorem{nota}[thm]{Notation}
\numberwithin{equation}{section}

\usepackage{hyperref}
\hypersetup{
    colorlinks=true,
    linkcolor=blue,
    filecolor=magenta,    
    urlcolor=cyan,
    citecolor=cyan,
}

\newcommand{\be }{\begin{equation}}
\newcommand{\ee }{\end{equation}}

\newcommand{\br}[1]{   [ \cdot,    \cdot  ]   }

\newcommand{\sgn}{\mathrm{sgn}}

\newcommand{\diag}{\mathrm{diag}}

\newcommand{\ad}{\tmop{ad}}
\newcommand{\da}{\tmop{da}}

\catcode`\>=\active \def>{
\fontencoding{T1}\selectfont\symbol{62}\fontencoding{\encodingdefault}}
\catcode`\<=\active \def<{
\fontencoding{T1}\selectfont\symbol{60}\fontencoding{\encodingdefault}}

\usepackage{ulem}

\newcommand{\assign}{:=}
\newcommand{\mathd}{\mathrm{d}}
\newcommand{\mathe}{\mathrm{e}}
\newcommand{\mathi}{\mathrm{i}}
\newcommand{\nin}{\not\in}

\newcommand{\tmop}[1]{\ensuremath{\operatorname{#1}}}

\newcommand{\tmmathbf}[1]{\ensuremath{\boldsymbol{#1}}}

\newcommand{\udots}{{\mathinner{
\mskip1mu\raise1pt\vbox{\kern7pt\hbox{.}}
\mskip2mu\raise4pt\hbox{.}
\mskip2mu\raise7pt\hbox{.}\mskip1mu}}}

\newenvironment{prf}
    {\proof}
    {\hspace*{\fill}$\Box$\medskip}
    
\title{The Boundary Condition for Some Isomonodromy Equations}
\author{Qian Tang and Xiaomeng Xu}
\date{}

\newcommand{\Addresses}{{
  \bigskip
  \footnotesize
\noindent \textsc{Department of Mathematics, The University of Hong Kong, Hong Kong 999077, China}\par\nopagebreak
  \textit{E-mail address}: \texttt{919130201@qq.com}
}\\
\\
\footnotesize
\noindent \textsc{
School of Mathematical Sciences \& Beijing International Center
for Mathematical Research, Peking University, Beijing 100871, China}\par\nopagebreak
  \textit{E-mail address}: \texttt{xxu@bicmr.pku.edu.cn}
}

\begin{document}

\maketitle

\begin{abstract}
In this article, we study a special class of Jimbo-Miwa-Mori-Sato isomonodromy equations, which can be seen as a higher-dimensional generalization of Painlev\'e VI. We first construct its convergent $n\times n$ matrix series solutions satisfying certain boundary condition. We then use the Riemann-Hilbert approach to prove that the resulting solutions are almost all the solutions. Along the way, we find a shrinking phenomenon of the eigenvalues of the submatrices of the generic matrix solutions in the long time behaviour. 
\end{abstract}


\tableofcontents

\section{Introduction}
\subsection{Main result}
In this paper, we study the differential equation, called
the \textbf{$n$-th isomonodromy equation},
for a 
$n\times n$
matrix valued function $\Phi(u)$
with respect to $u_1,\ldots,u_n$,
\begin{eqnarray}
\label{isoeq}
\frac{\partial}{\partial u_k}
\Phi(u)=
[{\rm ad}^{-1}_{u}
{\rm ad}_{E_{k}}\Phi(u)
,\Phi(u)],
\quad 
\text{for every $k=1,\ldots,n$},
\end{eqnarray} 
where 
\begin{itemize}
\item $E_k=\tmop{diag} (0, \ldots ,
\underset{k \text{-th}}{1} 
, \ldots, 0)$;

\item $u=\diag(u_1,\ldots,u_n)$
does not belong to
the fat diagonal
\[\Delta=
\{(u_1,\ldots,u_n)\in \mathbb{C}^n
~|~
u_i = u_j, \text{for some $i\neq j$} \};\]

\item for any matrix $A$, $\ad_u^{-1} A$ is the matrix with the $(i,j)$-entry
\begin{align}
\label{Nota:adu-1}
(\ad_u^{-1} A)_{ij}\assign
\left\{\begin{array}{ll}
    \frac{1}{u_i - u_j} A_{i j} & ; i \neq j\\
    0 & ; i = j
  \end{array}\right..
\end{align}
\end{itemize}

The equation \eqref{isoeq} first appeared in the work of Jimbo-Miwa-Mori-Sato \cite{JMMS}, as the isomonodromy deformation equation of the $n\times n$ linear differential system with Poncar\'{e} rank $1$
\begin{align}
\label{ConfluBegin}
\frac{\mathd}{\mathd \xi} F(\xi)
&=
\left( u +
\frac{\Phi(u)}{\xi}\right)\cdot F(\xi).
\end{align}
Then following Miwa \cite{Miwa}, 
the solutions $\Phi(u)$ of \eqref{isoeq}
have the strong Painlev\'{e} property: 
they are multi-valued meromorphic functions 
of $u_1,\ldots,u_n$ and 
the branching occurs 
when $u$ moves along a loop around 
the fat diagonal $\Delta$.
Thus, according to the original idea of Painlev\'{e}, 
they can be a new class of special functions. Later on, it was shown by Harnad \cite{Harnad} that there is a duality of the equation that relates to the Schlesinger equations, and thus for $n=3$ \eqref{isoeq} is equivalent to the Painlev\'e VI equation. Therefore, it can be seen as a higher rank Painlev\'e equation.

Since then, the sources of the interests in \eqref{isoeq} become quite diverse: 
the particular case (with skew-symmetric $\Phi(u)$) was studied by Dubrovin \cite{Dubrovin} in relation to the Gromov-Witten theory, and in general (with $\Phi(u)\in \mathfrak{gl}_n(\mathbb{C})$) by Boalch \cite{Boalch3} in relation to complex reflections and Poisson-Lie groups; the system \eqref{isoeq} appeared in the work \cite{BTL} of Bridgeland and Toledano Laredo in relation to the stability conditions; And it is a multi-time dependent Hamiltonian system with time variables $u_1,\ldots,u_n$, whose quantization, following the work \cite{Reshetikhin1992} of 
Reshetikhin, is related to the Knizhnik--Zamolodchikov equation in the conformal field theory. 

Despite the many applications, the problem of determining the behavior of its solutions $\Phi(u)$ at the fixed critical singularities and the monodromy problems of the associated linear system \eqref{ConfluBegin} are only studied recently by the second author. In \cite{Xu}, the asymptotic behavior of skew-Hermitian valued solutions was studied, and problem for the general $\mathfrak{gl}_n$ valued solution is left open. 

In this paper, we give the asymptotic expansion of generic solutions $\Phi(u)$ of \eqref{isoeq} at a fixed singularity. Along the way, we find the following shrinking phenomenon of solutions $\Phi(u)$ in long time $u_1,\ldots,u_n$ behaviour. 
To state our main theorem, let us denote $\delta_k A$ as the
the upper left $k\times k$ submatrix and 
the diagonal part of a $n\times n$ matrix $A$, i.e.
\begin{align}\label{deltak}
  (\delta_k A)_{i j} 
  & \assign
  \left\{\begin{array}{ll}
    A_{i j} & ; 1 \leqslant i, j \leqslant k \quad
    \text{or}\quad i = j\\
    0 & ; \text{otherwise}
  \end{array}\right..
\end{align} 
And let us denote by $z^{{\rm ad}(X)}Y:=z^XYz^{-X}$ for any $n\times n$ matrices $X,Y$ and complex number $z$. 
\begin{thm}
\label{Thm:termwise}
For almost all the solutions $\Phi(u)=\Phi_n(u)$ 
of the $n$-th isomonodromy equation \eqref{isoeq},
there exists a sequence of $n \times n$ matrix-valued functions $\Phi_k(u_1,\ldots,u_k)$ for $k=1,\ldots,n$ and a constant $n\times n$ matrix $\Phi_0$ such that 
for every $1\leqslant k \leqslant n$
\begin{align} 
\label{limit1}
\underset{u_k \rightarrow \infty}{\lim} 
\delta_{k - 1} \Phi_k  = 
\delta_{k - 1} \Phi_{k - 1},\quad
\underset{u_k \rightarrow \infty}{\lim} 
\left(
\frac{u_k-u_{k-1}}{u_{k-1}-u_{k-2}}
\right)^{\tmop{ad} 
\delta_{k - 1} \Phi_{k - 1}} \Phi_k  =  \Phi_{k - 1},
\end{align}
and in the meanwhile $\Phi_0$ satisfies the boundary condition
\begin{equation}
\label{boundcondtion}
| \tmop{Re} (
\lambda^{(k-1)}_{i} - 
\lambda^{(k-1)}_{j}
) | < 1,
\quad
\text{for every $1\leqslant i,j\leqslant k-1$},
\end{equation}
where 
$\{\lambda^{(k-1)}_{i}\}_{i=1,\ldots,k-1}$
are the eigenvalues of the
upper left $(k-1)\times (k-1)$ submatrix
of $\Phi_{0}$,
and by convention $u_0:=0,u_{-1}:=-1$.

Conversely, for any constant matrix $\Phi_0$ 
that satisfies the boundary condition \eqref{boundcondtion},
there exists a unique solution 
$\Phi(u;\Phi_0)=\Phi_n(u)$
and a sequence of matrix functions $\Phi_{n-1},\ldots,\Phi_{1}$
such that \eqref{limit1} holds
for every $1\leqslant k \leqslant n$.
Furthermore,
each $\Phi_k$ has the series expansion
\begin{align*}
\Phi_k(u_1,\ldots,u_k) = \sum_{m=0}^\infty 
\left(
\frac{u_k-u_{k-1}}{u_{k-1}-u_{k-2}}
\right)^{-m} \phi_{k,m}(u_1,\ldots,u_k),
\end{align*}
given by the recursive relation 
\eqref{Form:dk-1pkm+1}, \eqref{Form:pkm+1}
with the initial value $\phi_{k, 0} = \left(
\frac{u_k-u_{k-1}}{u_{k-1}-u_{k-2}}
\right)^{-\tmop{ad} \delta_{k-1} \Phi_{k-1}} \Phi_{k-1}$.
\end{thm}

The solutions $\Phi(u;\Phi_0)$ in the generic class of Theorem \ref{Thm:termwise} are called
\textbf{shrinking solutions}.
The regularized limit $\Phi_0$ is called the \textbf{boundary value} of $\Phi(u;\Phi_0)$ (at the limit $\frac{u_{k+1}-u_{k}}{u_{k}-u_{k-1}}\rightarrow \infty$ for $k=1,\ldots,n$), and \eqref{boundcondtion} is called the \textbf{boundary condition}. 
See Definition \ref{Def:GS} for more details.
Following Theorem \ref{Thm:termwise}, the boundary value $\Phi_0$ gives a parameterization of the shrinking solutions of \eqref{isoeq}. 
\[\text{$\Big\{ \text{The space $\mathfrak{c}_0$ of boundary value } \Phi_0 \Big\}$} 
\xleftrightarrow{\text{Theorem \ref{Thm:termwise}}}
\text{\Big\{The class of shrinking solutions $\Phi(u;\Phi_0)$ \Big\} }\]
It should be emphasized that the above correspondence 
implicitly requires the information of 
$\arg\left(
\frac{u_k-u_{k-1}}{u_{k-1}-u_{k-2}}
\right)$, 
otherwise \eqref{limit1} or 
the initial value $\phi_{k, 0}$ would be meaningless.

\begin{rmk}
See Theorem \ref{Thm:GenisGen} 
for a explicit characterization of the shrinking solutions 
by the monodromy data of the corresponding linear systems. 
\end{rmk}
Let us explain the name shrinking solutions.
Note that the $\lambda^{(k-1)}_i$ 
in Theorem \ref{Thm:termwise} 
are also the eigenvalues of 
the upper left $(k-1)\times (k-1)$ submatrix
of $\Phi_{k-1}$. 
From \eqref{limit1} we see that
the boundary conditions \eqref{boundcondtion}
correspond to the shrinking phenomenon of the eigenvalues of the submatrices of the generic solutions in the long time limit $u_1\ll u_2\ll \cdots \ll u_n$. See Section \ref{SubSect:SolAsy} for more details. For example, given any shrinking solution $\Phi(u)$, let $\lambda^{(n-1)}_1(u)$ and 
$\lambda^{(n-1)}_2(u)$ be
any two eigenvalues of 
the upper left $(n-1)\times (n-1)$ submatrix of $\Phi(u)$,
by Theorem \ref{Thm:termwise} we have the shrinking phenomenon in the long time  $u_n\rightarrow\infty$ behavior
\begin{align*}
\underset{u_n\rightarrow\infty}
{\lim}
|{\rm Re}(
\lambda^{(n-1)}_1(u)-
\lambda^{(n-1)}_2(u))|<1.
\end{align*}

\begin{rmk}\label{catpt}
Since the isomonodromy equation is a $(u_1,\ldots,u_n)$ multi-time dependent Hamiltonian system, one can think of Theorem \ref{Thm:termwise} studies the long time $u_n\gg u_{n-1}\gg \cdots \gg u_1$ behaviour of its solutions $\Phi(u)$.
Geometrically the limit $\frac{u_{k+1}-u_{k}}{u_{k}-u_{k-1}}\rightarrow \infty$ for $k=2,\ldots,n$ is understood as a point, called the caterpillar point, on the De Concini-Procesi space that is a certain conmpactification of the parameter space $u\in \mathbb{C}^n\setminus \Delta$. See \cite{Xu} for more details. 
\end{rmk}

\begin{ex}
Let us illustrate the shrinking phenomenon using the numerical simulation. Let us take the case $n=3$, and the $3\times 3$ matrix solution $\Phi(u_1,u_2,u_3)$ with the randomly chosen initial point $u_1=1+{\rm i}, u_2=3+{\rm i}, u_3=7.65$ and initial value
\[\Phi(1+{\rm i}, 3+{\rm i}, 7.65)= \left(\begin{array}{ccc}
    1.35+3.46 {\rm i} & -1.48+0.09{\rm i} & 0.91+3.92 {\rm i}  \\
    1.49+1.75 {\rm i} & 0.48+6.48 {\rm i} & 7.42+2.48 {\rm i} \\
    2.40+2.06 {\rm i} & 1.04+0.08 {\rm i} & 9.16+0.84 {\rm i}
  \end{array}\right).\]
Denote by $\lambda^{(2)}_2$ and $\lambda^{(2)}_1$ the eigenvalues of its upper left $2\times 2$ submatrix. Then the numerical simulation shows that how ${\rm Re}(\lambda^{(2)}_2-\lambda^{(2)}_1)$ gradually becomes less than one as the time variable $u_3$ goes to infinite along real axis.   
\begin{figure}
\centering
\includegraphics[width=0.9\textwidth]{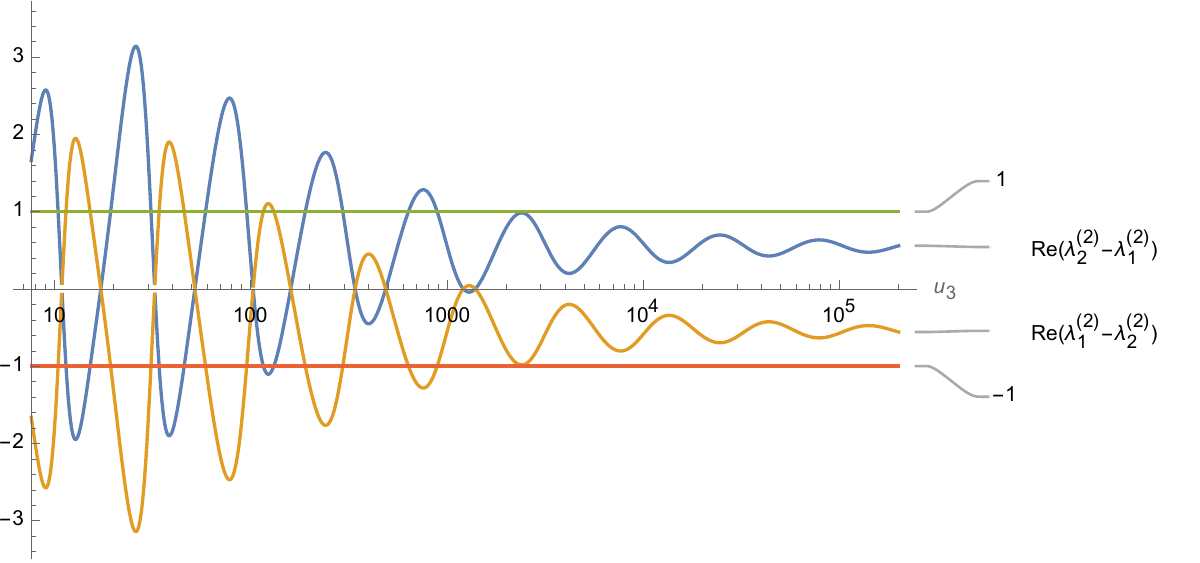}
\caption{Numerical simulation of the shrinking phenomenon when $n=3$}
\end{figure}
\end{ex}

From the numerical simulation of the typical example, we see that the difference of two eigenvalues has some oscillatory behavior, and the shrinking happens only after a quit long time in $u_3$. It implies that the shrinking phenomenon is a global behaviour, and indicates the complexity of Theorem \ref{Thm:termwise}. Accordingly, to study the second part of the theorem, i.e., the series expansion of solutions, we only use a local analysis. See Section \ref{Sect:Series}. And to derive the first part of the theorem, i.e., the boundary condition for generic solutions, we have to take the Riemann-Hilbert approach. We remark that the Riemann-Hilbert method is a powerful tool in the study of the global analysis of the nonlinear equation that can be realized as a compatibility condition of a linear system.

\subsection{The idea of Proof: the Riemann-Hilbert method}
Let us consider the $n\times n$ linear system of meromorphic differential equation for a function 
$F(z,u_1,\ldots,u_n)\in {\rm GL}_n(\mathbb{C})$ 
\begin{align}\label{introisoStokeseq1}
\frac{\partial F}{\partial z}&=\left( u +
\frac{\Phi(u;\Phi_0)}{z}\right)\cdot F(z,u),\\
\label{introisoStokeseq2}
\frac{\partial F}{\partial u_k}&=\left(E_kz+{\rm ad}^{-1}_u{\rm ad}_{E_{k}}\Phi(u;\Phi_0)\right)\cdot F, \ \text{for all} \ k=1,\ldots,n.
\end{align}
One checks that the compatibility condition of the linear system is the nonlinear equation \eqref{isoeq}. 

For any fixed $u$, the equation \eqref{introisoStokeseq1} 
has a unique formal solution $\hat{F}(z,u)$ 
around $z = \infty$. 
Then the standard resummation theory 
states that there exist certain sectorial regions around $z=\infty$, 
such that on each of these sectors there is a unique (therefore canonical) holomorphic solution with the prescribed asymptotics $\hat{F}(z,u)$. 
These solutions are in general different 
(that reflects the Stokes phenomenon), 
and the transition between them 
can be measured by a pair of Stokes matrices $S_\pm(u,\Phi(u;\Phi_0))$ (which are $(\mathe^{\mp\pi\mathi\delta\Phi}
S_{-\frac{\pi}{2}}^{\mp}(\mathi u,\Phi))^{\pm1}$ in the convention of Section \ref{subSect:BasicODE}, see Remark \ref{diffconv}).
Varying $u$, the Stokes matrices 
$S_\pm(u,\Phi(u;\Phi_0))\in{\rm GL}_n(\mathbb{C})$ of the system are locally constant (independent of $u$). 

Therefore, the Stokes matrices are the first integrals the nonlinear equation \eqref{isoeq}. Moreover, in terms of these integrals, the boundary values $\Phi_0$ of the solutions $\Phi(u;\Phi_0)$ of \eqref{isoeq} can be characterized by the following theorem (i.e., a solution to the Riemann-Hilbert problem).
We will provide its proof in Appendix \ref{Sect:Diag}.
\begin{thm}
\label{thm: introcatformula}\cite{Xu}
For all
purely imaginary parameters $u_1,\ldots,u_n$ 
with ${\rm Im}(u_1)<{\rm Im}(u_2)<\cdots <{\rm Im}(u_n)$
and generic $\Phi_0$, 
the sub-diagonals of the Stokes matrices $S_\pm(u,\Phi(u;\Phi_0))$ 
are given by
\begin{align*}
(S_+)_{k,k+1} &=
2\pi\mathrm{i}\cdot \mathe^{-\pi\mathrm{i}\cdot\frac{\small
{(\Phi_0)_{k k}+(\Phi_0)_{k+1, k+1}}}{2}} \\
&\quad \times 
\sum_{i=1}^k\frac{\prod_{l=1,l\ne i}^{k}
\Gamma(1+\lambda^{(k)}_i-\lambda^{(k)}_l)}{\prod_{l=1}^{k+1}
\Gamma(1+\lambda^{(k)}_i-\lambda^{(k+1)}_l)}\frac{\prod_{l=1,l\ne i}^{k}
\Gamma(\lambda^{(k)}_i-\lambda^{(k)}_l)}{\prod_{l=1}^{k-1}
\Gamma(1+\lambda^{(k)}_i-\lambda^{(k-1)}_l)}\cdot 
\det
(\lambda^{(k)}_iI-{\Phi_0})
^{1,\ldots,k-1,k}_{1,\ldots,k-1,k+1},\\
(S_-)_{k+1,k}& = 
-2\pi\mathrm{i}\cdot 
\mathe^{-\pi\mathrm{i}\cdot\frac{\small
{(\Phi_0)_{k+1, k+1}+(\Phi_0)_{k k}}}{2}}\\
&\quad \times 
\sum_{i=1}^k \frac{\prod_{l=1,l\ne i}^{k}\Gamma(1+\lambda^{(k)}_l-\lambda^{(k)}_i)}{\prod_{l=1}^{k+1}\Gamma(1+\lambda^{(k+1)}_l-\lambda^{(k)}_i)}\frac{\prod_{l=1,l\ne i}^{k}\Gamma(\lambda^{(k)}_l-\lambda^{(k)}_i)}{\prod_{l=1}^{k-1}\Gamma(1+\lambda^{(k-1)}_l-\lambda^{(k)}_i)}\cdot 
{\det
(\Phi_0-{\lambda^{(k)}_i}I)
^{1,\ldots,k-1,k+1}_{1,\ldots,k-1,k}}.
\end{align*}
where $k=1,\ldots,n-1$ and 
$A
^{a_1,\ldots,a_k}_{b_1,\ldots,b_k}$ is 
the $k\times k$ submatrix of $A$ formed 
by the $(a_1,\ldots,a_k)$ rows and 
$(b_1,\ldots,b_k)$ columns (here ${I}$ is the $n\times n$ identity matrix). 
Furthermore, the other entries are also given by explicit expressions.
\end{thm}
\begin{rmk}
From the viewpoint of integrable system, one can interpret the theorems as follows. The nonlinear equation \eqref{isoeq} is a multi-time dependent Hamiltonian system with the time variables $u_1,\ldots,u_n$, and has the Lax representation. The formula in Theorem \ref{thm: introcatformula} can be explicitly invertible on a dense subset of the space of Stokes matrices, that gives function $\Phi_0(S_{\pm})$, see in Section \ref{Sect:GeneralSol}. Thus Theorem \ref{Thm:termwise} and \ref{thm: introcatformula} state that the long time $u_1\ll u_2\ll \cdots \ll u_n$ behaviour of the solution $\Phi(u)$ of the nonlinear equation can be explicitly given by the "scattering data", i.e., Stokes matrices, at any given initial time. In literature, the long time behavior of a wide variety of integrable equations, like Korteweg de Vries (KdV), nonlinear Schr\"odinger, Painlev\'e transcendents, the Toda lattice, $tt^*$-equation and so on, have been studied by the Riemann-Hilbert
nonlinear steepest descent method of Deift and Zhou, see e.g., the works \cite{DZ1, DZ2, DIZ1, GIL1, GIL2}.
In this spirit, we solve the case for the nonlinear equation \eqref{isoeq}.
\end{rmk}

\begin{rmk}
Following \cite{Dubrovin}, the Witten-Dijkgraaf-Verlinde-Verlinde equations of 2D topological field theory or the semisimple Frobenius manifolds can be characterized by the system \eqref{introisoStokeseq1} and \eqref{introisoStokeseq2} with skew-symmetric matrix solutions $\Phi(u;\Phi_0)$. We expect that the theorems can find applications there, for example, in the study of Gamma conjecture, see e.g., \cite{Galkin2016},
\cite{cotti2018helix}.
\end{rmk}

The proof of our main Theorem \ref{Thm:termwise} essentially relies on the answer to the Riemann-Hilbert problem in Theorem \ref{thm: introcatformula}: starting from any constant $n\times n$ matrix $\Phi_0$ satisfying the condition \eqref{boundcondtion}, we construct a solution $\Phi(u;\Phi_0)$ given in a convergent series expansion form; to verify the solutions obtained in this way are almost all the solutions, we simply check that the corresponding Stokes matrices $S_\pm(u; \Phi(u;\Phi_0))$ already give almost all the possible Stokes matrices. 

To be more precise, the organization of the proof of Theorem \ref{Thm:termwise} or actually the paper, is as follows. 
In Section \ref{Sect:Series} we recursively construct 
series solutions $\Phi(u;\Phi_0)$ parameterized 
by constant matrices $\Phi_0$ that satisfy the boundary condition. We name the resulting series solutions $\Phi(u;\Phi_0)$ as
shrinking solutions. In Section \ref{Sect:GS},
we provide an effective criterion for the class of shrinking solutions, 
as well as the uniqueness part of Theorem \ref{Thm:termwise}.
In Section \ref{Sect:Cat}, 
we start from the boundary value $\Phi_0$ of 
a shrinking solution $\Phi(u;\Phi_0)$, 
and provide the explicit expression for the monodromy data of 
the associated linear system \eqref{introisoStokeseq1}-\eqref{introisoStokeseq2}, i.e., an analog of Theorem \ref{thm: introcatformula}. 
In Section \ref{Sect:GeneralSol}, 
based on the inverse of the explicit Riemann-Hilbert mapping 
introduced in Section \ref{Sect:Cat}, 
we solve the inverse monodromy problem, 
and prove that the monodromy data of the shrinking solutions 
from a dense subset of the monodromy data space. 
The using the property of the Riemann-Hilbert correspondence, we are able to show that the shrinking solutions give almost all the solutions (therefore are generic), and complete the proof of Theorem \ref{Thm:termwise}.

\[\begin{tikzcd}
{\text{Boundary values }\Phi_0} 
& [-13pt]
& [-13pt]{\text{Shrinking solutions }\Phi(u;\Phi_0)} \\
\\
& {\text{Open dense subset of monodromy data space}}
\arrow[dl,shift left,swap,
"{\text{Sect \ref{Sect:Series}, 
shrinking phenomenon \& series expansion}}"', 
from=1-1, to=1-3]
\arrow[ru,shift left,
"{\text{Sect \ref{Sect:GS}, criterion \& uniqueness}}", 
from=1-3, to=1-1]
\arrow["{\text{Sect \ref{Sect:GeneralSol}}}"', 
"\text{generic property}",
leftrightarrow, from=1-3, to=3-2]
\arrow["{\text{Sect \ref{Sect:Cat}}}", 
"\text{explicit RH mapping}"',
from=1-1, to=3-2]
\end{tikzcd}\]

\subsection{A new interpretatioin of Jimbo's formula for Painlev\'e VI}

We refer the reader to the book of Fokas, Its, Kapaev and Novokshenov \cite{Fokas2006} for a thorough introduction to the history and developments of the study of Painlev\'{e} equations. In particular, the sixth Painlev\'{e} equation (simply denoted by PVI or Painlev\'e VI) is the nonlinear differential equation
\begin{align*}
\frac{\mathd^2y}{\mathd x^2}
&=
\frac{1}{2}
\left(\frac{1}{y}+\frac{1}{y-1}+\frac{1}{y-x}\right)
\left(\frac{\mathd y}{\mathd x}\right)^2-
\left(\frac{1}{x}+\frac{1}{x-1}+\frac{1}{y-x}\right)
\frac{\mathd y}{\mathd x}\\
&\quad
+\frac{y(y-1)(y-x)}{x^2(x-1)^2}
\left(\alpha+\beta\frac{x}{y^2}+\gamma\frac{x-1}{(y-1)^2}
+\delta\frac{x(x-1)}{(y-x)^2}\right), 
\quad
\alpha,\beta,\gamma,\delta\in\mathbb{C}.
\end{align*}

It was shown by Harnad \cite{Harnad}, see also \cite{Mazzocco2002}, \cite[Section 3]{Boalch2005} and \cite{Degano2021} for a detailed way to do the Harnad duality, that Painlev\'{e} VI is equivalent to the equation \eqref{isoeq} with $n=3$ and suitable $\Phi(u)$. 
Therefore, the isomonodromy equation
can be seen as higher rank Painlev\'e transcendents. 
In particular, if
$\Phi(u_1, u_2, u_3)=(\Phi_{i j})_{3\times 3}$ is a solution
with eigenvalues $\lambda^{(3)}_1,
\lambda^{(3)}_2,
\lambda^{(3)}_3$, then 
\begin{align}\label{Phiy}
y(x) & = \frac{x}{x+(1-x)R(x)},
\quad x = \frac{u_2-u_1}{u_3-u_1},
\end{align}
with
\begin{subequations}
\begin{align}
\label{Form:R}
R(x) & = 
\frac{1}{v}\left(
\lambda^{(3)}_3-\Phi_{33}+
\frac{\Phi_{13}\Phi_{31}-
(\lambda^{(3)}_3-\Phi_{11})
(\lambda^{(3)}_3-\Phi_{33})}
{\lambda^{(3)}_3-\Phi_{11}+v}
\right),\\
\label{Form:PVIy}
v & = 
\frac{\Phi_{12}\Phi_{23} \Phi_{31} + 
(\lambda^{(3)}_2-\Phi_{22})
\Phi_{13} \Phi_{31}}
{(\lambda^{(3)}_2-\Phi_{11})
(\lambda^{(3)}_2-\Phi_{22})-
   \Phi_{12} \Phi_{21}},
\end{align}
\end{subequations}
is a solution of Painlev\'e VI with the parameters given by
\begin{subequations}
\begin{alignat}{4}\label{abcd}
\alpha & = 
\frac{(\theta_\infty-1)^2}{2},&\quad
\beta  & =-\frac{\theta_1^2}{2},&\quad
\gamma & = \frac{\theta_3^2}{2},&\quad
\delta & = \frac{1 -\theta_2^2}{2},\\
\label{PVI:theta}
\theta_\infty & = 
\lambda^{(3)}_1-\lambda^{(3)}_2,&\quad
\theta_1 & = \lambda^{(3)}_3 - \Phi_{11},&\quad
\theta_3 & = \lambda^{(3)}_3 - \Phi_{33},&\quad
\theta_2 & = \lambda^{(3)}_3 - \Phi_{22}.
\end{alignat}
\end{subequations}
According to Remark \ref{Rmk:reduced2var}, 
it can be verified that the left sides of 
\eqref{Form:R} and \eqref{Form:PVIy}
indeed only depend on $x$.

From the special case $n=3$ of 
Theorem \ref{kmainthm}, 
there exists a constant matrix
$\Phi_0=(\varphi_{i j})_{3\times 3}$ such that
\begin{subequations}
\begin{align}\label{3Phi02by2}
(u_2 - u_1)^{\tmop{ad} \delta_0 \Phi_0} \delta_2 \Phi_3 
& = \delta_2
\Phi_0 + O(x^{1-\sigma}),
\quad x\to 0,\\
\label{3Phi0}
x^{- \tmop{ad} \delta_2 \Phi_0} (u_2 - u_1)^{\tmop{ad} \delta_0 \Phi_0}\Phi_3 
& = \Phi_0 + O(x^{1-\sigma}),
\quad x\to 0,
\end{align}
\end{subequations}
where $\sigma$ is 
the difference $\lambda^{(2)}_1-\lambda^{(2)}_2$
of the eigenvalues of 
the $2\times 2$ upper left submatrix of $\Phi_0$,
and it satisfies
$0\leqslant {\rm Re} \sigma <1$.
Plugging \eqref{3Phi02by2}, \eqref{3Phi0} 
into the expression \eqref{Phiy}, by a direct computation we get
\begin{pro}\label{PhitoJ}
The solution $y(x)$ of Painlev\'e VI given by
$\Phi_3(u;\Phi_0)$ as in \eqref{Phiy}-\eqref{Form:PVIy}
has the following asymptotics
as $x=\frac{u_2-u_1}{u_3-u_1}\rightarrow 0$,
\begin{equation}
y(x)= Jx^{1-\sigma}+o(x^{1-\sigma}),\quad \sigma \neq 
0,
\theta_{\infty} + \theta_3,
\theta_{\infty} - \theta_3,
\end{equation}
where
\begin{align*}
(- \theta_3 + \theta_{\infty} - \sigma) (\theta_3 + \theta_{\infty} -
  \sigma) J 
& =  \varphi_{13} \varphi_{31} + \frac{2 (\varphi_{13}
  \varphi_{32} \varphi_{21} - \varphi_{12} \varphi_{23}
  \varphi_{31})}{\sigma}\\
&\quad\  - \frac{((\varphi_{11} - \varphi_{22}) \varphi_{13} + 2 \varphi_{12}
  \varphi_{23}) ((\varphi_{11} - \varphi_{22}) \varphi_{31} + 2 \varphi_{32}
  \varphi_{21})}{\sigma^2}.
\end{align*}
\end{pro}

Therefore, under the Harnad duality, the $n=3$ case of the Theorem \ref{Thm:termwise} recovers Jimbo's asymptotic formula \cite{Jimbo} for the generic solutions of Painlev\'{e} VI. In particular, the boundary condition \eqref{boundcondtion} for $n=3$ coincides with Jimbo's boundary condition $-1<\rm{Re} \sigma <1$ for Painlev\'{e} VI. 
Furthermore, under the Harnad duality, 
the Theorem \ref{thm: introcatformula}
is equivalent to 
Jimbo's monodromy formula \cite{Jimbo} for Painlev\'{e} VI, see our follow-up work \cite{WXX} for details.
It therefore justifies that the boundary value $\Phi_0$ in Theorem \ref{Thm:termwise} is the right parameterization of generic transcendents $\Phi(u)$, and it is very encouraging to see that the formula in Theorem \ref{thm: introcatformula} for $n=3$ case matches up with Jimbo's monodromy formula for Painlev\'{e} VI! 

We also remark that the correspondence from $\Phi_3$ to PVI $y$
given by \eqref{Form:R},\eqref{Form:PVIy} depends on the choices of the order of eigenvalues. The PVI $y$, corresponding to different choices of the ordering, differ by rational transformations, known as the Bäcklund transformations. 
The Bäcklund transformations of PVI \cite{inaba2004backlund,loray2016isomonodromic}
motivates us to introduce the Bäcklund transformations of the isomonodromy equation \eqref{isoeq}, see our next work.

As stressed in \cite{Fokas2006,IN,Conte2008}, 
the Painlev\'{e} transcendents are seen as nonlinear special functions, because they play the same role in nonlinear mathematical physics as that of classical special functions, like Airy functions, Bessel functions, etc., in linear physics. And it is the asymptotic formula at critical points, and expression of the monodromy of the associated $2\times 2$ or $3\times 3$ linear system, make Painlev\'{e} transcendents as efficient in applications as linear special functions (see \cite{Fokas2006, IN,Conte2008} and the references therein for more details). Therefore, we expect that in the higher rank cases, the asymptotic formula at critical points, and expression of the monodromy of the associated $n\times n$ linear system, i.e., Theorem \ref{Thm:termwise} and \ref{thm: introcatformula}, will find more applications, for example, in the study of special solutions of the equation \eqref{isoeq}. We remark that some important applications in representation theory and Poisson geometry have been found in \cite{Xu}.

\subsection{Further directions: the special solutions}
Let us discuss how to apply our results to find special solutions. 

\subsubsection*{Algebraic solutions}
Let us first recall the case $n=3$: based on Jimbo's formula, the algorithm in \cite{DM00, Boalch2} derives various algebraic solutions of the Painlev\'e VI. See \cite[Section 5]{Boalch2} for an important and detailed example of the algorithm. These in turn give algebraic solutions $\Phi(u)$ of \eqref{isoeq} in the case of rank $n=3$ under the Harnad duality. 

Parallelly, Theorem \ref{Thm:termwise} and \ref{thm: introcatformula} enable us to find the algebraic solutions of \eqref{isoeq} for general rank $n$ case in a systematic way. We know from Miwa’s theorem \cite{Miwa} that the matrix function $\Phi(u)$ is meromorphic on the universal covering of $\mathbb{C}^n\setminus \Delta$. Continuation along closed paths in the deformation space interchanges the branches of $\Phi(u)$. And such monodromy of the nonlinear isomonodromy equation is explicitly given in terms of the geometric terms, i.e., by an explicit braid group action on the corresponding Stokes matrices. See \cite{Dubrovin,BoalchG}.
In particular, a solution $\Phi(u)$ of the isomonodromy equation \eqref{isoeq} with the given Stokes matrix $S_\pm(u,\Phi(u))$ is an algebraic function, with branching along the diagonals $u_i=u_j$, if $S_\pm$ belong to a finite orbit of the action of
the pure braid group. See e.g., \cite[Appendix F]{Dubrovin}. 

Therefore, we have the following systematic way to find the algebraic solutions $\Phi(u)$,
\begin{equation}
    \Big\{\text{Stokes matrices $S_\pm(u,\Phi(u;\Phi_0))$}\Big\}\Longrightarrow \Big\{\text{boundary value} \ \Phi_0\Big\}\Longrightarrow \Big\{\text{solution} \ \Phi(u;\Phi_0) \Big\}.
\end{equation}
The algorithm is as follows: starting from such a pair of explicit $S_\pm$ that belong to a finite orbit, by Theorem \ref{thm: introcatformula} we can get the boundary value $\Phi_0$ of the series solution $\Phi(u;\Phi_0)$. Following Theorem \ref{Thm:termwise}, the boundary value $\Phi_0$ determines recursively the series expansion of the corresponding solution $\Phi(u;\Phi_0)$. Substituting back the leading term into the series expansion of the solutions of isomonodromy equation would determine, algebraically, any desired term in the multivariable Puiseux expansion, and then the solution itself. See our next work on the algebraic solutions $\Phi(u)$ for more details. In this paper, we give some examples of rational solutions whose corresponding Stokes matrices are identity matrix.

\subsubsection*{Non-generic solutions}
The shrinking solutions in Theorem \ref{Thm:termwise} are generic. However, there are non-generic solutions (solutions not in our generic class) that relate to other subjects. See Section \ref{subSect:NGS} for several examples of non-generic solutions of
the $n$-th isomonodromy equations. In particular, Proposition \ref{Pro:Rational} gives a new rational solution which is non-generic. 

Let us first recall that in $n=3$ case, there are important special cases of 
Painlev\'e VI transcendent $y(x)$, 
which are not covered by Jimbo's generic class, 
they were studied by 
Shimomura \cite{Sh1,Sh2}, 
Dubrovin-Mazzocco \cite{DM00}, 
Mazzocco \cite{Mazzocco2001}, 
Guzzetti \cite{Gu1,Gu2,Gu3}, 
Boalch \cite{Boalch3}, 
Kaneko \cite{Kaneko1,Kaneko2,Kaneko3} and 
Bruno-Goryuchkina \cite{Bruno2010} among many others.

Accordingly, in a next work, we will 
study the non-generic solutions $\Phi(u)$ of \eqref{isoeq}. 
Our starting point is Theorem \ref{thm: introcatformula}.
Note that the expression in Theorem \ref{thm: introcatformula} is an analytic function of $\Phi_0$ for
all $\Phi_0\in\mathfrak{gl}_n$ satisfying \eqref{boundcondtion}. 
Its poles along the boundary, when some of the inequalities in \eqref{boundcondtion} become equalities, is closely related to the non-generic solutions of the system \eqref{isoeq}. 
To be more precise, the boundary value $\Phi_0$ of 
a generic solution $\Phi(u)$ is understood as the initial value at 
the limit point $\frac{u_{k+1}-u_{k}}{u_{k}-u_{k-1}}\rightarrow \infty$. 
The formula, as a function of $\Phi_0$, 
has poles as $\lambda^{(k)}_l(\Phi_0)-\lambda^{(k)}_i(\Phi_0)=-1$, 
i.e., as the upper left $k\times k$ submatrix of $\Phi_0$ is resonant. 
These correspond to the non-generic solutions, 
whose initial condition posed at the limit point is not "smooth" 
but is a pole type. 
To distinguish the solutions that has a pole type, 
one needs to consider proper blow ups of the space 
$\mathfrak{c}_0$  of boundary values, i.e., 
certain blow ups along $\lambda^{(k)}_l-\lambda^{(k)}_i=-1$.
Thus we expect that there is a class of non-generic solutions $\Phi(u)$ that are geometrically parameterized by the boundary of a certain compactification $\widetilde{\mathfrak{c}_0 }$ of the space $\mathfrak{c}_0$. 
We expect that in $n=3$ case, it is equivalent to the compactification of 
the space of initial values in the study of Painlev\'e VI. 

More importantly,
the formula in Theorem \ref{thm: introcatformula} is 
deeply related to 
the linearization in Poisson geometry, 
Gelfand-Tsetlin theory and 
the representation theory of quantum groups 
(both crystal basis and at roots of unit) see \cite{Xu}, and 
the theory of spectral network and cluster algebras, see \cite{ANXZ}. 
We then expect the non-generic solutions, that correspond to the poles of the formula, 
have interesting and deep relations with the symplectic geometry of the coadjoint orbits and the representation theory of quantum groups at roots of unit.

\subsection*{Acknowledgements}
\noindent
Xiaomeng Xu is supported by the National Key Research and Development Program of China (No. 2021YFA1002000) and by the National Natural Science Foundation of China (No. 12171006).
Qian Tang is supported by the GRF grant no. 17303420 of the University Grants Committee of the Hong Kong SAR, China.

  \section{Series Solutions and the Shrinking Solution}
\label{Sect:Series}

To investigate the series solution of the isomonodromy equation, we introduce the following coordinates
\begin{align*}
z_1=\frac{u_1-u_0}{u_0-u_{-1}},\quad
z_2=\frac{u_2-u_1}{u_1-u_0},\quad
z_3=\frac{u_3-u_2}{u_2-u_1},\quad
\ldots\quad,
z_n=\frac{u_n-u_{n-1}}{u_{n-1}-u_{n-2}},
\end{align*}
where $u_0=0, u_{-1}=-1$.
Under this coordinates, 
the $n$-th isomonodromy equation \eqref{isoeq} 
becomes
\begin{equation}
\label{zisoeq}
  \frac{\partial}{\partial z_k} \Phi(z)
   =  
  [\tmop{ad}_u^{- 1}
  \tmop{ad}_{\partial u/\partial z_k} 
    \Phi(z)
  , \Phi(z)], \quad k = 1, \ldots, n,
\end{equation}
where the $n\times n$ diagonal matrix $u={\rm diag}(u_1,...,u_n)$ in the new coordinate is
\begin{subequations}
\begin{align}
\label{Expression:uz}
  u(z) 
  & = z_1 (I - \tilde{I}_0) + z_1 z_2 (I - \tilde{I}_1) + \cdots + z_1
  \cdots z_n (I - \tilde{I}_{n - 1}),\\
  \tilde{I}_m 
  & \assign 
  \tmop{diag} (\underset{m \text{-}
  \tmop{terms}}{\underbrace{1, \ldots, 1}}, 0, \ldots, 0).
\label{Def:tildeI}
\end{align}
\end{subequations}

In this section, for any $n\times n$ constant matrix $\Phi_0\in \mathfrak{c}_0$ the space of matrices satisfying the boundary condition \eqref{boundcondtion}, 
we will recursively construct 
a sequence of $n\times n$ matrix-valued functions
\[
\Phi_0\xrightarrow{\mathcal{s}_1}
\Phi_1(z_1)\xrightarrow{\mathcal{s}_2}
\Phi_2(z_1,z_2)\xrightarrow{\mathcal{s}_3}
\cdots \rightarrow 
\Phi_{n-1}(z_1,\ldots,z_{n-1})\xrightarrow{\mathcal{s}_n}
\Phi_n(z_1,\ldots,z_n),\]
where 
$\Phi_k$ satisfies the asymptotic behavior \eqref{limit1}
and the \textbf{$(k,n)$-equation}
\begin{eqnarray}
\label{kziso}
  \frac{\partial}{\partial z_j} \Phi_k(z_1,\ldots,z_k) 
  =
  [(\tmop{ad}_u^{- 1}
  \tmop{ad}_{\partial u/\partial z_j}
  - z_j^{- 1}) 
  \delta_k \Phi_k, \Phi_k],
  \qquad j=1,\ldots,k.
\end{eqnarray}
Here recall the operator $\delta_k$ is defined in \eqref{deltak}.
Note that for $k=n$, the equation \eqref{kziso} is just the $n$-th isomonodromy equation \eqref{zisoeq}. For $k<n$, it can be viewed as a degeneration of
$n$-th isomonodromy equation \eqref{zisoeq}: the upper left $k\times k$ submatrix of $\Phi_k$ simply satisfies the $k$-th isomonodromy equation.
\begin{rmk}
\label{Rmk:reduced2var}
In fact, by assign
$\tilde{\Phi}_k \assign
(z_1 z_2)^{\ad\delta\Phi_0}\Phi_k$
for $k>1$, 
we can eliminate the dependence on $z_1,z_2$. 
This is because 
the equation satisfied by $\tilde{\Phi}_k$ becomes
\begin{align*}
\label{tildekziso}
\frac{\partial}{\partial z_j} \tilde{\Phi}_k =
\left\{
\begin{alignedat}{2}
&0,&
& \quad j=1,2,\\
&[(\tmop{ad}_u^{- 1}
  \tmop{ad}_{\partial u/\partial z_j}
  - z_j^{- 1}) 
  \delta_k \tilde{\Phi}_k, \tilde{\Phi}_k],&
& \quad j>2,
\end{alignedat}
\right..
\end{align*}
\end{rmk}

To be more precise, the recursive step is as follows. Section \ref{step1} 
will start with an $n\times n$
constant matrix $\Phi_{k-1}$
as the leading term,
and then construct a formal series $\mathcal{s}_k(\Phi_{k-1})$
that satisfies the $(k,n)$-equation 
with respect to $z_k$. 
Section \ref{step2} will prove that for certain $\Phi_{k-1}$ (those satisfying condition \eqref{k-1bound}), the series solution $\Phi_k:=\mathcal{s}_k(\Phi_{k-1})$ is convergent. It also gives the explicit convergence radius of $\Phi_k=\mathcal{s}_k(\Phi_{k-1})$.
Section \ref{step3} will prove 
that when $\Phi_{k-1}$ further satisfies 
the $(k-1,n)$-equation, 
the constructed convergent series solution $\Phi_k=\mathcal{s}_k(\Phi_{k-1})$
will satisfy 
the entire $(k,n)$-equation. 
Finally, we will provide some examples and derive 
the following main theorem of this section
under Notation \ref{Not:Dk}.

\begin{thm}
\label{kmainthm}
For any $n\times n$ matrix valued solution
$\Phi_{k-1}$ of the $(k-1,n)$-equation that satisfies for all $i,j=1,...,k-1$,
\begin{equation}\label{k-1bound}
\left| \tmop{Re}(\lambda^{(k-1)}_i(\Phi_{k-1})-\lambda^{(k-1)}_j(\Phi_{k-1}))
\right|<1-\varepsilon,
\quad \text{with a real number } \varepsilon>0,
\end{equation}
the $(k,n)$-equation \eqref{kziso} has 
a convergent series solution $\Phi_k$ 
with respect to the variable $z_k$
\begin{equation*}
\Phi_k(z_1,\ldots,z_{k-1},z_k) = 
\sum_{m = 0}^{\infty} 
z_k^{- m} \phi_{k, m}(z_k;z_1,\ldots,z_{k-1}),
\end{equation*}
and $\phi_{k, m}$ has the following structure
\begin{align*}
\delta_{k-1}\phi_{k,m}(z_k)
& =
\sum_{d_m\in m\mathcal{D}_{k-1}(\Phi_{k-1})}
\sum_{t_m=0}^{2(\nu-1)m}
\beta_{d_m,t_m}(z_1,\ldots,z_{k-1})\cdot
z_k^{d_m} (\ln z_k)^{t_m},\\
z_k^{\ad \delta_{k-1}\Phi_{k-1}}\phi_{k,m} & =
\sum_{d_m\in m\mathcal{D}_{k-1}(\Phi_{k-1})}
\sum_{t_m=0}^{2(\nu-1)m}
\alpha_{d_m,t_m}(z_1,\ldots,z_{k-1})\cdot
z_k^{d_m} (\ln z_k)^{t_m}.
\end{align*}
\end{thm}
\begin{prf}
Combining the Proposition
\ref{Lem:TotallyFormalPhik},
\ref{pro:Convergent} and
\ref{Cor:ConstraPhikAna}
in Section \ref{step1}-\ref{step3},
\end{prf}

\subsection{Construction of Formal Series Solutions}
\label{step1} 

In this section, we
will start with an $n\times n$
constant matrix $\Phi_{k-1}$
as the leading term, and
construct a formal series $\mathcal{s}_k(\Phi_{k-1})$
that satisfies the $(k,n)$-equation only
with respect to $z_k$,
which is the ordinary differential equation
\begin{eqnarray}
\label{kzisoode}
  \frac{\mathd}{\mathd z_k} \Phi_k(z_k) =
  [(\tmop{ad}_u^{- 1}
  \tmop{ad}_{\partial u/\partial z_k}
  - z_k^{- 1}) 
  \delta_k \Phi_k(z_k), \Phi_k(z_k)].
\end{eqnarray}
\begin{nota}
\label{Not:Dk}
Denote the operators for 
any two $n\times n$ matrices $A$ and $X$
as follows
\begin{align*}
\tmop{ad}(A)X 
&\assign A X - X A ,\\
\tmop{da}(A)X 
&\assign A X + X A ,
\end{align*}
and introduce the following notations
\begin{subequations}
\begin{align}
\Phi^{[k]}
&\assign
\text{upper left $k\times k$ submatrix of $\Phi$},\\
\label{def:dkPhi}
\mathcal{D}_k(\Phi) 
& \assign
\{
\lambda^{(k)}_i-\lambda^{(k)}_j ~|~
\text{$\lambda^{(k)}_i, \lambda^{(k)}_j$
are the eigenvalues of $\Phi^{[k]}$
}\},\\
\label{def:mS}
m \mathcal{S}
& \assign
\{d_1+\cdots+d_m~|~ m\in \mathbb{N},
d_1,\ldots,d_m\in\mathcal{S}
\},\quad\text{for any $\mathcal{S}\subseteq\mathbb{C}$},
\end{align}
\end{subequations}
especially we have
$0 \mathcal{S} = \{0\}$.
\end{nota}
We will also 
frequently use the following properties
\begin{subequations} of the operators defined in \eqref{deltak}
\begin{align}
\label{dkFreq0}
\delta_{k_1} \delta_{k_2}
& =
\delta_{\min\{k_1,k_2\}},\\
\label{dkFreq1}
\delta_k ((\delta_k A_1) \cdot A_2)
& =
(\delta_k A_1) \cdot (\delta_k A_2),\\
\label{dkFreq2}
\delta_k ( A_1 \cdot (\delta_k A_2))
& =
(\delta_k A_1) \cdot (\delta_k A_2).
\end{align}
\end{subequations}

\begin{lem}
\label{Lem:FomExpand_adu}
There are formal power series expansions of the operators at $z_k=\infty$
\begin{subequations}
\begin{align}
\label{zexpank}
(\tmop{ad}_u^{- 1} 
\tmop{ad}_{\partial u/\partial z_k}
-z_k^{-1})
\delta_k
& = 
- z_k^{-1} \delta_{k-1}+
\sum_{m = 1}^{\infty} (- 1)^m 
z_k^{- (m + 1)}
\tmop{da} (U^m_{(k)}) \circ
(\delta_k - \delta_{k - 1}),\\
(\tmop{ad}_u^{- 1} 
\tmop{ad}_{\partial u/\partial z_j}
-z_j^{-1})
\delta_k
& = 
(\tmop{ad}_u^{- 1} 
\tmop{ad}_{\partial u/\partial z_j}
-z_j^{-1})
\delta_{k-1}
\nonumber\\
& \phantom{=}
- \sum_{m = 0}^{\infty} (- 1)^m 
z_k^{- (m + 1)}
\tmop{da}\left(
U^m_{(k)}\frac{\partial U_{(k)}}{\partial z_j}
\right) \circ
(\delta_k - \delta_{k - 1}),
\qquad j<k,
\label{zexpanj}
\end{align}
\end{subequations}
where 
\begin{align}
\nonumber
U_{(k)} & \assign
\tmop{diag} \left( 
\frac{u_{k - 1} - u_1}{u_{k - 1} - u_{k - 2}}
, \ldots, 
\frac{u_{k - 1} - u_{k - 3}}{u_{k - 1} - u_{k - 2}}
, 1, 0, 0, \ldots, 0 \right)\\
& =
\tilde{I}_{k-2} +
z_{k-1}^{-1} \tilde{I}_{k-3} +
\ldots +
z_{k-1}^{-1}\cdots z_3^{-1} \tilde{I}_{1},
\label{Def:Uk}
\end{align}
and $\tilde{I}_m$ is defined in \eqref{Def:tildeI}.
\end{lem}

\begin{prf}
From \eqref{Nota:adu-1},
for $1\leqslant j \leqslant k \leqslant n$ we have
\begin{eqnarray}\label{SpecPDE}
 \tmop{ad}_u^{- 1}\tmop{ad}_{E_j} \delta_k=
-\tmop{ad}_{D^{(n)}_j} \tmop{ad}_{E_j} \delta_k=
-\tmop{ad}_{D^{(k)}_j} \tmop{ad}_{E_j} \delta_k=
-\tmop{ad}_{D^{(k)}_j} \tmop{ad}_{E_j},
\end{eqnarray}
where
\begin{eqnarray*}
D^{(k)}_j \assign 
\tmop{diag} \left( \frac{1}{u_j - u_1}, \ldots, 0, \ldots,
\frac{1}{u_j - u_k}, 0, \ldots, 0 \right) .
\end{eqnarray*}
Since 
$\tmop{ad}_{E_m}\delta_k=0$ for $m>k$, 
from \eqref{Expression:uz} and \eqref{SpecPDE}
we have
\begin{eqnarray}
\label{TkExp}
\tmop{ad}_u^{- 1} 
\tmop{ad}_{\partial u/\partial z_k}
\delta_k
 = 
\frac{\partial u_k}{\partial z_k}\cdot
\tmop{ad}_u^{- 1} 
\tmop{ad}_{E_k} 
\delta_k
 = 
- z_1\cdots z_{k-1}\cdot
\tmop{ad}_{D_k^{(k)}} 
\tmop{ad}_{E_k} 
\delta_k.
\end{eqnarray}
Notice that
\begin{align}
\nonumber
z_1 \cdots z_{k - 1} D^{(k)}_k 
& = 
\tmop{diag} \left( 
\frac{1}{z_k + \frac{u_{k - 1} - u_1}{u_{k - 1} - u_{k - 2}}}
, \ldots, 
\frac{1}{z_k + \frac{u_{k - 1} - u_{k - 3}}{u_{k - 1} - u_{k - 2}}}, 
\frac{1}{z_k + 1},
\frac{1}{z_k}
, 0, \ldots, 0 \right)\\
\nonumber
& =
\tilde{I}_{k-1}/
(z_k^{-1}I + U_{(k)})\\
\label{z1zk-1Dkk:Expression}
& = 
z_k^{-1}\tilde{I}_{k - 1} + 
\sum_{m = 1}^{\infty} (- 1)^m z_k^{- (m + 1)} U_{(k)}^m.
\end{align}
Then replacing the corresponding term in \eqref{TkExp} by \eqref{z1zk-1Dkk:Expression}, and applying
\begin{align*}
  \tmop{ad} (U_{(k)}^m) \tmop{ad} (E_k) 
  & = 
  \tmop{da} (U_{(k)}^m) \tmop{ad}
  (\tilde{I}_{k - 1}) \tmop{ad} (E_k),\\
  \tmop{ad} (\tilde{I}_{k - 1}) \tmop{ad} (E_k) 
  & = 
  \delta_{k - 1} - \delta_k,
\end{align*}
we prove \eqref{zexpank}.
To prove \eqref{zexpanj}, 
first note that for $j<k$, 
we have
\begin{eqnarray}
\ad_u^{-1}\ad_{E_j} \delta_k
=
\ad_u^{-1}\ad_{E_j} \delta_{k-1} +
\frac{1}{u_k-u_j} \ad_{E_k}\ad_{E_j},
\end{eqnarray}
and hence
\begin{align}
\nonumber
\ad_u^{-1} 
\ad_{\partial u / \partial z_j} 
\delta_k 
& = 
\left( 
  \sum_{m = j}^{k - 1} 
  \frac{\partial u_m}{\partial z_j} 
  \tmop{ad}_u^{- 1} 
  \tmop{ad}_{E_m} \delta_k 
\right) + 
\frac{\partial u_k}{\partial z_j}
\tmop{ad}_u^{- 1} 
\tmop{ad}_{E_k} \delta_k\\
\nonumber
& =  
\sum_{m = j}^{k - 1} \left( \frac{\partial u_m}{\partial z_j}
  \tmop{ad}_u^{- 1} \tmop{ad}_{E_m} \delta_{k - 1} + \frac{1}{u_k - u_m}
  \frac{\partial u_m}{\partial z_j} \tmop{ad}_{E_k} \tmop{ad}_{E_m} \right) -
  \frac{\partial u_k}{\partial z_j} \tmop{ad}_{D^{(k)}_k} \tmop{ad}_{E_k}\\
\nonumber
& =  
(\tmop{ad}_u^{- 1} \tmop{ad}_{\partial u / \partial z_j} \delta_{k -
  1} - z_j^{- 1} \tmop{ad}_{\tilde{I}_{k - 1}} \tmop{ad}_{E_k}) -
  \frac{\partial u_k}{\partial z_j} \tmop{ad}_{D^{(k)}_k \tilde{I}_{j - 1}}
  \tmop{ad}_{E_k} + z_j^{- 1} \cdot \tmop{ad}_{\tilde{I}_{j - 1}}
  \tmop{ad}_{E_k}\\
\label{TjLeading}
& =  
(\tmop{ad}_u^{- 1} \tmop{ad}_{\partial u / \partial z_j} \delta_{k -
  1} + z_j^{- 1} (\delta_k - \delta_{k - 1})) + z_j^{- 1} \cdot \tmop{ad}_{(I
  - (u_k - u_{j - 1}) D^{(k)}_k) \tilde{I}_{j - 1}} \tmop{ad}_{E_k}.
\end{align}
Note that we also have the following identity from the compatibility condition
\begin{equation}\label{compkj}
\frac{\partial}{\partial z_k}
\left( \tmop{ad}_u^{- 1} 
\tmop{ad}_{\partial u/\partial z_j}
\delta_k \right)=
\frac{\partial}{\partial z_j} 
\left( \tmop{ad}_u^{- 1} 
\tmop{ad}_{\partial u/\partial z_k}
\delta_k \right).
\end{equation}
Plugging \eqref{TjLeading} into the left hand side of \eqref{compkj} and then comparing the coefficients of the $z_k^{-m}$ terms 
on both sides of \eqref{compkj} lead to \eqref{zexpanj}.
\end{prf}

Next, 
based on Lemma \ref{Lem:FomExpand_adu},
we introduce the formal series solution of \eqref{kzisoode}
in the form
$\Phi_k(z_k)=
\sum_{m = 0}^{\infty}
z_k^{- m} \phi_{k, m}(z_k)$.
Since each term $\phi_{k, m}(z_k)$
still depends on $z_k$, 
in order to provide 
the graded structure for this formal series, 
we need to introduce an additional 
free variable $\hbar$,
and consider the equation \eqref{hkzisoode}.

For simplicity, let us denote the operator 
\[\mathcal{T}_k(z_k) \assign
(\tmop{ad}_u^{- 1} \tmop{ad}_{\partial u/\partial z_k} -
z_k^{- 1}) \delta_k.\]
\begin{pro}
\label{Lem:TotallyFormalPhik}
For any
$n\times n$ constant matrix $\Phi_{k-1}$,
the ordinary differential equation
\begin{eqnarray}
\label{hkzisoode}
\frac{\hbar \mathd}{\mathd z_k} \Phi_k(z_k) = 
[\mathcal{T}_k(z_k/\hbar)
  \Phi_k(z_k), \Phi_k(z_k)],
\end{eqnarray}
has a formal series solution of the form
\begin{eqnarray}
\label{Series:Phik}
\Phi_k(z_k;\hbar) = 
\sum_{m = 0}^{\infty}
\hbar^{m}
z_k^{- m} \phi_{k, m}(z_k),
\end{eqnarray}
where the leading term
$\phi_{k, 0}(z_k) =
z_k^{-\tmop{ad} \delta_{k-1} \Phi_{k-1}} \Phi_{k-1}$,
and the coefficients are recursively given by 
\begin{subequations}
\begin{align}
\label{Form:dk-1pkm+1}
\delta_{k - 1} \phi_{k, m + 1}(z_k)
& =
z_k^{m+1} \int^{z_k} z_k^{-(m+1)}
b_{m+2}(z_k) 
\frac{\mathd z_k}{z_k},
\quad m \in \mathbb{N},\\
\label{Form:pkm+1}
z^{\tmop{ad} \delta_{k-1} \Phi_{k-1}}_k 
\phi_{k, m + 1} (z_k)
& =
z_k^{m+1} \int^{z_k} z_k^{-(m+1)} 
a_{m+2}(z_k) \frac{\mathd z_k}{z_k},
\quad m \in \mathbb{N}.
\end{align}
\end{subequations}
Here
\begin{subequations}
\begin{align}
\label{Form:bm}
b_{m + 2}(z_k) 
& \assign 
\sum_{p = 1}^{m+1} 
\sum_{i + j = m + 1 - p} (- 1)^p
\tmop{ad} (U_{(k)}^p) 
\Big(
(\delta_k-\delta_{k-1})\phi_{k,i}(z_k) \cdot
(\delta_k-\delta_{k-1})\phi_{k,j}(z_k)
\Big),\\
\label{Form:am}
a_{m + 2} (z_k)
& \assign 
\underset{i < m + 1}{\sum_{i + j = m + 1}}
z_k^{\tmop{ad} \delta_{k - 1} \Phi_{k - 1}} 
[\phi_{k, i}(z_k), \delta_{k - 1}
 \phi_{k, j}(z_k)]\nonumber\\
& \phantom{\assign} + 
\sum_{p = 1}^{m+1} 
\sum_{i + j = m + 1 - p} (- 1)^p
z_k^{\tmop{ad} \delta_{k - 1} \Phi_{k - 1}} 
\Big[\tmop{da} (U_{(k)}^p) 
(\delta_k-\delta_{k-1}) 
\phi_{k, i}(z_k), 
\phi_{k, j}(z_k)\Big].
\end{align}
\end{subequations}
\end{pro}

\begin{prf}
Suppose that the formal series \eqref{Series:Phik}
satisfies \eqref{hkzisoode}.
First we notice that,
from \eqref{dkFreq0} to \eqref{dkFreq2} we have
\begin{align}
\nonumber
\frac{\hbar\mathd}{\mathd z_k}
\delta_{k-1}\Phi_k 
& =
\delta_{k-1}
[\mathcal{T}_k(z_k/\hbar) \Phi_k,
\delta_k \Phi_k]\\
\nonumber
& =
\delta_{k-1}[\mathcal{T}_k(z_k/\hbar)
(\delta_k - \delta_{k-1}) \Phi_k,\delta_k \Phi_k]\\
& =
[\mathcal{T}_k(z_k/\hbar)
(\delta_k - \delta_{k-1}) \Phi_k,
(\delta_k - \delta_{k-1}) \Phi_k].
\label{Eq:dk-1pk}
\end{align}
Applying Lemma \ref{Lem:FomExpand_adu} and 
substituting the series \eqref{Series:Phik} into \eqref{Eq:dk-1pk},
and comparing the coefficients of 
the $\hbar^{m+2}$ terms,
we get
\begin{align}
\label{hseriesdk-10}
\frac{\mathd}{\mathd z_k}\delta_{k-1}\phi_{k,0}
&= 0,\\
\label{hseriesdk-1}
\frac{\mathd}{\mathd z_k}\left(
z_k^{-(m+1)}\delta_{k - 1} \phi_{k, m + 1}\right) 
&=
z_k^{-(m+2)}b_{m+2}.
\end{align}
Here the identity \eqref{hseriesdk-1} gives \eqref{Form:dk-1pkm+1}.

Similarly,
following Lemma \ref{Lem:FomExpand_adu}
\begin{align}
\nonumber
z_k^{-\ad\delta_{k-1}\Phi_{k-1}}
\frac{\hbar\mathd}{\mathd z_k}\left(
z_k^{\ad\delta_{k-1}\Phi_{k-1}}
\Phi_k \right) & = 
\hbar z_k^{-1}[\delta_{k-1}\Phi_{k-1}, \Phi_k] +
\frac{\hbar\mathd}{\mathd z_k} \Phi_k\\
\nonumber
& = 
\hbar z_k^{-1}[\Phi_k,\delta_{k-1}(\Phi_k-\Phi_{k-1})]\\
& \phantom{=}
+\sum_{m=1}^\infty (-1)^m 
\hbar^{m+1} z_k^{-(m+1)}
[\da(U_{(k)}^m)(\delta_k-\delta_{k-1})
\Phi_k,\Phi_k].
\label{Eq:zkadpk}
\end{align}
By substituting
series \eqref{Series:Phik} into \eqref{Eq:zkadpk} and then
comparing the coefficients of 
the $\hbar^{m+2}$ terms, 
we obtain
\begin{align}
\label{hseriesad0}
\frac{\mathd}{\mathd z_k}\left(
z_k^{\ad\delta_{k-1}\Phi_{k-1}}
\phi_{k,0}\right) 
&= 0,\\
\label{hseriesad}
\frac{\mathd}{\mathd z_k}\left(
z_k^{-(m+1)+\ad\delta_{k-1}\Phi_{k-1}}
\phi_{k, m + 1} \right) 
&=
z_k^{-(m+2)}a_{m+2}.
\end{align}
Here the identity \eqref{hseriesad} 
gives \eqref{Form:pkm+1}.
\end{prf}

\begin{rmk}
According to the equations \eqref{hseriesdk-10} and \eqref{hseriesad0}, 
all the formal solutions of \eqref{hkzisoode} 
in the forms of \eqref{Series:Phik} 
are given in this way.
In order to obtain 
a definite expression for \eqref{Series:Phik},
we adopt the following convention 
in \eqref{Form:dk-1pkm+1} and \eqref{Form:pkm+1}
to avoid the $\hbar^m z^0$ term
\begin{eqnarray}
\label{IntConvention}
\int^z z^s (\ln z)^t \frac{\mathd z}{z}
\assign
\left\{
\begin{alignedat}{2}
&\left(
\frac{\mathd}{\mathd s} \right)^t 
\left( \frac{z^{s}}{s} \right),&
& \quad t\in\mathbb{N}, s\neq 0,\\
&\frac{1}{t+1}
  (\ln z)^{t+1},&
& \quad t\in\mathbb{N}, s = 0,
\end{alignedat}
\right..
\end{eqnarray}
Note that the different choices of 
the lower limit of integration 
are equivalent to each other.
\end{rmk}

Letting $\hbar=1$ in 
Proposition \ref{Lem:TotallyFormalPhik} we introduce 
\begin{defi}
\label{Def:sk}
For any
$n\times n$ constant matrix $\Phi_{k-1}$,
denote the corresponding formal series solution of
the $(k, n)$-equation with respect to $z_k$, obtained in 
Proposition \ref{Lem:TotallyFormalPhik}, by
\begin{align}
\mathcal{s}_k(\Phi_{k-1}) & \assign
\sum_{m = 0}^{\infty}
z_k^{- m} \phi_{k, m}(z_k;\Phi_{k-1}),\\
\mathcal{s}_{k,\hbar}(\Phi_{k-1}) & \assign
\sum_{m = 0}^{\infty}
\hbar^m
z_k^{- m} \phi_{k, m}(z_k;\Phi_{k-1}).
\end{align}
\end{defi}

Let us analyze the structure of the series solution.
Roughly speaking, 
if we expand the terms
$\delta_{k-1}\phi_{k,m}(z_k)$ and
$z_k^{\ad \delta_{k-1}\Phi_{k-1}}
\phi_{k,m}(z_k)$
into a finite linear combination of complex powers of $z_k$,
we can verify that the powers of $z_k$ will 
fall into the subset $m\mathcal{D}_{k-1}(\Phi_{k-1})$ of complex number.
Therefore, 
if we impose 
the additional boundary condition \eqref{Cond:ConstraEigen}
on $\Phi_{k-1}$, then
the formal series
$\delta_{k-1}
\mathcal{s}_k(\Phi_{k-1})$ and
$z_k^{\ad \delta_{k-1}\Phi_{k-1}}
\mathcal{s}_k(\Phi_{k-1})$,
as complex power series of $z_k$, 
have a nicer structure:
\begin{itemize}
    \item the real part of the non-zero power 
    will be negative and decreasing to $-\infty$;

    \item the convention \eqref{IntConvention} 
    can be realized as a convergent definite integral 
    from $\infty$ to $z_k$.
\end{itemize}
Base on this structure, 
we prove
\begin{lem}
\label{Lem:FormalExpandPhik}
For any
$n\times n$ constant matrix $\Phi_{k-1}$
that satisfies
the following boundary condition
\begin{eqnarray}
\label{Cond:ConstraEigen}
\left| \tmop{Re}
\mathcal{D}_{k-1}(\Phi_{k-1})
\right|<1-\varepsilon,
\quad \exists\varepsilon>0,
\end{eqnarray}
under the notation 
in Proposition \ref{Lem:TotallyFormalPhik}
and convention \eqref{IntConvention}
we have
\begin{subequations}
\begin{align}
\label{Form:dk-1pkm+1finite}
\delta_{k - 1} \phi_{k, m + 1}(z_k)
& =
z_k^{m+1} \int^{z_k}_{\infty} z_k^{-(m+1)}
b_{m+2}(z_k) 
\frac{\mathd z_k}{z_k},\\
\label{Form:pkm+1finite}
z^{\tmop{ad} \delta_{k-1} \Phi_{k-1}}_k 
\phi_{k, m + 1} (z_k)
& =
z_k^{m+1} \int^{z_k}_{\infty} z_k^{-(m+1)} 
a_{m+2}(z_k) \frac{\mathd z_k}{z_k},
\end{align}
\end{subequations}
and $\phi_{k,m}(z_k)$ has 
the following structure
\begin{subequations}
\begin{align}
\label{Series:dk-1pkm}
\delta_{k-1}\phi_{k,m}(z_k)
& =
\sum_{d_m\in m\mathcal{D}_{k-1}(\Phi_{k-1})}
\sum_{t_m=0}^{2(\nu-1)m}
\beta_{d_m,t_m}\cdot
z_k^{d_m} (\ln z_k)^{t_m},\\
\label{Series:pkm}
z_k^{\ad \delta_{k-1}\Phi_{k-1}}
\phi_{k,m}(z_k) 
& =
\sum_{d_m\in m\mathcal{D}_{k-1}(\Phi_{k-1})}
\sum_{t_m=0}^{2(\nu-1)m}
\alpha_{d_m,t_m}\cdot
z_k^{d_m} (\ln z_k)^{t_m},
\end{align}
\end{subequations}
where $\nu$ is 
the nilpotency of 
the nilpotent part of $\delta_{k-1}\Phi_{k-1}$,
and $\beta_{d_m,t_m}$, $\alpha_{d_m,t_m}$ are 
the $n\times n$ constant matrices.
Here we recall that
$\mathcal{D}_{k-1}(\Phi_{k-1}),
m\mathcal{D}_{k-1}(\Phi_{k-1})$
is defined in \eqref{def:dkPhi}, \eqref{def:mS}.
\end{lem}
\begin{prf}
Let us denote 
the $\mathfrak{gl}_n(\mathbb{C})$-free module 
\begin{align*}
\mathcal{P}_m \assign
\tmop{span}\{
z_k^{d_m} (\ln z_k)^{t_m} ~|~
d_m\in m\mathcal{D}_{k-1}(\Phi_{k-1}),
0\leqslant t_m \leqslant 2(\nu-1)m
\}.
\end{align*}
Then we have
$\mathcal{P}_{m_1}\mathcal{P}_{m_2}
\subseteq
\mathcal{P}_{m_1+m_2}$.
We inductively assume that
for $s\leqslant m$,
\begin{subequations}
\begin{align}
\label{Control:dkpm}
\delta_{k-1}\phi_{k,s}(z_k)
&\in
\mathcal{P}_s,
\\
\label{Control:addkpm}
z_k^{\ad\delta_{k-1}\Phi_{k-1}}
\delta_{k-1}\phi_{k,s}(z_k)
&\in
\mathcal{P}_s,
\\
\label{Control:pm}
z_k^{\ad\delta_{k-1}\Phi_{k-1}}
\phi_{k,s}(z_k)
&\in
\mathcal{P}_s.
\end{align}
\end{subequations}
We will complete the induction and thus finish the proof
in the following way
\[
\begin{tikzcd}
&  & 
{\eqref{Control:dkpm},s\leqslant m+1}  \\
&  & 
{\eqref{Control:addkpm},s\leqslant m+1}\\
{\eqref{Control:pm},s\leqslant m} 
\arrow[rruu, Rightarrow]
\arrow[rru, Rightarrow] 
\arrow[rr, Rightarrow] 
&  & 
{\eqref{Control:pm},s\leqslant m+1}
\end{tikzcd}
\]
which is trivially valid 
for the base case $m=0$.

First we rewrite \eqref{Form:bm} as
\begin{align*}
b_{m + 2} 
& = 
\sum_{p = 1}^{m+1} 
\sum_{i + j = m + 1 - p} (- 1)^p
\tmop{ad} (U_{(k)}^p)\\
&\quad
z_k^{-\ad\delta_{k-1}\Phi_{k-1}}
\left(
(\delta_k-\delta_{k-1})
z_k^{\ad\delta_{k-1}\Phi_{k-1}}\phi_{k, i} \cdot
 (\delta_k-\delta_{k-1})
z_k^{\ad\delta_{k-1}\Phi_{k-1}}\phi_{k, j}
\right).
\end{align*}
According to the induction hypothesis
on \eqref{Control:pm} and
\begin{align*}
(\delta_k-\delta_{k-1})
z_k^{\ad\delta_{k-1}\Phi_{k-1}}\phi_{k, i} \cdot
 (\delta_k-\delta_{k-1})
z_k^{\ad\delta_{k-1}\Phi_{k-1}}\phi_{k, j}
& \in  \tmop{im}\delta_{k-1},\\
z_k^{-\ad\delta_{k-1}\Phi_{k-1}} 
\delta_{k-1}\mathcal{P}_m
& \subseteq 
\mathcal{P}_{m+1},
\end{align*}
we verify that
\begin{align}
\label{binP}
b_{m+2}(z_k)\in\mathcal{P}_{m+1}.
\end{align}
Hence by boundary condition \eqref{Cond:ConstraEigen}, 
the definite integral 
in \eqref{Form:dk-1pkm+1finite} 
absolutely converges, 
and \eqref{Control:dkpm} holds for $s=m+1$.

Next we rewrite 
\eqref{Form:dk-1pkm+1finite}, \eqref{Form:bm} as
\begin{align}
\label{Form:addphikm+1}
z_k^{\ad \delta_{k-1}\Phi_{k-1}}
\delta_{k-1}\phi_{k,m+1}
& =
z_k^{m+1+\ad \delta_{k-1}\Phi_{k-1}}
\int^{z_k}_\infty
z_k^{-(m+1+\ad \delta_{k-1}\Phi_{k-1})}
\left(
z_k^{\ad \delta_{k-1}\Phi_{k-1}} b_{m+2} \right)
\frac{\mathd z_k}{z_k},\\
z_k^{\ad \delta_{k-1}\Phi_{k-1}} b_{m+2} 
& =
\nonumber
\sum_{p = 1}^{m+1} 
\sum_{i + j = m + 1 - p} (- 1)^p
\tmop{ad} (z_k^{\ad\delta_{k-1}\Phi_{k-1}}U_{(k)}^p)\\
&\quad \left(
(\delta_k-\delta_{k-1})
z_k^{\ad\delta_{k-1}\Phi_{k-1}}\phi_{k, i} \cdot
 (\delta_k-\delta_{k-1})
z_k^{\ad\delta_{k-1}\Phi_{k-1}}\phi_{k, j}
\right).
\label{Form:adbm+2}
\end{align}
According to the induction hypothesis
on \eqref{Control:pm} and note that
\begin{eqnarray}
\label{adUinP1}
z_k^{\ad\delta_{k-1}\Phi_{k-1}}U_{(k)}^p
& \in &
\mathcal{P}_1,
\end{eqnarray}
we verify that
\begin{align}
\label{adbinP}
z_k^{\ad\delta_{k-1}\Phi_{k-1}} b_{m+2}(z_k)
\in\mathcal{P}_{m+1}.
\end{align}
Note that the following identity holds when the definite integral absolutely converges
\begin{align}
\nonumber
  z^{\mu + \tmop{ad} A} 
  \int_{\infty}^z z^{- (\mu + \tmop{ad} A)} \left(
  X z^s (\ln z)^t \right) \mathd z
  & = 
  \left( \frac{\mathd}{\mathd s} \right)^t 
  \left( z^{\mu + \tmop{ad} A}
  \int_{\infty}^z 
  z^{- (\mu + \tmop{ad} A)} (X z^s) \mathd z \right)\\
  & = 
  \left( \frac{\mathd}{\mathd s} \right)^t 
  \left( \frac{1}{s - \mu -
  \tmop{ad} A} (X z^s) \right),
\label{Form:ad(int)}
\end{align}
therefore we have
\begin{eqnarray}
  b_{m + 2} (z_k), z_k^{\tmop{ad} \delta_{k - 1} \Phi_{k - 1}} b_{m + 2} (z_k)
  \in \mathcal{P}_{m + 1} & \Rightarrow & z_k^{\tmop{ad} \delta_{k - 1}
  \Phi_{k - 1}} \delta_{k - 1} \phi_{k, m + 1} (z_k) \in \mathcal{P}_{m + 1},
\end{eqnarray}
and \eqref{Control:addkpm} holds for $s=m+1$.

Finally we rewrite \eqref{Form:am} as
\begin{align}
a_{m + 2} & =
\underset{i < m + 1}{\sum_{i + j = m + 1}}
[z_k^{\tmop{ad} \delta_{k - 1} \Phi_{k - 1}} \phi_{k, i}, 
 z_k^{\tmop{ad} \delta_{k - 1} \Phi_{k - 1}} 
 \delta_{k - 1}\phi_{k, j}]\nonumber
 + \sum_{p = 1}^{m+1} \sum_{i + j = m + 1 - p} \\
& \quad\
(- 1)^p
\Big[\tmop{da}(z_k^{\tmop{ad}\delta_{k-1}\Phi_{k-1}} U_{(k)}^p) 
(\delta_k-\delta_{k - 1}) 
z_k^{\tmop{ad} \delta_{k - 1} \Phi_{k - 1}} \phi_{k, i}, 
z_k^{\tmop{ad} \delta_{k - 1} \Phi_{k - 1}} \phi_{k, j}\Big].
\label{Form:Altam}
\end{align}
By using the induction hypothesis
on \eqref{Control:pm} for $s\leqslant m$,
on \eqref{Control:addkpm} for $s\leqslant m+1$,
and notice that \eqref{adUinP1},
we verify that $a_{m+2}(z_k)\in\mathcal{P}_{m+1}$. 
Hence by boundary condition \eqref{Cond:ConstraEigen}, 
the definite integral 
in \eqref{Form:pkm+1finite} 
absolutely converges, 
and \eqref{Control:pm} holds for $s=m+1$.
Thus we finish the proof.
\end{prf}
\begin{rmk}
Compared to $\delta_{k-1}
\mathcal{s}_k(\Phi_{k-1})$ and
$z_k^{\ad \delta_{k-1}\Phi_{k-1}}
\mathcal{s}_k(\Phi_{k-1})$,
the formal series
$\mathcal{s}_k(\Phi_{k-1})$ itself 
does not have the same nice structure, 
because of the subset of complex number
\begin{eqnarray*}
\mathcal{D}_n(\delta_{k-1}\Phi_{k-1})
\neq
\mathcal{D}_{k-1}(\delta_{k-1}\Phi_{k-1})
=
\mathcal{D}_{k-1}(\Phi_{k-1}).
\end{eqnarray*}
Moreover, 
the complex powers of $z_k$
in the expansion of 
$\delta_{k-1}
\mathcal{s}_k(\Phi_{k-1})$ and
$z_k^{\ad \delta_{k-1}\Phi_{k-1}}
\mathcal{s}_k(\Phi_{k-1})$
are specifically given by
the finite sum of the elements in the following set,
\begin{align}
\{-1+\lambda^{(k-1)}_i-\lambda^{(k-1)}_j~|~
\text{$\lambda^{(k-1)}_i, \lambda^{(k-1)}_j$
are the eigenvalues of $\Phi_{k-1}^{[k-1]}$
}\}.
\end{align}
When $\delta_{k-1}\Phi_{k-1}$ is diagonalizable, 
the $\ln z_k$ terms do not appear, thus the structure of the series expansion is simplified.
\end{rmk}

\begin{ex}
Suppose that $\Phi_{k-1}$ 
satisfies 
the boundary condition \eqref{Cond:ConstraEigen}.
In this case, the sub-leading term in the expansion of
$\Phi_k=\mathcal{s}_k(\Phi_{k-1})$ takes the explicit form
\begin{align}
\delta_{k-1}\phi_{k,1}(z_k) 
& =
\ad(U_{(k)})
\frac{z_k^{-\ad\delta_{k-1}\Phi_{k-1}}}
{1+\ad\delta_{k-1}\Phi_{k-1}}
((\delta_k-\delta_{k-1})\Phi_{k-1})^2,\\
\nonumber
z_k^{\ad\delta_{k-1}\Phi_{k-1}}\phi_{k,1}(z_k) 
& =
\Big[\Phi_{k-1},\ad\left(
\frac{z_k^{\ad\delta_{k-1}\Phi_{k-1}}}
{-1+\ad\delta_{k-1}\Phi_{k-1}}
U_{(k)}\right)
\frac{1}{1+\ad\delta_{k-1}\Phi_{k-1}}
((\delta_k-\delta_{k-1})\Phi_{k-1})^2\Big]\\
&\quad -\Big[\tmop{da}\left(
\frac{z_k^{\ad\delta_{k-1}\Phi_{k-1}}}
{-1+\ad\delta_{k-1}\Phi_{k-1}}
U_{(k)}\right)
(\delta_k-\delta_{k-1})
\Phi_{k-1},\Phi_{k-1}\Big].
\end{align}
\end{ex}

\subsection{Convergence of the Formal Series Solutions}
\label{step2}
We have seen that for those $\Phi_{k-1}$ satisfying the condition \eqref{Cond:ConstraEigen}, series $\mathcal{s}_k(\Phi_{k-1})$ is convergent.
In order to provide
the explicit convergence radius of 
the series $\mathcal{s}_k(\Phi_{k-1})$, 
we need the following lemma.

\begin{lem}
\label{Lem:RadofSeries}
Suppose that the complex numbers
sequences $(A_m)_{m\in\mathbb{N}}, (B_m)_{m\in\mathbb{N}^+}$ 
and complex number $\mu\in\mathbb{C}$
satisfy the following recurrence relation
\begin{align}
\label{AmDef}
(m+1)A_{m+1} & \leqslant 
C_a \sum_{i=1}^m A_i B_{m+1-i} +
{Q}   \sum_{p=1}^{m+1}\sum_{i=0}^{m+1-p}
\mu^p A_i A_{m+1-p-i},\\
\label{BmDef}
(m+1)B_{m+1} & \leqslant
C_b \sum_{p=1}^{m+1}\sum_{i=0}^{m+1-p}
\mu^p A_i A_{m+1-p-i},
\end{align}
for every $m\geqslant 1$, and assume that $C_aC_b\leqslant A_0^{-1}{Q}$
or $C_aC_b\leqslant A_0^{-2}$ holds, then we have
\begin{align*}
A_m & \leqslant 2^{m-1} \mu^m 
\max\{A_0,A_0^{m+1} {Q}^m\},\\
B_m & \leqslant C_b \cdot 2^{m-1} \mu^m 
\max\{A_0^2,A_0^{m+1} {Q}^{m-1}\}.
\end{align*}
\end{lem}

\begin{prf}
Define the sequences 
$(a_m)_{m\in\mathbb{N}}, (b_m)_{m\in\mathbb{N}^+}$ 
by assign
\begin{align*}
A_m & = a_m \mu^m
A_0 \max\{1,A_0 {Q}\}^m,\\
B_m & = b_m \mu^m
C_b A_0^2 \max\{1,A_0 {Q}\}^{m-1}.
\end{align*}
The we can reduce \eqref{AmDef}, \eqref{BmDef} to
\begin{align*}
(m+1)a_{m+1} & \leqslant
\sum_{i=1}^m a_i b_{m+1-i} +
\sum_{p=1}^{m+1}\sum_{i=0}^{m+1-p}
a_i a_{m+1-p-i},\\
(m+1)b_{m+1} & \leqslant
\sum_{p=1}^{m+1}\sum_{i=0}^{m+1-p}
a_i a_{m+1-p-i},
\end{align*}
with $a_0=1$.
Next we introduce the majorizing sequence
$(y_p)_{p\in\mathbb{N}}$ by
\begin{align*}
y_0 & = 1,\\
(m + 1) y_{m + 1} & = 
\sum_{i = 1}^m y_i y_{m + 1 - i} + 
\sum_{p = 1}^{m+1} \sum_{i = 0}^{m + 1 - p} 
y_i y_{m + 1 - p - i}.
\end{align*}
For $m\geqslant 1$ we check that 
\begin{eqnarray*}
2^{m-1}=y_m\geqslant\max\{a_m,b_m\}.
\end{eqnarray*}
Thus we finish the proof.
\end{prf}



Let us now provide an explicit bound
for the convergent radius of $\mathcal{s}_k(\Phi_k)$.
For this purpose,
we need to introduce
the Jordan-Chevalley decomposition of
$\delta_{k-1}\Phi_{k-1}$. Let us denote $\mathcal{N}$ as the nilpotent part 
of the $n\times n$ matrix $\delta_{k-1}\Phi_{k-1}$
with the nilpotency $\nu$, and denote $\mathcal{S}$ as the semisimple part
of $\delta_{k-1}\Phi_{k-1}$
with the spectral decomposition
\begin{eqnarray}
\mathcal{S}=
\lambda_1 P_1 +\cdots \lambda_{k-1}P_{k-1}
+ a_{k k} P_k +\cdots a_{n n} P_n,
\end{eqnarray}
where $P_j=E_j$ for $j\geqslant k$,
and $a_{jj}$ is the $(j,j)$-entry of $\Phi_{k-1}$.

\begin{pro}
\label{pro:Convergent}
For any
$n\times n$ constant matrix $\Phi_{k-1}$
that satisfies the boundary condition
\eqref{Cond:ConstraEigen},
the formal series $\mathcal{s}_k(\Phi_{k-1})$ 
is absolutely convergent when
\begin{align}\label{conradius}
|z_k|^\varepsilon 
& > 
\max\left\{
2\|U_{(k)}\|
\left(1+
\left(\frac{2 R \|\Phi_{k-1}\|}
{\varepsilon^\nu
 \varepsilon_0^{\nu-1}}\right)^2\right)
\left(1+\frac{1}{R \|\Phi_{k-1}\|}
\right),M
\right\},\quad
\arg z_k 
\in [\theta_1,\theta_2],
\end{align}
where $M$ is any given real number bigger than $\mathe^{2(\nu-1)}$,
matrix $U_{(k)}$ is defined in \eqref{Def:Uk}
from Lemma \ref{Lem:FomExpand_adu}, and
\begin{align*}
\varepsilon_0 & \assign
(1-\varepsilon)-
\max\left| \tmop{Re}
\mathcal{D}_{k-1}(\Phi_{k-1})
\right|,\\
R & \assign
\left(\sum_{i=1}^{n} \lVert P_i \rVert
\right)^2 \cdot
\left(\nu - \frac{1}{2}\right)
\left(1+
\left(4\nu \lVert \mathcal{N} \rVert
\right)^{2(\nu-1)}
\right) \cdot
\frac{1}
{1- \frac{2(\nu-1)}{\ln M} }.
\end{align*}
\end{pro}

\begin{prf}
For a fixed sufficiently large $N$ and 
constant $\varepsilon>0$ from the boundary condition
\eqref{Cond:ConstraEigen},
we introduce the following norm
for $n\times n$ matrix-valued function
$f(z_k) =
\sum_{i\in I} \alpha_i z_k^{d_i} (\ln z_k)^{t_i}$,
\begin{align}
\lVert f \rVert_m \assign
\inf\{C\geqslant 0 ~|~
\sum_{i\in I} \lVert\alpha_i \rVert
z_k^{d_i} (\ln z_k)^{t_i}
\leqslant C\cdot z_k^{m(1-\varepsilon)},
\text{ for all $z_k>M$}\},
\end{align}
where $I$ is a finite index set,
and $(d_i,t_i)\neq (d_j,t_j)$ for $i\neq j$.
It is direct to verify that
\begin{align}
\lVert f_1 f_2 \rVert_{m_1+m_2}
& \leqslant
\lVert f_1 \rVert_{m_1}\cdot
\lVert f_2 \rVert_{m_2}.
\end{align}
Define the sequences
$(A_m)_{m\in\mathbb{N}}$,
$(B_m)_{m\in\mathbb{N}}$ by
\begin{align*}
A_m & \assign 
\lVert
z_k^{\ad\delta_{k-1}\Phi_{k-1}}\phi_{k,m}(z_k)
\rVert_m,\\
B_m & \assign
\lVert
z_k^{\ad\delta_{k-1}\Phi_{k-1}}
\delta_{k-1}\phi_{k,m}(z_k)
\rVert_m.
\end{align*}
From Lemma \ref{Lem:FormalExpandPhik},
for every $m\in\mathbb{N}$
we have $A_m,B_m<\infty$ and
\begin{align}
\lVert
z_k^{\ad\delta_{k-1}\Phi_{k-1}}\phi_{k,m}(z_k)
\rVert & \leqslant
CA_m \cdot |z_k|^{m(1-\varepsilon)},
\quad\text{ for all $|z_k|^\varepsilon>M,
\arg z_k\in[\theta_1,\theta_2]$},\\
\lVert
z_k^{\ad\delta_{k-1}\Phi_{k-1}}
\delta_{k-1}\phi_{k,m}(z_k)
\rVert & \leqslant
CB_m \cdot |z_k|^{m(1-\varepsilon)},
\quad\text{ for all $|z_k|^\varepsilon>M,
\arg z_k\in[\theta_1,\theta_2]$}.
\end{align}
We will complete the proof 
by estimating the sequences
$(A_m)_{m\in\mathbb{N}}$,
$(B_m)_{m\in\mathbb{N}}$.

Let us denote
\begin{align}
z_k^{\ad\delta_{k-1}\Phi_{k-1}}
b_{m+2}(z_k) 
& =
\sum_{s\in I_{m+1}}
X_s \cdot z_k^{d_s} (\ln z_k)^{t_s},\quad
\text{$(d_i,t_i)\neq (d_j,t_j)$ for $i\neq j$},\\
a_{m+2}(z_k) 
& =
\sum_{s\in I_{m+1}}
Y_s \cdot z_k^{d_s} (\ln z_k)^{t_s},\quad
\text{$(d_i,t_i)\neq (d_j,t_j)$ for $i\neq j$}.
\end{align}
From \eqref{binP} and \eqref{adbinP} we know that
\begin{align}
\label{tsBound}
t_s &\leqslant (m+1)\cdot 2(\nu-1),\\
\label{dsSet}
d_s &\in (m+1)\mathcal{D}_{k-1}(\Phi_{k-1}),\\
\label{dsSet2}
d_s+\lambda_j-\lambda_i 
&\notin (m+1)\mathcal{D}_{k-1}(\Phi_{k-1})
\Rightarrow
P_i X_s P_j=0.
\end{align}
Thus by
\eqref{Form:addphikm+1} we have
\begin{align}
\nonumber
& \phantom{==}
z_k^{\tmop{ad} \delta_{k-1} \Phi_{k-1}} 
\delta_{k - 1} \phi_{k, m + 1}(z_k)\\
\nonumber & = 
\sum_{s \in I_{m+1} }
\left.
\left(\frac{\mathd}{\mathd x}\right)^{t_s}
\left(\frac{1}
{x-\ad\delta_{k-1}\Phi_{k-1}-(m+1)}
X_s z_k^x\right)
\right|_{x=d_s}\\
\nonumber & = 
\sum_{s \in I_{m+1} } 
\sum_{p=0}^{2(\nu-1)}
\left.
\left(\frac{\mathd}{\mathd x}\right)^{t_s}
\left(\frac{ (\ad\mathcal{N})^{p} }
{ (x-\ad\mathcal{S}-(m+1))^{p+1} }
X_s z_k^x\right)\right|_{x=d_s}\\
& = 
\sum_{s \in I_{m+1} } 
\underset{0\leqslant q\leqslant t_s}
{\sum_{0\leqslant p\leqslant 2(\nu-1)}}
\frac{(-1)^q (p+q)!}{p!q!}
\frac{t_s!}{(t_s-q)!(\ln z_k)^q}
\left(\frac{ (\ad\mathcal{N})^{p} }
{ (d_s-\ad\mathcal{S}-(m+1))^{p+q+1} }
X_s\right) 
z_k^{d_s}(\ln z_k)^{t_s}.
\end{align}
Following \eqref{tsBound}-\eqref{dsSet2} and setting
\begin{align}
R_\mathcal{S}\assign
(\sum_{i=1}^{n} \lVert P_i \rVert)^2,\quad
R_\mathcal{N} \assign
\left(\nu - \frac{1}{2}\right)
\left(1+
\left(4\nu \lVert \mathcal{N} \rVert
\right)^{2(\nu-1)}
\right),\quad
R_M \assign
\frac{1}
{1- \frac{2(\nu-1)}{\ln M} },
\end{align}
we can verify that for $|z_k|^\varepsilon>M>\mathe^{2(\nu-1)}$,
\begin{align}
\nonumber
& \phantom{\leqslant}
\underset{0\leqslant q\leqslant t_s}
{\sum_{0\leqslant p\leqslant 2(\nu-1)}}
\frac{(p+q)!}{p!q!}
\frac{t_s!}{(t_s-q)!(\ln z_k)^q}
\left\lVert\frac{ (\ad\mathcal{N})^{p} }
{ (d_s-\ad\mathcal{S}-(m+1))^{p+q+1} }
X_s\right\rVert\\
\nonumber
& \leqslant
\frac{ R_\mathcal{S} }{\varepsilon(m+1)}
\underset{0\leqslant q\leqslant 2(\nu-1)(m+1)}
{\sum_{0\leqslant p\leqslant 2(\nu-1)}}
\left(2\lVert \mathcal{N} \rVert
\frac{1+2(\nu-1)(m+1)}
{\varepsilon (m+1)}\right)^p
\left(\frac{2(\nu-1)(m+1)}
{\varepsilon (m+1)\ln z_k}\right)^q
\lVert X_s \rVert\\
\nonumber
& \leqslant
\frac{ R_\mathcal{S} }{\varepsilon(m+1)}
\underset{0\leqslant q\leqslant 2(\nu-1)(m+1)}
{\sum_{0\leqslant p\leqslant 2(\nu-1)}}
\left(
\frac{4\nu \lVert \mathcal{N} \rVert}
{\varepsilon}\right)^p
\left(\frac{2(\nu-1)}
{\ln M}\right)^q
\lVert X_s \rVert\\
\label{bNorm}
& \leqslant
\frac{ R_\mathcal{S} R_\mathcal{N} R_M}
{\varepsilon^{2\nu-1}(m+1)}
\lVert X_s \rVert.
\end{align}
Besides, from $1-(\varepsilon+\varepsilon_0)=
\max\left| \tmop{Re}
\mathcal{D}_{k-1}(\Phi_{k-1})
\right|$
and $\|\ln z_k\|_c\leqslant (c\cdot\mathe)^{-1}$,
we have
\begin{eqnarray}
\label{Uknorm}
\|z_k^{\ad\delta_{k-1}\Phi_{k-1}} U_{(k)}^p\|_{1}
\leqslant
\frac{R_\mathcal{S} R_\mathcal{N}}
{\varepsilon_0^{2(\nu-1)}}
\cdot\|U_{(k)}\|^p,
\end{eqnarray}
Combine \eqref{Form:adbm+2} with
\eqref{bNorm}, \eqref{Uknorm},
we deduce that
for every $m\in\mathbb{N}$
\begin{eqnarray}
\label{Est:zkad_dk-1_phikm+1}
\lVert
z_k^{\ad\delta_{k-1}\Phi_{k-1}}
\delta_{k-1}\phi_{k,m+1}(z_k)
\rVert_{m+1}
\leqslant
\frac{2R_\mathcal{S}^2 R_\mathcal{N}^2 R_M}{
\varepsilon^{2\nu-1}
\varepsilon^{2(\nu-1)}_0
(m+1)} 
\cdot \sum_{p=1}^{\infty} \sum_{i+j=m+1-p}
\lVert U_{(k)}\rVert^p A_i A_j,
\end{eqnarray}
and we finish the estimate of 
the sequence $(B_m)_{m\in\mathbb{N}}$.

By \eqref{Form:pkm+1finite} we have
\begin{align}
\nonumber
z_k^{\tmop{ad} \Phi_{k-1}} 
\delta_{k - 1} \phi_{k, m + 1}(z_k)
& = 
\sum_{s \in I_{m+1} }
\left.
\left(\frac{\mathd}{\mathd x}\right)^{t_s}
\left(\frac{X_s z_k^x}{x-(m+1)} \right)
\right|_{x=d_s}\\
\nonumber & = 
\sum_{s \in I_{m+1} } 
\sum_{q=0}^{t_s}
\frac{(-1)^q t_s!}{(t_s-q)!(\ln z_k)^q}
\frac{X_s}{(d_s-m-1)^{q+1}}
z_k^{d_s} (\ln z_k)^{t_s}.
\end{align}
Thus, we can verify that
\begin{eqnarray}
\| z_k^{\tmop{ad} \delta_{k - 1} 
\Phi_{k - 1}} \phi_{k, m + 1} (z_k) \|_{m+1}
\leqslant
\frac{R_M}{\varepsilon(m+1)}\cdot
\| a_{m+2} (z_k) \|_{m+1},
\end{eqnarray}
Combine \eqref{Form:Altam} with \eqref{Uknorm},
we have
\begin{align}
\| a_{m + 2} (z_k) \|_{m+1}
& \leqslant 
2 \|\Phi_{k - 1}\| \cdot 
\|
z^{\tmop{ad} \delta_{k - 1} \Phi_{k - 1}}_k \delta_{k - 1} \phi_{k, m + 1} (z_k) \|_{m+1} 
\nonumber\\
&\quad  +
\left( 2 \underset{0 < i < m + 1}{\sum_{i + j = m + 1}} 
A_i B_j + 
\frac{4 R_\mathcal{S} R_\mathcal{N}}
{\varepsilon_0^{2(\nu-1)}}
  \sum_{p = 1}^{\infty} \sum_{i + j = m + 1 - p} 
  \| U_{(k)} \|^p
  A_i A_j \right) .
\label{Est:a}
\end{align}
By substituting 
\eqref{Est:zkad_dk-1_phikm+1} into
\eqref{Est:a},
for every $m\in\mathbb{N}$
\begin{align}
\label{Bm}
(m + 1) B_{m + 1} 
& \leqslant
\frac{2R_\mathcal{S}^2 R_\mathcal{N}^2 R_M }{
\varepsilon^{2\nu-1}
\varepsilon^{2(\nu-1)}_0}
\sum_{p = 1}^{m+1} \sum_{i = 0}^{m+1-p}
\lVert U_{(k)}\rVert^p A_i A_{m+1-p-i},\\
\nonumber
(m + 1) A_{m + 1} 
& \leqslant
\frac{2 R_M}{\varepsilon} \sum_{i = 1}^m 
A_i B_{m + 1 - i} \\
&\quad + 
\frac{4 R_\mathcal{S} R_\mathcal{N} R_M }
{\varepsilon\cdot \varepsilon_0^{2(\nu-1)}}
\left(1+
\frac{R_\mathcal{S} R_\mathcal{N} R_M \|\Phi_{k-1}\|}
{\varepsilon^{2\nu-1}(m+1)}
\right)
\sum_{p = 1}^{m+1} \sum_{i = 0}^{m+1-p} 
\| U_{(k)} \|^p A_i A_{m + 1 - p - i}.
\label{Am}
\end{align}
In the end, let us use Lemma \ref{Lem:RadofSeries}: taking 
\begin{equation*}
C_a = \frac{2 R_M}{\varepsilon},\ \
C_b = \frac{2R_\mathcal{S}^2 R_\mathcal{N}^2 R_M}{
\varepsilon^{2\nu-1}
\varepsilon^{2(\nu-1)}_0},\ \
Q = \frac{4 R_\mathcal{S} R_\mathcal{N} R_M }
{\varepsilon\cdot \varepsilon_0^{2(\nu-1)}}
\left(1+
\frac{R_\mathcal{S} R_\mathcal{N} R_M \|\Phi_{k-1}\|}
{\varepsilon^{2\nu-1}}
\right),\ \
\mu = \| U_{(k)} \|,
\end{equation*}
in Lemma \ref{Lem:RadofSeries}
and noting that $C_a C_b \leqslant A_0^{-1} Q$ enable us to finish the proof.
\end{prf}

\begin{cor}
\label{Cor:AnaByPhi}
If $\Phi_{k-1}$ satisfies 
the boundary condition \eqref{Cond:ConstraEigen}, 
then $\Phi_k := \mathcal{s}_k(\Phi_{k-1})$ 
is analytic with respect to $\Phi_{k-1}$,
and we can recover
$\Phi_{k-1}$ from 
$\Phi_k = \mathcal{s}_k(\Phi_{k-1})$
by taking the limit
\begin{subequations}
\begin{align}
\label{Lim:dPk-1}
\lim_{z_k\to\infty}
\delta_{k-1}\Phi_k(z_k) 
& =
\delta_{k-1}\Phi_{k-1},\\
\label{Lim:Pk-1}
\lim_{z_k\to\infty}
z_k^{\ad\delta_{k-1}\Phi_{k-1}}
\Phi_k(z_k) 
& =
\Phi_{k-1}.
\end{align}
\end{subequations}
\end{cor}

\begin{prf}
It is straightforward to see that
$\phi_{k,m}(z_k)$ is analytic with respect to $\Phi_{k-1}$. Then according to Lemma \ref{pro:Convergent}, 
$\Phi_k = \mathcal{s}_k(\Phi_{k-1})$, 
as a series with respect to $\Phi_{k-1}$,
is locally uniformly convergent,
thus analytic with respect to $\Phi_{k-1}$.
The limits 
\eqref{Lim:dPk-1}, 
\eqref{Lim:Pk-1} are directly given by 
combining
\eqref{Series:dk-1pkm}, 
\eqref{Series:pkm} in Lemma \ref{Lem:FormalExpandPhik}
with Proposition \ref{pro:Convergent}.
\end{prf}

\subsection{Definition of the Shrinking Solutions}
\label{step3}
In this subsection, we will show that 
if $\Phi_{k-1}(z_1,\ldots,z_{k-1})$ 
is an analytic solution 
of the $(k-1,n)$-equation
and satisfies 
the boundary condition \eqref{Cond:ConstraEigen},
then the construed series $\mathcal{s}_k(\Phi_{k-1})$ 
will formally satisfy the entire $(k,n)$-equation,
not just with respect to $z_k$.
According to the convergence proved 
in Proposition \ref{pro:Convergent},
the series
$\mathcal{s}_k(\Phi_{k-1})$ eventually provides 
an analytic solution 
for the $(k,n)$-equation.
We will also provide 
the definition of the shrinking solution 
based on this construction process.

\begin{pro}
\label{Cor:ConstraPhikAna}
For any solution $\Phi_{k-1}$ of 
the $(k-1,n)$-equation that satisfies the boundary condition
\eqref{Cond:ConstraEigen},
the convergent series $\mathcal{s}_k(\Phi_{k-1})$
defined in Definition \ref{Def:sk}
is an analytic solution of 
the $(k,n)$-equation \eqref{kziso}.
\end{pro}

\begin{prf}
Denote the following operators 
for $1\leqslant j \leqslant k$ and
$\hbar \in \mathbb{C}\setminus\{0\}$,
\begin{align*}
\mathcal{T}_j(z_k) & \assign
(\ad_u^{-1}\ad_{\partial u/\partial z_j}-z_j^{-1})\delta_k,\\
\mathcal{T}_{j,\hbar}(z_k) & \assign
\mathcal{T}_j(z_k/\hbar).
\end{align*}
According to the compatibility of 
the $(k,n)$-equation \eqref{kziso}, 
which is
\begin{align*}
\frac{\partial \mathcal{T}_{k,\hbar}}{\partial z_j}
=
\frac{\hbar\partial \mathcal{T}_{j,\hbar}}{\partial z_k},\quad
\mathcal{T}_{j,\hbar}[\mathcal{T}_{k,\hbar} X,X] +
\mathcal{T}_{k,\hbar}[X,\mathcal{T}_{j,\hbar} X]
=
[\mathcal{T}_{k,\hbar} X,\mathcal{T}_{j,\hbar} X],
\end{align*}
we can verify that 
$\Phi_{k}=
\mathcal{s}_{k,\hbar}(\Phi_{k-1})$ satisfy
\begin{align*}
\frac{\hbar\partial}{\partial z_k} 
\left( \frac{\partial \Phi_k}{\partial z_j}
\right) 
& = 
\frac{\partial}{\partial z_j} 
\left( \frac{\hbar\partial}{\partial z_k} \Phi_k 
\right) \nonumber\\
& = 
\left[ \frac{\partial \mathcal{T}_{k,\hbar}}{\partial z_j} 
\Phi_k, \Phi_k \right] +
\left[ 
\mathcal{T}_{k,\hbar}
\frac{\partial \Phi_k}{\partial z_j}, \Phi_k \right] + 
\left[ 
\mathcal{T}_{k,\hbar}
\Phi_k, \frac{\partial \Phi_k}{\partial z_j} \right], \\
\frac{\hbar\partial}{\partial z_k} 
\left[\mathcal{T}_{j,\hbar}
\Phi_k, \Phi_k\right]
& = 
\left[
\frac{\hbar\partial \mathcal{T}_{j,\hbar}}
{\partial z_k} \Phi_k, \Phi_k \right] + 
\left[\mathcal{T}_{j,\hbar}
\left[\mathcal{T}_{k,\hbar}
\Phi_k,\Phi_k \right], \Phi_k\right] + 
[\mathcal{T}_{j,\hbar} \Phi_k, 
[\mathcal{T}_{k,\hbar} \Phi_k, \Phi_k]] \nonumber\\
& = 
\left[ \frac{\partial \mathcal{T}_{k,\hbar}}{\partial z_j} 
\Phi_k, \Phi_k \right] +
[\mathcal{T}_{k,\hbar} [\mathcal{T}_{j,\hbar}
\Phi_k, \Phi_k], \Phi_k] + 
[\mathcal{T}_{k,\hbar} \Phi_k, [\mathcal{T}_{j,\hbar} 
\Phi_k, \Phi_k]].
\end{align*}
Therefore, both
$\frac{\partial \Phi_k}{\partial z_j}$
and
$[\mathcal{T}_{j,\hbar}\Phi_k,\Phi_k]$,
as formal series $X(z_k) = 
\sum_{m=0}^{\infty} \hbar^m
z_k^{-m} x_m(z_k)$
with respect to $\hbar$,
satisfy the following equation formally
\begin{align}
\label{PdiffEq}
\frac{\hbar\partial}{\partial z_k} X
 =
\left[\frac{\partial \mathcal{T}_{k,\hbar}}{\partial z_j}
\Phi_k,\Phi_k\right] +
[\mathcal{T}_{k,\hbar} X
,\Phi_k] +
[\mathcal{T}_{k,\hbar}\Phi_k,X].
\end{align}
Note that $\Phi_{k-1}$ satisfies 
the $(k-1,n)$-equation,
from \eqref{zexpanj} in Lemma \ref{Lem:FomExpand_adu}
we have
\begin{align*}
\frac{\partial}{\partial z_j} \delta_{k - 1} \phi_{k, 0} 
& =
[\mathcal{T}_{j,0}
\phi_{k, 0}, \delta_{k - 1} \phi_{k, 0}], \\
\frac{\partial}{\partial z_j} 
z_k^{\tmop{ad} \delta_{k - 1} \phi_{k, 0}}
\phi_{k, 0} 
& =
[\mathcal{T}_{j,0} \phi_{k, 0}, 
z_k^{\tmop{ad} \delta_{k - 1} \phi_{k, 0}} \phi_{k, 0}].
\end{align*}
Therefore both
$\frac{\partial \Phi_k}{\partial z_j}$ and 
$[\mathcal{T}_{j,\hbar}\Phi_k,\Phi_k]$
have the same $0$-th term
\begin{eqnarray}
\label{y0is0}
x_0=
\frac{\partial \phi_{k, 0}}{\partial z_j}
= 
[\mathcal{T}_{j,0} \phi_{k, 0}, \phi_{k, 0}],
\end{eqnarray}
and satisfy the same equation \eqref{PdiffEq}.

Next we introduce the difference
\begin{align*}
Y(z_k;\hbar)=
\frac{\partial \Phi_k}{\partial z_j} -
[\mathcal{T}_{j,\hbar}\Phi_k,\Phi_k] =
\sum_{m=0}^{\infty} \hbar^m
z_k^{-m} y_m(z_k).
\end{align*}
From \eqref{PdiffEq}, we know that
$Y(z_k;\hbar)$ satisfy
the following equation formally
\begin{eqnarray}
\frac{\hbar\partial}{\partial z_k} Y  =  
[\mathcal{T}_{k,\hbar} Y, \Phi_k] + 
[\mathcal{T}_{k,\hbar} \Phi_k, Y].
\end{eqnarray}
Denote 
$\tilde{\mathcal{T}}_{k,\hbar}=
\hbar z_k^{- 1} \delta_{k - 1} + 
\mathcal{T}_{k,\hbar}$,
note that
\begin{align}
\frac{\hbar\partial}{\partial z_k} Y
& = 
[\tilde{\mathcal{T}}_{k,\hbar} Y, \Phi_k] +
[\tilde{\mathcal{T}}_{k,\hbar} \Phi_k, Y] - 
\hbar z_k^{- 1} 
\left([\delta_{k - 1} Y, \Phi_k] +
[\delta_{k - 1} \Phi_k, Y]\right), \\
\frac{\hbar\partial}{\partial z_k} \delta_{k - 1} Y
& = 
\delta_{k - 1}
\left(
[\tilde{\mathcal{T}}_{k,\hbar} Y, \Phi_k] + 
[\tilde{\mathcal{T}}_{k,\hbar} \Phi_k, Y]\right).
\label{RecdY}
\end{align}
According to \eqref{zexpank},
the lowest of $\tilde{\mathcal{T}}_{k,\hbar}$ is $\hbar^2$,
therefore from the equation \eqref{RecdY}
$\delta_{k - 1} y_{m + 1}$ can 
be recursively determined by
$y_m, \ldots, y_0$. Note that
\begin{align*}
{z_k^{- \tmop{ad} \delta_{k - 1} \phi_{k,0}}}  
\frac{\hbar \partial}{\partial z_k}
\left( {z_k^{\tmop{ad} \delta_{k - 1} \phi_{k,0}}}  Y \right) 
= [\mathcal{T}_{k,\hbar} Y, \Phi_k] + 
[\tilde{\mathcal{T}}_{k,\hbar} \Phi_k, Y], 
\end{align*}
thus $y_{m + 1}$ can 
be recursively determined by 
$\delta_{k - 1} y_{m + 1},
y_m,\ldots, y_0$.
Since $\Phi_{k-1}$ 
satisfies condition \eqref{Cond:ConstraEigen},
Lemma \ref{Lem:FormalExpandPhik} ensures that 
the term $\hbar^m z_k^0$ does not appear in $Y(z_k;\hbar)$, 
and \eqref{y0is0} ensures that $y_0=0$ (therefore recursively $y_1=0$ and so on),
thus we finally verify that formally
\begin{eqnarray}
\label{FormPhikisSol}
\frac{\partial \Phi_k}{\partial z_j}
=
[\mathcal{T}_{j,\hbar}\Phi_k,\Phi_k],\qquad
1\leqslant j <k,
\end{eqnarray}
for the series $\Phi_k=
\mathcal{s}_{k,\hbar}(\Phi_{k-1})$.

Denote $\mathcal{T}_{j,\hbar}=
\sum_{p=0}^\infty
\hbar^p \mathcal{T}_{j,\hbar}^{(p)}$
for $j<k$.
From \eqref{FormPhikisSol},
\eqref{zexpanj} and Lemma \ref{Lem:FormalExpandPhik}
we have
\begin{align*}
z_k^{\ad \delta_{k-1}\Phi_{k-1}}
\frac{\partial}{\partial z_j}
\phi_{k,m}(z_k) 
& = 
z_k^{\ad \delta_{k-1}\Phi_{k-1}}
\sum_{p+i_1+i_2=m}
[\mathcal{T}_{j,\hbar}^{(p)}\phi_{k,i_1}(z_k),
    \phi_{k,i_2}(z_k)]\\
& =
O(z_k^{(m+2)(1-\varepsilon)}),
\quad
z_k\to\infty.
\end{align*}
Therefore the series
$\sum_{m=0}^\infty z_k^{-m}
\frac{\partial \phi_{k,m}
}{\partial z_j}$
is also absolutely convergent
for sufficiently large $z_k$
by Proposition \ref{pro:Convergent}. It proves that $\Phi_k$ is analytical of 
the $(k,n)$-equation \eqref{kziso}.
\end{prf}

In fact, 
the eigenvalues of 
the upper left $(k-1)\times (k-1)$ submatrix 
of $\Phi_{k-1}$ are 
the trivial first integrals 
of the $(k-1,n)$-equation. 
Therefore, in Proposition \ref{Cor:ConstraPhikAna} the boundary condition \eqref{Cond:ConstraEigen} 
only needs to hold for certain 
$z_1,\ldots,z_{k-1}$. 
More precisely, we have the following property
\begin{pro}
\label{Pro:TrivalFirst}
The $(k, n)$-equation \eqref{kziso}
has the following first integrals
for solution $\Phi_k$,
\begin{itemize}
  \item entries of the diagonal part $\delta \Phi_k$;
  
  \item eigenvalues $(\lambda^{(m)}_j)_{j = 1, \ldots, m}$ of 
  $\Phi_k^{[m]}$, where $m \geqslant k$;
  
  \item entries 
  $((z_1\cdots z_k)^{\ad \delta\Phi_k}\Phi_k)_{i j}$, 
  where $i, j \geqslant k + 1$.
\end{itemize}
If $\Phi_{k-1}$ 
satisfies the boundary condition \eqref{Cond:ConstraEigen},
then $\Phi_k = \mathcal{s}_k(\Phi_{k-1})$ 
has the following properties
\begin{itemize}
\item $\delta \Phi_{k} =\delta \Phi_{k-1}$;

\item $\sigma(\delta_m\Phi_k)=\sigma(\delta_m\Phi_{k-1})$
for $m\geqslant k$;

\item $(z_k^{\ad\delta\Phi_{k}}\Phi_k)_{i j}=
(\Phi_{k-1})_{i j}$
for $i,j\geqslant k+1$.
\end{itemize}
\end{pro}
\begin{prf}
We only need to notice that for $m\geqslant k$,
the equation satisfied by $\delta_m\Phi_k$ 
\begin{eqnarray}
\label{kmziso}
\frac{\partial}{\partial z_j} 
X  = 
[(\tmop{ad}_u^{- 1}
  \tmop{ad}_{\partial u/\partial z_j}
  - z_j^{- 1}) 
  \delta_k X, X],
  \qquad j=1,\ldots,k,
\end{eqnarray}
is in the form of Lax pair,
thus the eigenvalues of $\delta_m\Phi_k$ 
are first integrals. 
The remaining results follow from a direct verification.
\end{prf}\\

It can be seen that when  $\Phi_{k-1}$ is a solution of the $(k-1,n)$-equation that satisfies the boundary condition \eqref{Cond:ConstraEigen} and $\Phi_k=\mathcal{s}_k(\Phi_{k-1})$, the the eigenvalues of 
the upper left $(k-1)\times (k-1)$ submatrix 
of $\Phi_{k}$ equal to the ones of $\Phi_{k-1}$.
Therefore,
based on Theorem \ref{kmainthm} and
Proposition \ref{Pro:TrivalFirst}, 
we can introduce the following definition.
\begin{defi}
\label{Def:GS}
A solution $\Phi(u;\Phi_0)=
\Phi_n(z)$ of 
the $n$-th isomonodromy equation \eqref{isoeq} 
is called a \textbf{shrinking solution} if 
\[\Phi_n=
\mathcal{s}_n\circ \mathcal{s}_{n-1}\circ\cdots \circ
\mathcal{s}_1(\Phi_0)\]
for some 
$n\times n$ constant matrix $\Phi_0$ 
that satisfies 
the \textbf{boundary conditions} \eqref{boundcondtion},
i.e.,
\begin{equation}
\label{boundcondtionNew}
| \tmop{Re} (
\lambda^{(k)}_{i}(\Phi_0) - 
\lambda^{(k)}_{j}(\Phi_0)
) | < 1,
\quad
\text{for every $1\leqslant i,j\leqslant k\le n-1$},
\end{equation}
and $\Phi_0$ is called 
the \textbf{boundary value}
of $\Phi(u;\Phi_0)$
at the boundary point
$z_1,\ldots,z_n\to\infty$.
\end{defi}

\begin{ex}
\label{Ex:PhikSeries}
Boundary condition \eqref{Cond:ConstraEigen}
always holds for $\Phi_{k-1}$ automatically for $k=1,2$,
thus
the series solutions are convergent for 
$(1,n)$-equation and 
$(2,n)$-equation. They have the explicit expression
\begin{align}
\Phi_1(z_1)=
\mathcal{s}_1(\Phi_0)
& = 
z_1^{-\ad\delta_0\Phi_0}\Phi_0,\\
\Phi_2 (z_1,z_2)=
\mathcal{s}_2(\Phi_1)
& = 
z_2^{-\ad\delta_1\Phi_1(z_1)}\Phi_1(z_1)
  = (z_1z_2)^{-\ad\delta_0\Phi_0}\Phi_0.
\end{align}
They also
provide all the solutions to those equations. 
\end{ex}

\begin{ex}[Rational Solution \cite{Dubrovin}]
\label{Ex:Rational}
Suppose that $\delta_{n-1}\Phi_{0}=0$, and
\begin{eqnarray}
\sum_{j=1}^{n-1}
(\Phi_{0})_{n j}
(\Phi_{0})_{j n} = 1.
\end{eqnarray}
Then there exists two vectors
\[ a = \left(\begin{array}{c}
     a_1\\
     \vdots\\
     a_{n - 1}\\
     1
   \end{array}\right), \quad b = \left(\begin{array}{cccc}
     b_1 & \cdots & b_{n - 1} & 1
   \end{array}\right), \]
satisfying $b a = 0$, 
such that $\Phi_{0} = \tmop{ad}_{E_n} (a b)$.
It can be verified that we have 
$\Phi_k = \Phi_{0}$ for $k\leqslant n-1$.
Furthermore, we have
\begin{align*}
  \Phi_n (u) & =  \frac{1}{b u a} \tmop{ad}_u (a b),\\
  \Phi_n (z) & =  \frac{1}{1 - (b U_{(n)} a) \cdot z_n^{- 1}} \cdot
  (\tmop{ad}_{E_n} (a b) - z_n^{- 1} \tmop{ad}_{U_{(n)}} (a b)),
\end{align*}
and $\phi_{n, m}$ is independent of $z_n$,
\begin{align*}
  \phi_{n, 0} & =  \tmop{ad}_{E_n} (a b),\\
  \phi_{n, m} & =  (b U_{(n)} a)^m \tmop{ad}_{E_n} (a b) - (b U_{(n)} a)^{m -
  1} \tmop{ad}_{U_{(n)}} (a b),
  \quad m\geqslant 1.
\end{align*}
\end{ex}

\section{Uniqueness of shrinking solutions}
\label{Sect:GS}
In this section,
we will prove the uniqueness of $\Phi_k$ in
\eqref{limit1}
for a given $\Phi_{k-1}$. 
This property also provides an effective criterion 
for the shrinking solution. 

\begin{thm}
\label{Thm:AsySol}
Suppose that
$\Psi_n(z)$ is a solution of
the $n$-th isomonodromy equation \eqref{zisoeq},
and there exists
$\Psi_{n-1},\ldots,\Psi_{0}$
such that for every $1\leqslant k \leqslant n$
we locally uniformly have
\begin{subequations}
\begin{align} 
\label{limit1S4}
\underset{z_k \rightarrow \infty}{\lim} 
\delta_{k - 1} \Psi_k & = 
\delta_{k - 1} \Psi_{k - 1},\\ 
\label{limit2S4}
\underset{z_k \rightarrow \infty}{\lim} 
z_k^{\tmop{ad} 
\delta_{k - 1} \Psi_{k - 1}} \Psi_k & =  \Psi_{k - 1},
\end{align}
\end{subequations}
with $\Psi_0\in\mathfrak{c}_0$,
then $\Psi_n(z)$ must coincide with 
the shrinking solution
with the boundary value $\Psi_0$ 
constructed in Section \ref{Sect:Series}.
\end{thm}

\begin{prf}
It follows from Proposition \ref{Cor:GSolutionAna}
and the definition of the shrinking solution.
\end{prf}

\begin{nota}\label{bcs}
For the solution $\Psi_k(z_k)$ of 
the $(k, n)$-equation \eqref{kzisoode}
with respect to $z_k$, denote
\begin{subequations}
\begin{align}
\label{lkdP}
    \mathcal{l}_k (\delta_{k - 1} \Psi_k) 
    & \assign 
    \underset{z_k \rightarrow \infty}{\lim} 
    \delta_{k - 1} \Psi_k, \\
    \mathcal{l}_k (\Psi_k) 
    & \assign 
    \underset{z_k \rightarrow \infty}{\lim}
    z_k^{\tmop{ad} \mathcal{l}_k 
    (\delta_{k - 1} \Psi_k)} \Psi_k,
\label{lkP}
\end{align}
\end{subequations}
if the limits exists
\footnote{
It is clear that 
$\mathcal{l}_k (\delta_{k - 1} \Psi_k) =
\delta_{k - 1} \mathcal{l}_k(\Psi_k),$
when 
$\mathcal{l}_k (\Psi_k)$ exists 
(this statement implies that 
$\mathcal{l}_k (\delta_{k - 1} \Psi_k)$ exists).}.
For various solution spaces,
we denote
\begin{itemize}
\item $\mathfrak{c}_k$ as
the space of the
$\mathfrak{gl}_n(\mathbb{C})$-valued functions
$\mathfrak{c}_k \assign
\mathcal{s}_k\circ
\mathcal{s}_{k-1}\circ
\cdots\circ
\mathcal{s}_1(\mathfrak{c}_{0})$,
where $\mathfrak{c}_0\subset \mathfrak{gl}_n(\mathbb{C})$ 
is the set of matrices
that satisfies the boundary conditions
\eqref{boundcondtionNew}.
Specifically, $\mathfrak{c}_n$
is the space of the shrinking solution;

\item $\mathfrak{s}_k$ as 
    the space of solutions of the $(k,n)$-equation \eqref{kziso}
    that have locally uniform limit
    \eqref{lkdP}, \eqref{lkP}, 
    and satisfy the \textbf{shrinking condition}
\begin{align}
\label{shrink}
\underset{z_k \to \infty}
{\overline{\lim}}
\left| \tmop{Re}
\mathcal{D}_{k-1}(\Psi_{k})
\right|<1-\varepsilon,
\quad \varepsilon>0.
\end{align}
Denote $\mathfrak{s}_k(\bm{z}_{k-1})$ as
the space of solutions of the $(k,n)$-equation \eqref{kzisoode}
with respect to $z_k$
that have limits
\eqref{lkdP}, \eqref{lkP},
with fixed values
$\bm{z}_{k-1} = (z_1,\ldots,z_{k-1})$;

\item $\mathfrak{b}_{k-1}$ as 
the space of solutions of the $(k-1,n)$-equation 
that satisfy the boundary condition
\eqref{Cond:ConstraEigen}.
Denote $\mathfrak{b}_{k-1} (\tmop{pt})$ as 
the set of $n\times n$ constant matrix
that satisfy condition
\eqref{Cond:ConstraEigen}.
\end{itemize}
\end{nota}
\noindent
From \eqref{Lim:dPk-1}, \eqref{Lim:Pk-1}
in Corollary \ref{Cor:AnaByPhi},
one sees that
$\mathcal{l}_k\mathcal{s}_k=
\tmop{id}_{\mathfrak{b}_{k-1}(\mathrm{pt})}$,
therefore $\mathfrak{b}_{k-1}(\tmop{pt})
\xrightarrow{\mathcal{s}_k}
\mathfrak{s}_k(\bm{z}_{k-1})$ is an injection.

Section \ref{SubSect:SolAsy} 
will give a characterization of 
the solution space $\mathfrak{s}_k(\bm{z}_{k-1})$.
Section \ref{SubSect:EffCri} 
will prove that 
$\mathfrak{s}_k(\bm{z}_{k-1})
\xrightarrow{\mathcal{l}_k}
\mathfrak{b}_{k-1}(\tmop{pt})$
is also an injection, 
thus induces the inverse of 
$\mathfrak{b}_{k-1} \xrightarrow{\mathcal{s}_k}
\mathfrak{s}_k$, 
which implies 
the uniqueness in Theorem \ref{Thm:termwise}.
As an application, 
we will see that every skew-Hermitian solution 
is a shrinking solution, 
which is parameterized by 
the space of $n\times n$ skew-Hermitian matrices.

\subsection{Solution with 
Shrinking Phenomenon}
\label{SubSect:SolAsy}

The purpose of this section is 
to provide weaker conditions
\eqref{Cond:BounddkPk},
\eqref{Cond:BoundPk} rather than
\eqref{lkdP}, \eqref{lkP}
under the shrinking condition \eqref{shrink}.
It can be used to determine 
the elements of the solution space
$\mathfrak{s}_k(\bm{z}_{k-1})$. 

\begin{lem}
\label{Lem:AnaControl}
Suppose that there exist real constants $N>0$, $\theta_2>\theta_1$, such that 
$|z| > N$ and 
$\arg z\in [\theta_1,\theta_2]$. We have
\begin{align*}
\| g (z) \| +
\| z^{-\tmop{ad} A} g (z) \|
& \leqslant
C_\varepsilon |z|^{-\varepsilon},
\quad \forall \varepsilon<\varepsilon_0, 
\end{align*}
then for sufficiently small $h$ we have
\begin{align*}
\left\|
z^{\tmop{ad} A} \int_{\infty}^z z^{-\tmop{ad} A}
g(z) \frac{\mathd z}{z}
\right\|
& \leqslant
C_A \frac{ C_{\varepsilon+h} }
{\varepsilon h^{4(\nu-1)} } |z|^{-\varepsilon},
\quad \forall \varepsilon<\varepsilon_0, 
\end{align*}
where $C_A$ depends on $A,\theta_1, \theta_2$, 
and $\nu$ is the nilpotency of
the nilpotent part of $A$.
\end{lem}

\begin{prf}
For $|z| > N \geqslant \mathe^{-\theta}$ and 
$\arg z\in [\theta_1,\theta_2]
\subseteq[-\theta,\theta]$ we have
\begin{subequations}
\begin{align}
\label{Norm:z^a}
|z^\alpha| &\leqslant
\mathe^{ |\tmop{Im}\alpha| \theta } \cdot
|z|^{ \tmop{Re}\alpha },
\quad \forall \alpha \in \mathbb{C}, \\
\label{Norm:lnz^m}
\frac{|\ln z|^m}{m!}  &\leqslant
r^{-m} \mathe^{r \theta} \cdot
|z|^r,
\quad \forall r>0,m\geqslant 0.
\end{align}
\end{subequations}
Consider the Jordan-Chevalley decomposition 
$\mathcal{S}+\mathcal{N}$ of $A$, where
\begin{itemize}
\item $\mathcal{N}$ is the nilpotent part of $A$
with the nilpotency $\nu$;

\item $\mathcal{S}$ is the semisimple part of $A$
with the spectral decomposition
$\mathcal{S}=
\lambda_1 P_1 +\cdots \lambda_{n}P_{n}$.
\end{itemize}
From \eqref{Norm:lnz^m} we have
the operator norm
\begin{align}
\label{Norm:z^adN}
\|z^{\pm \ad \mathcal{N} }\| \leqslant
2\nu \left(
2 \|\mathcal{N}\|
\right)^{2(\nu-1)}
\mathe^{2\|\mathcal{N}\| \theta} \cdot
r^{-2(\nu-1)} |z|^r,\quad \forall r\in(0,2 \|\mathcal{N}\|).
\end{align}
Note that
\begin{align*}
P_i f(z) P_j &=
z^{\lambda_j - \lambda_i - \ad \mathcal{N} }(
P_i ( z^{\ad A} f(z) ) 
P_j),\\
z^{\tmop{ad} A} 
\int_{\infty}^z f(t) \frac{\mathd t}{t} &=
\sum_{1\leqslant i,j \leqslant n}
z^{\lambda_i - \lambda_j + \ad \mathcal{N}}
\int_{\infty}^z P_i f(t) P_j
\frac{\mathd t}{t}.
\end{align*}
Thus, from \eqref{Norm:z^a} and \eqref{Norm:z^adN},
for sufficiently small $r$
we have
\begin{align*}
\| P_i f (z) P_j \|
& \leqslant  
\left\{
\begin{array}{ll}
\|P_i\| \cdot \|P_j\| \cdot
\frac{ C_A' C_\varepsilon }{r^{2(\nu-1)}} 
|z|^{-(\varepsilon-r) - 
\tmop{Re} (\lambda_i - \lambda_j)}& ; 
\tmop{Re} (\lambda_i - \lambda_j) \geqslant
      0\\
\|P_i\| \cdot \|P_j\| \cdot
C_\varepsilon |z|^{-\varepsilon}& ;
\tmop{Re} (\lambda_i - \lambda_j) \leqslant 0
\end{array}\right.,\\
\left\| z^{\tmop{ad} A} 
\int_{\infty}^z f(t) \frac{\mathd t}{t}
\right\| & \leqslant
C_A \cdot
\frac{C_\varepsilon}{ \varepsilon r^{4(\nu-1)} }
|z|^{-\varepsilon + 2r }.
\end{align*}
And we finish the proof.
\end{prf}

In the following, suppose that $f(z)$ is a matrix-valued function, if there exists a constant $C$ such that
  \begin{align*}
    \| f (z) \| \leqslant C \cdot | z |^{\alpha}, \quad
    \alpha \in \mathbb{R},
  \end{align*}
for some $| z | > N$ and
$\arg z\in [\theta_1,\theta_2]$,
then we denote
\begin{align*}
f = O (z^{\alpha}), \quad z \rightarrow \infty .
\end{align*}

\begin{pro}
\label{Lem:ControlisEnough}
Suppose that $\Psi_k(z_k)$ is the solution of 
the $(k, n)$-equation \eqref{kzisoode}
with respect to $z_k$,
and satisfies the shrinking condition
\begin{align}
\underset{z_k \rightarrow \infty}{\overline{\lim}} 
\tmop{max}
| \tmop{Re}
\mathcal{D}_{k - 1} (\Psi_k) | 
& = 
1 - \varepsilon_0 < 1
\label{Cond:BoundEigenLim}.
\end{align}
  \begin{itemize}
    
    \item[(1).] If we assume that
\begin{align}
((\delta_k - \delta_{k - 1}) \Psi_k)^2 
& = 
O (z_k^{1 - \varepsilon}),
\quad z_k \rightarrow \infty, \quad \forall \varepsilon < \varepsilon_0,
\label{Cond:BounddkPk}
    \end{align}
    then $\mathcal{l}_k (\delta_{k - 1} \Psi_k)$ 
    exists and satisfy that
    \begin{align*}
      \delta_{k - 1} \Psi_k & = \mathcal{l}_k (\delta_{k - 1} \Psi_k) + O
      (z_k^{- \varepsilon}), \quad z_k \rightarrow \infty, \quad \forall
      \varepsilon < \varepsilon_0 ; 
    \end{align*}
    
    \item[(2).] If we assume that
    $\mathcal{l}_k (\delta_{k - 1} \Psi_k)$ exists,
    then $\Psi_k\in \mathfrak{s}_k(\bm{z}_{k-1})$
    if and only if
    \begin{align}
      z_k^{\tmop{ad} \mathcal{l}_k (\delta_{k - 1} \Psi_k)} \Psi_k & = 
      O(1), \quad z_k \rightarrow \infty;
      \label{Cond:BoundPk}
    \end{align}

    \item[(3).] If we assume that 
    $\Psi_k\in \mathfrak{s}_k(\bm{z}_{k-1})$,
    then \eqref{Cond:BounddkPk}, \eqref{Cond:BoundPk} hold
    and
    \begin{align*}
      z_k^{\tmop{ad} \delta_{k - 1} \mathcal{l}_k (\Psi_k)} \Psi_k & =
      \mathcal{l}_k (\Psi_k) + O (z_k^{- \varepsilon}), \quad z_k \rightarrow
      \infty, \quad \forall \varepsilon < \varepsilon_0 . 
    \end{align*}
  \end{itemize}
\end{pro}

\begin{prf}
From \eqref{TkExp}, \eqref{z1zk-1Dkk:Expression}
and \eqref{Eq:dk-1pk}, we have
\begin{align*}
(\tmop{ad}_u^{- 1} 
 \tmop{ad}_{\partial u / \partial z_k} - z_k^{- 1})
\delta_k 
& = 
- z_k^{- 1} \delta_{k - 1} 
+ z_k^{- 2} \tmop{ad}_{U_{(k)} /
(\tilde{I}_{k - 1} + z_k^{- 1} U_{(k)})} 
\tmop{ad}_{E_k} \delta_k. 
\end{align*}
Therefore we can verify that
\begin{align*}
    \frac{\partial}{\partial z_k} \delta_{k - 1} \Psi_k 
    & = 
    \delta_{k - 1}
    [\tmop{ad}_u^{- 1} \tmop{ad}_{\partial u / \partial z_k} (\delta_k -
    \delta_{k - 1}) \Psi_k, \delta_k \Psi_k] \nonumber\\
    & = 
    z_k^{- 2} [\tmop{ad}_{U_{(k)} / (\tilde{I}_{k - 1} + z_k^{- 1}
    U_{(k)})} \tmop{ad}_{E_k} (\delta_k - \delta_{k - 1}) \Psi_k, (\delta_k -
    \delta_{k - 1}) \Psi_k] \nonumber\\
    & = 
    z_k^{- 2} \left[ \frac{U_{(k)}}{\tilde{I}_{k - 1} + z_k^{- 1}
    U_{(k)}}, ((\delta_k - \delta_{k - 1}) \Psi_k)^2 \right]. 
\end{align*}
From \eqref{Cond:BounddkPk} we know that 
  $\mathcal{l}_k (\delta_{k - 1} \Psi_k)$ 
  exists and
\begin{align}
    \delta_{k - 1} \Psi_k 
    & = 
    \mathcal{l}_k (\delta_{k - 1} \Psi_k) +
    \int^{z_k}_{\infty} 
    z_k^{- 1} 
    \left[ \frac{U_{(k)}}
    {\tilde{I}_{k - 1} + z_k^{- 1} U_{(k)}}, 
    ((\delta_k - \delta_{k - 1}) \Psi_k)^2 \right]
    \frac{\mathd z_k}{z_k} \nonumber\\
    & = 
    \mathcal{l}_k (\delta_{k - 1} \Psi_k) + O (z_k^{- \varepsilon}),
    \quad z_k \rightarrow \infty, \quad \forall \varepsilon < \varepsilon_0 . 
\label{Pic:dPk}
\end{align}
  Next we notice that from \eqref{Cond:BoundEigenLim}, \eqref{Cond:BoundPk},
  \begin{align*}
    z_k^{- 1 + \tmop{ad} 
    \mathcal{l}_k (\delta_{k - 1} \Psi_k)} 
    \left[
    \frac{U_{(k)}}{\tilde{I}_{k - 1} + z_k^{- 1} U_{(k)}}, 
    ((\delta_k -
    \delta_{k - 1}) \Psi_k)^2 \right] 
    & = 
    O (z_k^{- \varepsilon}), \quad z_k
    \rightarrow \infty. 
  \end{align*}
According to Lemma \ref{Lem:AnaControl},
  \begin{align*}
    z_k^{\tmop{ad} \mathcal{l}_k (\delta_{k - 1} \Psi_k)} 
    \delta_{k - 1} \Psi_k \nonumber
    & = 
    \mathcal{l}_k (\delta_{k - 1} \Psi_k) + z_k^{\tmop{ad}
    \mathcal{l}_k (\delta_{k - 1} \Psi_k)} 
    \int^{z_k}_{\infty} z_k^{- 1}
    \left[ \frac{U_{(k)}}{\tilde{I}_{k - 1} + z_k^{- 1} U_{(k)}}, 
    ((\delta_k - \delta_{k - 1}) \Psi_k)^2 \right] 
    \frac{\mathd z_k}{z_k} \nonumber\\
    & = 
    \mathcal{l}_k (\delta_{k - 1} \Psi_k) + 
    O (z_k^{- \varepsilon}),
    \quad z_k \rightarrow \infty . 
  \end{align*}
  From \eqref{Cond:BoundEigenLim}, \eqref{Cond:BoundPk} we also have
  \begin{align*}
    z_k^{\tmop{ad} \mathcal{l}_k (\delta_{k - 1} \Psi_k)} 
    \tmop{ad}_{U_{(k)}/(\tilde{I}_{k - 1} + z_k^{- 1} U_{(k)})} 
    \tmop{ad}_{E_k} \Psi_k 
    & = 
    O(z_k^{1 - \varepsilon}), \quad z_k \rightarrow \infty. 
  \end{align*}
Therefore
  \begin{align}
    \frac{\partial}{\partial z_k} z_k^{\tmop{ad} \mathcal{l}_k (\delta_{k - 1}
    \Psi_k)} \Psi_k 
    & = 
    - z_k^{- 1 + \tmop{ad} \mathcal{l}_k (\delta_{k - 1}
    \Psi_k)} [\delta_{k - 1} \Psi_k -\mathcal{l}_k (\delta_{k - 1} \Psi_k),
    \Psi_k] \nonumber\\
    & \quad
    + z_k^{- 2 + \tmop{ad} \mathcal{l}_k (\delta_{k - 1} \Psi_k)}
    [\tmop{ad}_{U_{(k)} / (\tilde{I}_{k - 1} + z_k^{- 1} U_{(k)})}
    \tmop{ad}_{E_k} \delta_k \Psi_k, \Psi_k] \nonumber\\
    & = 
    O (z_k^{- 1 - \varepsilon}), \quad z_k \rightarrow \infty. 
\label{Pic:Pk}
  \end{align}
Thus we finish the proof.
\end{prf}

\subsection{Effective Criterion from
Asymptotic Behavior}
\label{SubSect:EffCri}

In this section, 
we will prove that $\mathcal{l}_k$
is injective on $\mathfrak{s}_k(\bm{z}_{k-1})$
through Picard iteration.
Thus the set 
$\mathcal{l}_k(\mathfrak{s}_k) = \mathfrak{b}_k$, 
and $\mathcal{l}_k$ 
are also injective on $\mathfrak{s}_k$. 
Therefore, to determine whether $\Phi_k \in \mathfrak{c}_k$, 
it suffices to verify conditions
\eqref{Cond:BoundEigenLim},
\eqref{Cond:BounddkPk}, 
and \eqref{Cond:BoundPk} for $\Phi_k$
and then inductively verifies whether
$\mathcal{l}_k( \Phi_k ) \in \mathfrak{c}_{k-1}$.

\begin{pro}
\label{Cor:GSolutionAna}
Suppose that 
$\Psi_k$ is the solution of 
the $(k, n)$-equation \eqref{kziso}
that satisfies 
\eqref{Cond:BoundEigenLim},
\eqref{Cond:BounddkPk}, and 
\eqref{Cond:BoundPk}. 
Denote $\Psi_{k-1} =\mathcal{l}_k (\Psi_k)$, then $\Psi_k=\mathcal{s}_k(\Psi_{k-1})
\in\mathfrak{s}_k$, and
$\Psi_{k-1}$ is the solution of 
the $(k-1, n)$-equation.
\end{pro}

\begin{prf}
For $\Psi_{k-1} = 
\mathcal{l}_k (\Psi_k)$
with fixed $\bm{z}_{k-1} =
(z_1,\ldots,z_{k-1})$, 
it can be derived
from \eqref{Pic:dPk} and \eqref{Pic:Pk} that
$z_k^{\ad \delta_{k-1} \Psi_{k-1}} \Psi_k$
satisfies 
the following integral equation
\begin{subequations}
\begin{align}
\label{Eq:IntIso}
X(z_k) & = \Psi_{k-1} +
\int_{\infty}^{z_k} R(z_k, X(z_k) )
\frac{\mathd z_k}{ z_k },\\
R(z_k, X ) & =
-[ \mathcal{T}_{\delta_{k-1}\Psi_{k-1}}
[U(z_k),((\delta_k-\delta_{k-1})X)^2],X] +
[ \ad_{U(z_k)} \ad_{E_k} \delta_k X, X],
\end{align}
\end{subequations}
where we denote
\begin{align*}
(\mathcal{T}_A f)(z) & \assign
z^{\ad A} \int_\infty^z z^{-\ad A}
f(z) \frac{\mathd z}{ z },\\
U(z_k) & = 
z_k^{-1 + \ad \delta_{k-1} \Psi_{k-1}}
\frac{U_{(k)}}
{\tilde{I}_{k - 1} + z_k^{- 1} U_{(k)}}.
\end{align*}
For any $|z_k|>N$ and $\arg z_k \in 
[\theta_1,\theta_2]$,
from \eqref{Cond:BoundEigenLim}
we can verify that
for sufficiently small $h>0$ we have
\begin{align*}
\| U(z_k) \| & \leqslant C'_h |z_k|^{-h}.
\end{align*}
Now if both
$X^{(1)}$, $X^{(2)}$ satisfy 
the integral equation \eqref{Eq:IntIso} and 
$\|X^{(1)}\|_\infty, \|X^{(2)}\|_\infty \leqslant C$,
then according to Lemma \ref{Lem:AnaControl}, 
we have
\begin{align*}
\| X^{(1)} - X^{(2)} \|_\infty \leqslant
C_h N^{-h}
\| X^{(1)} - X^{(2)} \|_\infty,
\end{align*}
where $C_h$ is independent of $N$.
Therefore $\Psi_k$ is unique.
Combining with Proposition \ref{Lem:ControlisEnough},
we finish the proof.
\end{prf}

As an important case, 
we prove that 
the all the skew-Hermitian valued solutions are 
shrinking solutions.

\begin{lem}
\label{Lem:Skew}
If the $(k,n)$-equation \eqref{kziso}
has a solution $\Psi_k(z_1,\ldots,z_k)$,
with skew-Hermitian initial value at
\begin{eqnarray}
z^{(0)}=(z_1^{(0)},\ldots,z_k^{(0)})\in
(\mathbb{R}\setminus\{0\})^k,
\end{eqnarray}
then for all 
$z=(z_1,\ldots,z_k)\in(\mathbb{R}\setminus\{0\})^k$
on the same connected component as $z^{(0)}$,
the solution $\Psi_k(z)$ is skew-Hermitian and 
the Frobenius norm $\|\Psi_k(z)\|_F$ is constant.
\end{lem}

\begin{prf}
It is direct to see that
\begin{eqnarray}
-\overline{\Psi}_k(z)^{\intercal}
\assign
-\overline{\Psi_k(\overline{z})}^{\intercal}
\end{eqnarray}
is also a solution to 
the $(k,n)$-equation \eqref{kziso}.
According to the Picard–Lindelöf theorem, 
we have 
\begin{eqnarray}
-\overline{\Psi}_k(z)^{\intercal}
=\Psi_k(z).
\end{eqnarray}
Since the eigenvalues of $\Psi_k(z)$ are constant, 
the Frobenius norm
\begin{eqnarray}
\|\Psi_k(z)\|_F \assign
\tmop{tr}(
\Psi_k(z)
\cdot
\overline{\Psi_k(z)}^{\intercal})
=
-\tmop{tr}(
\Psi_k(z)
\cdot
\Psi_k(\overline{z}))
\end{eqnarray}
remains constant when $z\in(\mathbb{R}\setminus\{0\})^k$
on the same connected component as $z^{(0)}$.
\end{prf}

Although in general on $\mathbb{C}^k$, 
the solution $\Psi_k$ in Lemma \ref{Lem:Skew}
may not necessarily be skew-Hermitian, 
for convenience, 
we still simply refer to 
the \textbf{skew-Hermitian solutions} 
as those solutions that 
have skew-Hermitian initial values 
at some point in
$(\mathbb{R}\setminus\{0\})^k$.
\begin{cor}
The skew-Hermitian solutions of
the $n$-th isomonodromy equation \eqref{zisoeq}
are all shrinking solutions.
\end{cor}
\begin{prf}
Consider the skew-Hermitian solution $\Psi_k$ 
of the $(k,n)$-equation for $k\leqslant n$ inductively.
Since the eigenvalues of the 
upper left $(k-1)\times (k-1)$ submatrix of $\Psi_k$ 
are pure imaginary, 
thus the 
contraction condition \eqref{Cond:BoundEigenLim} holds
on $(\mathbb{R}\setminus\{0\})^k$. 
According to Lemma \ref{Lem:Skew}, 
the control condition \eqref{Cond:BounddkPk} 
also holds on $(\mathbb{R}\setminus\{0\})^k$, 
therefore
$\mathcal{l}_k(\delta_{k-1}\Psi_{k})$
exists and is also skew-Hermitian, 
thus the control condition \eqref{Cond:BoundPk} holds
on $(\mathbb{R}\setminus\{0\})^k$.

Finally, from Proposition \ref{Cor:GSolutionAna}, 
there is a skew-Hermitian solution
$\Psi_{k-1}=\mathcal{l}_k(\Psi_k)$
of the $(k-1,n)$-equation, 
and the proof can be completed by induction and
Theorem \ref{Thm:AsySol}.
\end{prf}

  \section{Explicit Expressions of Monodromy Data}
\label{Sect:Cat}

In Section \ref{subSect:BasicODE}, 
we will introduce the canonical fundamental solution, 
the monodromy factor $\nu_{d}(u,A)$,
the Stokes matrices $S_d^\pm(u,A)$ and
the central connection matrices $C_d(u,A)$ of 
the linear system \eqref{Conflu}.
The diagonal elements of $u$ can have multiplicity.
In Section \ref{subSect:Fac}, given a shrinking solution $\Phi(u)$ with the boundary value $\Phi_0$, 
we will give the explicit expression of 
the monodromy factor $\nu_{d}(u,\Phi(u))$
in Theorem \ref{Thm:Cat},
as our main result.


\subsection{Confluent Hypergeometric Systems}
\label{subSect:BasicODE}

This subsection will organize 
the basic notations and properties 
of the confluent hypergeometric that people already know.
Most of the results and their proofs can be found in 
\cite{Balser}.

The \textbf{confluent hypergeometric system} 
is the equation for an $n\times n$ matrix-valued function $F(\xi)$
\begin{eqnarray}
\label{Conflu}
\frac{\mathd}{\mathd \xi} F(\xi) =
\left(u+\frac{A}{\xi}
\right)\cdot F(\xi),
\end{eqnarray}
where $A$ is an $n\times n$ constant matrix and $u=\diag(u_1 I_{n_1},\ldots, u_\nu I_{n_\nu})$. Here $u_1,\ldots,u_\nu$ are distinct with multiplicity $n_1,\ldots,n_\nu$. Accordingly, the equation \eqref{Conflu} is divided into
    $(n_1,\ldots,n_\nu)$-blocks,
    and we denote $A=(A^{[i j]})_{\nu\times\nu}$
   the block matrix.

\begin{defi}
Let $\sigma(A)$ denote the set of
eigenvalues of matrix $A$, if
\begin{align*}
\left(
\sigma(A)-\sigma(A)
\right) \cap \mathbb{Z}
= \{0\},
\end{align*}
i.e. the difference between two eigenvalues of $A$
does not take non-zero integers,
we called $A$ \textbf{non-resonant}.
Otherwise, $A$ is called \textbf{resonant}.
The set of $n\times n$ non-resonant matrix
is denoted as $\mathfrak{gl}_n^{\tmop{nr}}(\mathbb{C})$.
\end{defi}

\subsubsection{Stokes matrices}

\begin{pro}
If $A^{[11]},\ldots,A^{[\nu\nu]}$
are non-resonant, then the equation \eqref{Conflu} has 
a unique formal fundamental solution of the following form
\begin{subequations}
\begin{align*}
F^{[\tmop{irr}]}(\xi) 
& =
H^{[\tmop{irr}]}(\xi)  \cdot 
\mathe^{u \xi} \xi^{
\delta_u A
}\\
H^{[\tmop{irr}]}(\xi)
& =
I+\sum_{p=1}^{\infty} H_p^{[\tmop{irr}]} \xi^{-p},
\end{align*}
\end{subequations}
which is divergent in general,
and this is called the
\textbf{canonical formal fundamental solution of 
system \eqref{Conflu} at $\xi=\infty$},
where 
\[\delta_u A:=
\diag(A^{[11]},\ldots,A^{[\nu\nu]})\] is called the 
\textbf{formal monodromy matrix}.
\end{pro}

Although the canonical formal fundamental solution 
$F^{[\tmop{irr}]}(\xi)$
at $\xi=\infty$ is divergent, 
it is known that under Borel resummation method, 
the formal solution $F^{[\tmop{irr}]}(\xi) $ 
corresponds to analytic solutions 
in certain sectorial regions. 
Let us introduce the \textbf{anti-Stokes directions},
which are the directions on the logarithmic Riemann surface 
with the arguments
\begin{align*}
\tmop{aS}(u)\assign
\{
{-\arg (u_i-u_j)
+2k\pi:k\in\mathbb{Z},i\neq j}
\}.
\end{align*}
We arrange these directions/arguments into 
a strictly monotonically increasing sequence
\begin{align*}
\cdots< \tau_{-1} < \tau_{0} < \tau_1 < \cdots.
\end{align*} 
For $d\in(\tau_j,\tau_{j+1})$, 
the \textbf{Stokes sector} $\tmop{Sect}_d$ is denoted as
the sectorial region
\begin{align*}
\tmop{Sect}_d\assign
\{\xi\in\tilde{\mathbb{C}}:
\arg \xi\in \left(\tau_j-\frac{\pi}{2},\tau_{j+1}+\frac{\pi}{2}\right)
\}.
\end{align*}
The following proposition is well known. See e.g., \cite{Balser}.
\begin{pro}
On each $\tmop{Sect}_d$, 
there is a unique
analytic solution $F_d(\xi)$
with the prescribed asymptotics 
$F^{[\tmop{reg}]}(\xi)$.
The solution is called the
\textbf{canonical fundamental solution of 
system \eqref{Conflu} at $\xi=\infty$
with respect to direction $d$}.
\end{pro}
These canonical fundamental solutions at $\xi=\infty$
are in general different 
(that reflects the Stokes phenomenon), 
and the transition between them 
can be measured by Stokes matrices 
$S_d^\pm(u,A)$.

\begin{defi}
\label{Def:MonoData}
For every $d\notin \tmop{aS}(u)$,
we define the \textbf{Stokes matrices} as
\begin{align*}
S_d^{\pm}(u,A) \assign 
F_{d\pm\pi}(\xi)^{-1} F_d(\xi),
\end{align*}
and we denote the \textbf{monodromy factor}
of $F_d(\xi)$ as
\begin{align*}
\nu_{d}(u,A)\assign
F_d(\xi)^{-1} F_d(\xi\mathe^{2\pi\mathi}).
\end{align*}
\end{defi}

The following relations on Stokes matrices are known.
\begin{pro}
\label{AppPro:FormSC}
If $A^{[11]}, \ldots, A^{[\nu \nu]}$ are non-resonant, then we have
\begin{subequations}
\begin{align}
  \label{LUDecomp}
    \nu_{d} (u,A) 
    & = 
    S^-_d (u, A)^{- 1} \cdot \mathe^{2 \pi \mathi
    \delta A} \cdot S^+_d (u, A), \\
    S^{\pm}_{d + 2 k \pi} (u, A) 
    & = 
    \mathe^{- 2 k \pi \mathi \delta
    A} \cdot S^{\pm}_d (u, A) \cdot \mathe^{2 k \pi \mathi \delta
    A}, \quad k \in \mathbb{Z}, \\
    S^{\pm}_{d \mp \pi} (u, A) 
    & = 
    S^{\mp}_d (u, A)^{- 1},
\end{align}
\end{subequations}
and for $c\in\mathbb{C}\setminus\{0\}$ we have
\begin{subequations}
\begin{align}
S_d^\pm(c u, A) & =
c^{\delta A}
S_{d+\arg c}^\pm(u, A)
c^{- \delta A},\\
S_d^\pm(\overline{u}, \overline{A})  =
\overline{S_{-d}^\mp(u, A)},\quad
S_d^\pm(-u^\top, -A^\top) & =
S_{d}^\pm(u, A)^{-\top},\quad
S_d^\pm(-u^\dagger, -A^\dagger)  =
S_{-d}^\mp(u, A)^{-\dagger}.
\end{align}
\end{subequations}
\end{pro}

\begin{pro}
\label{AppPro:StoTri}
  Suppose that $d \nin \tmop{aS} (u)$.
  \begin{itemize}
    \item[(a).] If we have
    \begin{equation}
    \label{Cond:+UpperTri}
      \tmop{Im} (u_1 \mathe^{\mathi d}) > \cdots > \tmop{Im} (u_{\nu}
      \mathe^{\mathi d}),
    \end{equation}
    then the Stokes matrices $S_d^+, S_d^-$ are block upper triangular matrix
    and lower triangular matrix, respectively, with the diagonal block
    identity;
    
    \item[(b).] If we have
    \begin{equation}
    \label{Cond:+LowerTri}
      \tmop{Im} (u_1 \mathe^{\mathi d}) < \cdots < \tmop{Im} (u_{\nu}
      \mathe^{\mathi d}),
    \end{equation}
    then the Stokes matrices $S_d^-, S_d^+$ are block upper triangular matrix
    and lower triangular matrix, respectively, with the diagonal block
    identity.
  \end{itemize}
\end{pro}
\noindent
It can be seen that 
when condition \eqref{Cond:+UpperTri} holds, 
\eqref{LUDecomp} gives 
the LU decomposition of $\nu_{d}(u,A)$; 
when condition \eqref{Cond:+LowerTri} holds, 
\eqref{LUDecomp} gives
the LU decomposition of $\nu_{d}(u,A)^{-1}$. 

Let us introduce the explicit expressions of 
the Stokes matrices $S^\pm_d(u,A)$ 
for confluent hypergeometric systems \eqref{Conflu}
with coefficients $(u,A)=(E_k,\delta_kA)$.
We will use the notations in Notation \ref{Nota:Mat}, 
for example
\begin{align*}
  S_d^{\pm} (E_k, \delta_k A) 
   =  
  \left(\begin{array}{cc}
    S_d^{\pm} (E_k, A^{[k]}) & 0\\
    0 & {I}_{n - k}
  \end{array}\right),
\end{align*}
The following propositions are known. See \cite{Balser88} or \cite{LX}.

\begin{pro}
\label{Pro:S}
Suppose that $A^{[k - 1]}$ is non-resonant and 
  ${q} \in \mathbb{Z}$. For
  $d\in((2{q} - 1) \pi , 2{q} \pi)$ we have
  \begin{align*}
    S_d^+ (E_k, A^{[k]}) 
     = 
    \left(\begin{array}{cc}
      {I} & \mathe^{2{q} \pi \mathi (a_{k k} {I} - A^{[k - 1]})} s_{\hat{k} k}\\
      0 & 1
    \end{array}\right), \quad
    S_d^- (E_k, A^{[k]}) 
     = 
    \left(\begin{array}{cc}
      {I} & 0\\
      - s_{k \hat{k}} \mathe^{(2{q} - 1) \pi \mathi (A^{[k - 1]} - a_{k k}
      {I})} & 1
    \end{array}\right), 
  \end{align*}
  for $d\in(2{q} \pi , (2{q} + 1) \pi)$ we have
  \begin{align*}
    S_d^+ (E_k, A^{[k]}) 
     = 
    \left(\begin{array}{cc}
      {I} & 0\\
      s_{k \hat{k}} \mathe^{(2{q} + 1) \pi \mathi (A^{[k - 1]} - a_{k k}{I})}
      & 1
    \end{array}\right), \quad
    S_d^- (E_k, A^{[k]}) 
     = 
    \left(\begin{array}{cc}
      {I} & - \mathe^{2{q} \pi \mathi (a_{k k} {I} - A^{[k - 1]})} s_{\hat{k}
      k}\\
      0 & 1
    \end{array}\right), 
  \end{align*}
  where
  \begin{align*}
    s_{\hat{k} k} 
     = 
    2 \pi \mathi \cdot \left( \prod_{p, j} \frac{(A^{[k
    - 1]} - \lambda^{(k - 1)}_p {I}) !}{(A^{[k - 1]} - \lambda^{(k)}_j {I}) !}
    \right) A^{[k]}_{\hat{k} k}, \quad
    s_{k \hat{k}} 
     = 
    2 \pi \mathi \cdot A^{[k]}_{k \hat{k}} \left(
    \prod_{p, j} \frac{(\lambda^{(k - 1)}_p {I} - A^{[k - 1]})
    !}{(\lambda^{(k)}_j {I} - A^{[k - 1]}) !} \right).
  \end{align*}
We denote $x!\assign \Gamma(x+1)$.
\end{pro}


\subsubsection{Connection matrices}

\begin{pro}
If $A$ is non-resonant,
then the system \eqref{Conflu} has 
a unique (formal) fundamental solution of the following form 
\begin{align*}
F^{[\tmop{reg}]}(\xi)
& =
H^{[\tmop{reg}]}(\xi) \cdot \xi^{A}\\
H^{[\tmop{reg}]}(\xi) 
& =
I+\sum_{p=1}^{\infty} H_p^{[\tmop{reg}]} \xi^{p},
\end{align*}
which is convergent,
and this is called the
\textbf{canonical fundamental solution of 
system \eqref{Conflu} at $\xi=0$}. 
\end{pro}

\begin{defi}
\label{Def:ConnData}
Suppose that $A$ is non-resonant, then we
define the \textbf{central connection matrices} (with respect to linear system \eqref{Conflu}) as
\begin{align*}
C_d(u,A)\assign
F_d(\xi)^{-1} F^{[\tmop{reg}]}(\xi).
\end{align*}
\end{defi}

The following relations on connection matrices are known.

\begin{pro}
If $A$ and 
$A^{[11]}, \ldots, A^{[\nu \nu]}$ are non-resonant, then we have
\begin{subequations}
\begin{align}
\label{Md=ConjPhi}
\nu_d (u, A) 
& = 
C_d (u, A) \cdot 
\mathe^{2\pi\mathi A} 
\cdot C_d (u, A)^{- 1}, \\
    C_{d + 2 k \pi} (u, A) 
    & = 
    \mathe^{- 2 k \pi \mathi \delta
    A} \cdot C_d (u, A) \cdot \mathe^{2 k \pi \mathi A}, \quad k
    \in \mathbb{Z}, \\
    C_{d \pm \pi} (u, A) 
    & = 
    S^{\pm}_d (u, A) \cdot C_d (u,
    A),
\end{align}
\end{subequations}
and for $c\in\mathbb{C}\setminus\{0\}$ we have
\begin{subequations}
\begin{align}
\label{PropC:cu}
C_d(c u, A) & =
c^{\delta A}
C_{d+\arg c}(u, A)
c^{- A},\\
C_d(\overline{u}, \overline{A})  =
\overline{C_{-d}(u, A)},\quad
C_d(-u^\top, -A^\top) & =
C_{d}(u, A)^{-\top},\quad
C_d(-u^\dagger, -A^\dagger)  =
C_{-d}(u, A)^{-\dagger}.
\label{PropC:uT}
\end{align}
\end{subequations}
\end{pro}

\subsection{Decoupling of the linear system in long time behaviour}
\label{subSect:Fac}
Let us consider the $n\times n$ linear system \eqref{introisoStokeseq1}-\eqref{introisoStokeseq2}, i.e.,
$F(\xi,u_1,\ldots,u_n)\in {\rm GL}_n(\mathbb{C})$
\begin{equation}
\label{FxiuSys}
\left\{
\begin{alignedat}{3}
\frac{\partial}{\partial \xi}
F(\xi,u) 
& =
\left(u+
\frac{\Phi(u;\Phi_0)}{\xi}
\right)\cdot F(\xi,u)\\
\frac{\partial}{\partial u_k}
F(\xi,u) 
& =
\left(E_k \xi+
\ad_u^{-1}\ad_{E_k}\Phi(u;\Phi_0)
\right)\cdot F(\xi,u)
\end{alignedat}
\right..
\end{equation}
One checks that the compatibility condition of the linear system is the nonlinear equation \eqref{isoeq}. 

For any fixed $u$ and any direction $d\notin \tmop{aS}(u)$, let $F_d(\xi,u)$ be the canonical fundamental solution of the first equation in \eqref{FxiuSys}. The following proposition is known. 

\begin{pro}\cite{JMMS, JMU}
The $n\times n$ matrix function $F_d(\xi,u)$ is a solution of the linear system \eqref{FxiuSys}. As a corollary, the Stokes matrices $S_d^\pm(u,\Phi(u;\Phi_0))$ are locally constant.   
\end{pro}
\begin{rmk}
The monodromy factor $\nu_{d}(u,\Phi(u))=F_d(\xi,u)^{-1}
F_d(\xi\mathe^{2\pi\mathi},u)$
of the canonical fundamental solution $F_d(\xi,u)$ is also locally preserved.
It should be emphasized that the matrix
$\nu_{d}(u, \Phi(u))$ also depends on 
the arguments of $u_1,\ldots,u_n$.
When $u_k$ is varied to 
$u_k\mathe^{2\pi\mathi}$, 
even though the position of $u$ remain unchanged, 
but the solution $\Phi(u)$ may no longer be the same,
thus $\nu_{d}(u,\Phi(u))$ may be different along the continuation.
\end{rmk}
\noindent
In particular, 
given any direction $d$, let $\mathcal{R}_d\in\mathbb{C}^n\setminus \Delta$ be the set consisting of all $u$ such that $d\notin \tmop{aS}(u)$, then the Stokes matrices $S_d^\pm(u,\Phi(u;\Phi_0))$ are constant for all $u$ in each connected component of $\mathcal{R}_d$. 
Let us take the following connected components $\mathcal{R}_d(J)$ of $\mathcal{R}_d$ labelled by the subsets 
$J\subseteq \{1,\ldots,n-1\}$: $u\in \mathcal{R}_d(J)$ if
\begin{align}
\tmop{Im}(u_{k+1} \mathe^{\mathi d}) &<
\min_{1\leqslant j \leqslant k} 
\tmop{Im}(u_j \mathe^{\mathi d}), \
\text{ for } k\in J,\\
\tmop{Im}(u_{k+1} \mathe^{\mathi d}) &>
\max_{1\leqslant j \leqslant k} 
\tmop{Im}(u_j \mathe^{\mathi d}), \ 
\text{ for } k\notin J.
\end{align}
See Appendix \ref{Sect:Diag}
for more details.

The following theorem gives a factorization of the monodromy $\nu_{d}(u,\Phi)=S^-_d (u, \Phi)^{- 1} \mathe^{2 \pi \mathi
    \delta \Phi} S^+_d (u, \Phi)$.
\begin{thm}
\label{Thm:Cat}
Suppose that $\Phi(u;\Phi_0)$ is 
a shrinking solution of 
the $n$-th isomonodromy equation \eqref{isoeq},
then we have
\begin{align}
\label{Cat}
\nu_{d}(u,\Phi(u;\Phi_0)) =
\tmop{Ad}\left(
\overrightarrow{\prod^n_{k=1}}
C_{d+\arg(u_k- u_{k-1} )}(E_k,\delta_k\Phi_0)
\right)
\mathe^{2\pi\mathi
\Phi_0}, \ \text{ for all } u\in \mathcal{R}_d(J).
\end{align}
\end{thm}
The rest of this subsection devotes to the proof of Theorem \ref{Thm:Cat}. The proof is based on Lemma \ref{decouple1}-\ref{Lem:PhikCat}, which state that in the long time $u_1\gg u_2\gg \cdots \gg u_n$ asymptotics, the system can be decoupled into a sequence of (exactly solvable) simpler and lower rank systems. After introducing these lemmas, the proof of Theorem \ref{Thm:Cat} is given in the end of this subsection.

Let us start by supposing that 
$\Phi_n(z_n)\in \mathfrak{s}_{n}(\bm{z}_{n-1})$, i.e., $\Phi_n(z_n)$ is a solution of
the $n$-th isomonodromy equation \eqref{kzisoode}
with respect to $z_n$
and satisfies the shrink condition
    \eqref{Cond:BoundEigenLim} and control conditions
    \eqref{Cond:BounddkPk}, \eqref{Cond:BoundPk}.
From Proposition \ref{Lem:ControlisEnough},
there is an $n\times n$ matrix 
$\Phi_{n-1} \in \mathfrak{b}_{n-1}(\rm{pt})$
such that
\begin{subequations}
\begin{align}
\label{dPnrmvar}
\delta_{n-1}
\Phi_n(z_n) &=
\delta_{n-1}
\Phi_{n-1} +
O(z_n^{-\varepsilon}),
\quad z_n\to\infty,
\\
\label{Pnrmvar}
z_n^{\ad\delta_{n-1}
\Phi_{n-1}}
\Phi_n(z_n) &=
\Phi_{n-1} +
O(z_n^{-\varepsilon}),
\quad z_n\to\infty
\end{align}
\end{subequations}
Fixed $\bm{z}_{n-1} = (z_1,\ldots,z_{n-1})$, 
then the
linear differential systems
\eqref{FxiuSys} reduces to
\begin{equation}
\label{Sysxizn}
\left\{
\begin{alignedat}{3}
\frac{\partial}{\partial \xi}
F(\xi,z_n) 
& =
\left(u+
\frac{\Phi_n(z_n)}{\xi}
\right)\cdot F(\xi,z_n) \\
\frac{\partial}{\partial z_n}
F(\xi,z_n) 
& =
\left(\frac{\partial u}{\partial z_n}\xi+
\ad_u^{-1}
\ad_{\partial u/\partial z_n}\Phi_n(z_n)
\right)\cdot F(\xi,z_n) 
\end{alignedat}
\right..
\end{equation}
Let us introduce the coordinate transformation
$\left\{
\begin{alignedat}{3}
x     & = \xi\\
\zeta_n & = \xi z_n
\end{alignedat}
\right.$
to prove Theorem \ref{Thm:Cat}.
\begin{lem}\label{decouple1}
In terms of the new coordinates $x,\zeta_n$, the equation \eqref{Sysxizn} takes the form
\begin{equation}
\label{xzetaSys}
\left\{
\begin{alignedat}{3}
\frac{\partial}{\partial x} 
F (x, \zeta_n)
& =
A_1(x,\zeta_n)
\cdot F (x, \zeta_n) \\
\frac{\partial}{\partial \zeta_n} 
\left(\left(
\frac{\zeta_n}{x}\right)^
{\delta_{n-1}\Phi_{n-1}}
F(x,\zeta_n)\right)
& =
A_2(x,\zeta_n)
\cdot 
\left(
\frac{\zeta_n}{x}\right)
^{\delta_{n - 1}\Phi_{n - 1}} 
F (x, \zeta_n)
\end{alignedat}
\right.,
\end{equation}
where the coefficient matrices satisfy as $x\to 0,\zeta_n\to\infty$
\begin{align}
A_1(x,\zeta_n) &=
u^{(n)} +
\frac{\delta_{n - 1} \Phi_{n - 1}}{x} + 
O (x^{-1+\varepsilon} \zeta_n^{- \varepsilon})\\
A_2(x,\zeta_n) &=
\frac{\partial u_n}{\partial z_n} E_n +
\frac{\Phi_{n - 1}}{\zeta_n} + 
O (x^{\varepsilon} \zeta_n^{-1-\varepsilon}),
\end{align}
with $u^{(n)}:=\diag(
u_1,\ldots,u_{n-1},u_{n-1}
)$.
\end{lem}

\begin{prf}
From \eqref{TkExp}, \eqref{z1zk-1Dkk:Expression}
and \eqref{Eq:dk-1pk}, we have
\begin{align*}
\tmop{ad}_u^{- 1} 
\tmop{ad}_{\partial u / \partial z_n} 
\Phi_n (z_n) & =
z_n^{- 1} \delta_{n - 1}^{\bot} 
\Phi_n (z_n) + 
z_n^{-2}
\ad_{U_{(n)}/(\Tilde{I}_{n-1}+z_n^{-1}U_{(n)})}
\ad_{E_n}\Phi_n(z_n).
\end{align*}
Therefore
from \eqref{dPnrmvar} and \eqref{Pnrmvar}, we get
\begin{align}
\nonumber
A_1(x,\zeta_n) &=
\frac{\partial \xi}{\partial x} 
\left( u + \frac{\Phi_n (z_n)}{\xi}
\right) + 
    \frac{\partial z_n}{\partial x} \left( \frac{\partial
    u}{\partial z_n} \xi + \tmop{ad}_u^{- 1} \tmop{ad}_{\partial u / \partial
    z_n} \Phi_n (z_n) \right)\\
\nonumber
& = 
u^{(n)} + 
\frac{\delta_{n-1}\Phi_n(z_n)}{x}- 
\frac{1}{\zeta_n}
\ad_{U_{(n)}/(\Tilde{I}_{n-1}+z_n^{-1}U_{(n)})}
\ad_{E_n}\Phi_n(z_n)\\
\label{A1expan}
& = u^{(n)} + \frac{\delta_{n - 1} \Phi_{n - 1}}{x} + 
    O (x^{-1+\varepsilon}
    \zeta_n^{- \varepsilon}), 
\quad 
x \rightarrow 0, \ 
\zeta_n \rightarrow \infty,
\end{align}
  and
\begin{align}
\nonumber
A_2(x,\zeta_n) 
&=
\left(
\frac{\zeta_n}{x}\right)^
{\ad\delta_{n-1}\Phi_{n-1}} 
\frac{\partial
    z_n}{\partial \zeta_n} \left( \frac{\partial u}{\partial z_n} \xi +
    \tmop{ad}_u^{- 1} \tmop{ad}_{\partial u / \partial z_n} \Phi_n (z_n)
    \right) + \frac{\delta_{n - 1} \Phi_{n - 1}}{\zeta_n}\\
\nonumber
& =  
\left(
\frac{\zeta_n}{x}\right)^
{\ad\delta_{n-1}\Phi_{n-1}}
\left(
    \frac{\partial u_n}{\partial z_n} E_n + \frac{1}{x} \tmop{ad}_u^{- 1}
    \tmop{ad}_{\partial u / \partial z_n} \Phi_n (z_n) \right) +
    \frac{\delta_{n - 1} \Phi_{n - 1}}{\zeta_n}\\
\label{A2expan}
& = 
\frac{\partial u_n}{\partial z_n} E_n + 
\frac{\Phi_{n -1}}{\zeta_n} + O (x^{\varepsilon} \zeta_n^{- 1 - \varepsilon}),
\quad 
x \rightarrow 0, \ 
\zeta_n \rightarrow \infty.
\end{align}
Thus we conclude the proof.
\end{prf}

To proceed, let us introduce the following notations 
for some canonical fundamental solutions:
\begin{itemize}
\item 
$F(\xi,z_n) =
H(\xi,z_n)
\xi^{\Phi_n(z_n)}$ as 
the canonical fundamental solution of 
system \eqref{Sysxizn}
for fixed $z_n$
at $\xi = 0$,
and
$F_d(\xi,z_n) =
H_d(\xi,z_n)
\xi^{\delta\Phi}
\mathe^{u \xi}$ as 
the canonical fundamental solution
for fixed $z_n$
at $\xi = \infty$ 
with respect to direction $d$;

\item 
$F^{[0]} (x, \zeta_n) =
H^{[0]} (x, \zeta_n)
x^{\delta_{n-1}\Phi_{n-1}}$ as 
the canonical fundamental solution of 
systems \eqref{xzetaSys} 
for fixed $\zeta_n$
at $x = 0$;

  
\item $T^{[0]} (x)=
Y^{[0]} (x)
x^{\delta_{n-1}\Phi_{n-1}}$
as the canonical fundamental solution of 
the following system at $x = 0$
  \begin{align}
    \frac{\partial}{\partial x} T (x)  =  \left( u^{(n)} + \frac{\delta_{n -
    1} \Phi_{n - 1}}{x} \right) \cdot T (x),
  \end{align}
and 
$T^{[0]}_d(x)=Y^{[0]}_d(x)
x^{\delta\Phi}
\mathe^{u^{(n)}x}$ as 
the canonical fundamental solution 
at $x = \infty$ 
with respect to direction $d$;

\item $T^{[n]} (\zeta_n)=
Y^{[n]} (\zeta_n)
\zeta_n^{\Phi_{n-1}}$
as the canonical fundamental solution of 
the following system at $\zeta_n = 0$
  \begin{align}
  \label{Eq:T2}
    \frac{\partial}{\partial \zeta_n} T (\zeta_n)  =  \left( \frac{\partial
    u_n}{\partial z_n} E_n + \frac{\Phi_{n - 1}}{\zeta_n} \right) \cdot T
    (\zeta_n).
  \end{align}
and
$T^{[n]}_d (\zeta_n) =
Y^{[n]}_d (\zeta_n)
\zeta_n^{\delta_{n-1}\Phi_{n-1}}
\mathe^{(u_{n-1}-u_{n-2})\zeta_n E_n}$ as 
the canonical fundamental solution 
at $\zeta_n = \infty$ 
with respect to the direction $d$.
\end{itemize}

\begin{lem}
\label{Lem:CatPic}
The canonical fundamental solutions have the following locally uniform limits
\begin{itemize}
\item
when $x\to 0$
\begin{align}
\label{Lim:adH1}
\left(
\frac{\zeta_n}{x}
\right)^{\ad\delta_{n-1}\Phi_{n-1}} 
H^{[0]}(x,\zeta_n) 
&\xrightarrow{x\to 0} {I},\\
\label{Lim:addH1}
\frac{\partial }
{\partial \zeta_n}
\left(
\frac{\zeta_n}{x}
\right)^{\ad\delta_{n-1}\Phi_{n-1}} 
H^{[0]}(x,\zeta_n)
&\xrightarrow{x\to 0} 0,\\
\label{Lim:adH}
\left(
\frac{\zeta_n}{x}
\right)^{\ad\delta_{n-1}\Phi_{n-1}} 
H(x,\zeta_n/x)
&\xrightarrow{x\to 0}
Y^{[n]}(\zeta_n);
\end{align}

\item
when $\zeta_n\to\infty$
\begin{align}
\label{Lim:H1}
H^{[0]}(x,\zeta_n)
&\xrightarrow{\zeta_n\to \infty}
Y^{[0]}(x);
\end{align}

\item
when $\zeta_n\to\infty$
with $\arg \zeta_n\in
(d_{\zeta_n}-\varepsilon,
 d_{\zeta_n}+\varepsilon)$,
\begin{align}
\label{Lim:Hd}
H_{d_x}(x,\zeta_n/x)
&\xrightarrow{\zeta_n\to \infty}
Y_{d_x}^{[0]}(x),\\
\label{Lim:adYd}
\zeta_n^{-\ad
\delta_{n-1}\Phi_{n-1}}
Y_{d_{\zeta_n}}^{[n]}(x,\zeta_n)
&\xrightarrow{\zeta_n\to \infty}{I},
\end{align}
where $d_x\assign \arg x$.
\end{itemize}
\end{lem}

\begin{prf}
For convenience,
we will only provide the proofs for 
\eqref{Lim:adH1}, \eqref{Lim:addH1}.
The approach for the rest
are similar, 
see \cite{Xu} for more details.

Denote
$\Tilde{H}^{[0]}(x,\zeta_n)=
\left(
\frac{\zeta_n}{x}
\right)^{\ad\delta_{n-1}\Phi_{n-1}} 
H^{[0]}(x,\zeta_n)$.
Note that 
$\Tilde{H}^{[0]}(x,\zeta_n)$
satisfies the equation
\begin{align*}
\frac{\partial}{\partial x} \Tilde{H}^{[0]}(x,\zeta_n) 
& = \left( 
\left(
\frac{\zeta_n}{x}
\right)^{\ad\delta_{n-1}\Phi_{n-1}}
A_1 (x,\zeta_n) -
\frac{\delta_{n-1}\Phi_{n-1}}{x}
\right) \cdot 
\Tilde{H}^{[0]}(x,\zeta_n).
\end{align*}
For $|x|< h$
with bounded argument,
from \eqref{A1expan}
we can verify that
there exists a constant $C$ 
independent of $h$, 
such that
\begin{align*}
\left\| 
\left(
\frac{\zeta_n}{x}
\right)^{\ad\delta_{n-1}\Phi_{n-1}}
A_1 (x,\zeta_n) -
\frac{\delta_{n-1}\Phi_{n-1}}{x}
\right\| & \leqslant  
C \cdot 
|x|^{-1 +\varepsilon},\\
\left\|
\frac{\partial}
{\partial \zeta_n}
\left(
\frac{\zeta_n}{x}
\right)^{\ad\delta_{n-1}\Phi_{n-1}}
A_1 (x,\zeta_n)\right\| 
& \leqslant  
C \cdot 
|x|^{-1+\varepsilon}.
\end{align*}
Let us consider the Picard iteration of the corresponding equation with $P_0(x,\zeta_n)=I$
  \begin{align*}
    P_m (x, \zeta_n) & = 
    {I} + \int_0^x \left( 
\left(
\frac{\zeta_n}{x}
\right)^{\ad\delta_{n-1}\Phi_{n-1}}
A_1 (x,\zeta_n) -
\frac{\delta_{n-1}\Phi_{n-1}}{x}
    \right) \cdot P_{m - 1} (x, \zeta_n) \mathd x,\\
    \nonumber
    \frac{\partial P_m (x, \zeta_n)}{\partial \zeta_n} 
    & = \int_0^x \left(
\left(
\frac{\zeta_n}{x}
\right)^{\ad\delta_{n-1}\Phi_{n-1}}
A_1 (x,\zeta_n) -
\frac{\delta_{n-1}\Phi_{n-1}}{x}
    \right) \cdot \frac{\partial P_{m -
    1} (x, \zeta_n)}{\partial \zeta_n} \mathd x\\
    &\phantom{=}
    + \int_0^x \left( 
\frac{\partial}
{\partial \zeta_n}
\left(
\frac{\zeta_n}{x}
\right)^{\ad\delta_{n-1}\Phi_{n-1}}
A_1 (x,\zeta_n)
    \right) \cdot P_{m - 1}
    (x, \zeta_n) \mathd x.
  \end{align*}
Let us introduce the following norm
  \begin{align*}
    \| F \|_{\varepsilon,h} 
    & \assign 
    \underset{0 < x < h}{\sup} \| x^{-
    \varepsilon} F (x) \|,
  \end{align*}
then there exists a constant $C_0$ 
  independent of $h$ 
  such that
  \begin{align*}
\| P_{m+1}-P_m \|_{\varepsilon,h}
& \leqslant 
C_0 |h|^{\varepsilon} 
\cdot 
\| P_m-P_{m-1} \|_{\varepsilon,h},\\
\left\| 
\frac{\partial P_{m+1}}{\partial \zeta_n} - 
\frac{\partial P_{m  }}{\partial \zeta_n}
\right\|_{\varepsilon,h}
& \leqslant C_0 |h|^{\varepsilon} 
\cdot
\left(
\| P_{m  }-P_{m-1} \|_{\varepsilon,h}
+
\left\| 
\frac{\partial P_{m  }}{\partial \zeta_n} - 
\frac{\partial P_{m-1}}{\partial \zeta_n} 
\right\|_{\varepsilon,h}
\right).
\end{align*}
Therefore
$P_m$ and
$\frac{\partial P_{m  }}{\partial \zeta_n}$
are locally uniformly convergent to
$\Tilde{H}^{[0]}(x,\zeta_n)$ and 
$\frac{\partial \Tilde{H}^{[0]}(x,\zeta_n)}
{\partial \zeta_n}$, 
and we have
\begin{align*}
\Tilde{H}^{[0]}(x,\zeta_n)
& = {I}+O(x^\varepsilon),
\quad x\to 0,\\
\frac{\partial \Tilde{H}^{[0]}(x,\zeta_n)}
{\partial \zeta_n}
& = O(x^\varepsilon),
\quad x\to 0.
\end{align*}
Thus we prove
\eqref{Lim:adH1} and \eqref{Lim:addH1}.
\end{prf}

\begin{lem}
\label{Lem:PreCat}
Suppose that
\begin{align}
\label{Cond:NonDengStokesS}
\tmop{Im}(u_n \mathe^{\mathi d_x}) <
\min_{1\leqslant k < n} \tmop{Im}(u_k \mathe^{\mathi d_x})
\quad \text{ or } \quad
\tmop{Im}(u_n \mathe^{\mathi d_x}) >
\max_{1\leqslant k < n} \tmop{Im}(u_k \mathe^{\mathi d_x}),
\end{align}
then we have
\begin{align}
\nonumber
F_{d_x}(x,\zeta_n/x) \cdot
C_{d_x} (u,\delta_{n-1}\Phi_{n-1})
&=
F^{[0]}(x,\zeta_n) \cdot 
\zeta_n^{-\delta_{n-1} \Phi_{n-1}}
T^{[n]}_{d_{\zeta_n}} (\zeta_n)\\
\label{PreCat}
&=
F_{d_{\zeta_n}}^{[n]}(x,\zeta_n) \cdot 
x^{-\delta_{n-1} \Phi_{n-1}}
T^{[0]} (x),
\end{align}
where $d_x \assign \arg x$,
$d_{\zeta_n} \assign \arg \zeta_n$.
\end{lem}

\begin{prf}
Denote
$\Tilde{H}^{[0]}(x,\zeta_n)=
\left(
\frac{\zeta_n}{x}
\right)^{\ad \delta_{n - 1}\Phi_{n - 1}} 
H^{[0]}(x,\zeta_n)$.
Suppose that the system \eqref{xzetaSys} 
has a solution of the form
$F^{[0]} (x, \zeta_n) 
\zeta_n^{-\delta_{n - 1}\Phi_{n - 1}}
X(\zeta_n)$, 
then 
\begin{align*}
\frac{\partial}{\partial \zeta_n}  
(\Tilde{H}^{[0]} X)
= \frac
{\partial 
\Tilde{H}^{[0]}}
{\partial \zeta_n} 
X + 
\Tilde{H}^{[0]}
\frac
{\partial X}
{\partial \zeta_n}
= A_2 \cdot
\Tilde{H}^{[0]} X, 
\end{align*}
which can be rewritten as
\begin{align}
\label{T2eqLem}
\frac{\partial X}{\partial \zeta_n} \cdot 
X^{- 1} 
& = 
(\Tilde{H}^{[0]})^{- 1} A_2 
\Tilde{H}^{[0]} - 
(\Tilde{H}^{[0]})^{- 1} 
\frac{\partial \Tilde{H}^{[0]}}
{\partial \zeta_n}.
\end{align}
From \eqref{Lim:adH1} and \eqref{Lim:addH1}
in Lemma \ref{Lem:CatPic},
we can verify that $X(\zeta_n)$
satisfies the system \eqref{Eq:T2}
by taking the limit $x\to 0$
in \eqref{T2eqLem}.
Therefore, the function 
$F^{[0]}(x,\zeta_n) \cdot 
\zeta_n^{-\delta_{n-1} \Phi_{n-1}}
T^{[n]}_{d_{\zeta_n}} (\zeta_n)$ is a solution of
the system \eqref{xzetaSys}. Recall that by definition, $F_{d_x}(x,\zeta_n/x)$ is also a solution of \eqref{xzetaSys}.

Next, we fix $x$ and 
take the limit $\zeta_n\to\infty$. 
Condition \eqref{Cond:NonDengStokesS} 
ensures that
the Stokes sector of $F_{d_x}$ 
containing the ray $\arg x=d_x$ will not degenerate.
From \eqref{Lim:Hd}, \eqref{Lim:H1}
and \eqref{Lim:adYd}
in Lemma \ref{Lem:CatPic},
we locally uniformly have
the following asymptotic expansion
\begin{align*}
F_{d_x}(x,\zeta_n/x) \cdot
C_{d_x} (u,\delta_{n-1}\Phi_{n-1})
&\sim
({I}+o(1;\zeta_n))
\zeta_n^{\delta_{n-1}\Phi_{n-1}}
\mathe^{
\frac{\partial u_n}{\partial z_n}
\zeta_n E_n}\cdot 
x^{-\delta_{n-1} \Phi_{n-1}}
T^{[0]} (x),\\
F^{[0]}(x,\zeta_n) \cdot 
\zeta_n^{-\delta_{n-1} \Phi_{n-1}}
T^{[n]}_{d_{\zeta_n}} (\zeta_n)
&\sim
({I}+o(1;\zeta_n))
\zeta_n^{\delta_{n-1}\Phi_{n-1}}
\mathe^{
\frac{\partial u_n}{\partial z_n}
\zeta_n E_n}\cdot 
x^{-\delta_{n-1} \Phi_{n-1}}
T^{[0]} (x).
\end{align*}
Since the solution of the linear system 
is uniquely determined by 
its asymptotic leading term, we have
\begin{align*}
F_{d_x}(x,\zeta_n/x) \cdot
C_{d_x} (u,\delta_{n-1}\Phi_{n-1})
&=
F_{d_{\zeta_n}}^{[n]}(x,\zeta_n) \cdot 
x^{-\delta_{n-1} \Phi_{n-1}}
T^{[0]} (x),\\
F^{[0]}(x,\zeta_n) \cdot 
\zeta_n^{-\delta_{n-1} \Phi_{n-1}}
T^{[n]}_{d_{\zeta_n}} (\zeta_n)
&=
F_{d_{\zeta_n}}^{[n]}(x,\zeta_n) \cdot 
x^{-\delta_{n-1} \Phi_{n-1}}
T^{[0]} (x).
\end{align*}
We conclude the proof.
\end{prf}

\begin{rmk}
\label{Rmk:PreCat}
Condition \eqref{Cond:NonDengStokesS} can
be rewritten as
\begin{align*}
u_n &\notin
\tmop{conv}\{u_1,\ldots,u_{n-1}\},\\
d_x &\notin
-\arg(u_n -
\tmop{conv}\{u_1,\ldots,u_{n-1}\})+
\mathbb{Z}\pi.
\end{align*}
It can be seen in \eqref{PreCat} that 
both sides of the identity should depend on 
both $\arg \zeta_n = d_{\zeta_n}$ 
and $\arg x = d_{x}$
at the same time, 
but the identity itself does not explicitly show it.
In fact
\begin{itemize}

\item when we fix $x$ and
$\zeta_n$ is large enough,
if $d_{\zeta_n}+
\arg(\frac{\partial u_n}{\partial z_n})=
\arg(\frac{\partial u_n}{\partial z_n}\zeta_n)$ 
crosses $\mathbb{Z}\pi$,
the Stokes sector of $F_{d_x}(x,\zeta_n/x)$
which contains the ray $\arg x=d_x$ 
on the left-hand side will change,
therefore the $F_{d_x}$ 
will not remain the same,
and the right-hand side
will naturally change in these crossing 
by the definition;

\item when we fix $\zeta_n$, 
the right-hand side
does not depend on the choice of $d_x$.
But under the condition \eqref{Cond:NonDengStokesS}, 
both $F_{d_x}$ and $C_{d_x}$ 
in the left-hand side 
will simultaneously change 
with the choice of $d_x$.
\end{itemize}
\end{rmk}


\begin{lem}
\label{Lem:PhikCat}
Suppose that $\Phi_{k-1}\in
\mathfrak{b}_{k-1}(\mathrm{pt})$ and
$\Phi_k=\mathcal{s}_k(\Phi_{k-1})$, then
\begin{align}
\label{CdDeng}
\lim_{z_k\to\infty}
C_{d}(u,\delta_k\Phi_k(z_k))
z_k^{-\delta_{k-1}\Phi_{k-1}}
z_k^{\delta_{k}\Phi_{k-1}}
& =
C_{d}(u,\delta_{k-1}\Phi_{k-1})
C_{d+\arg z_k}\left(
\frac{\partial u_k}{\partial z_k}E_k,
\delta_k\Phi_{k-1}\right),
\end{align}
where the argument of $z_k$ is 
restricted in a sufficiently small range.
\end{lem}

\begin{prf}
Since $\delta_k\Phi_k(z_k)$ 
satisfies the $k$-th isomondromy equation
with respect to $z_k$,
we can assume that $k=n$
without loss of generality.
Denote $z_n = \zeta_n/x$,
from \eqref{Lim:adH} and \eqref{Lim:adH1}
we have
\begin{align*}
  z_n^{\delta_{n - 1} \Phi_{n - 1}} F (x, \zeta_n / x) z_n^{- \delta_{n - 1}
  \Phi_{n - 1}} z_n^{\Phi_n} 
  & \xrightarrow{x \rightarrow 0} 
  T^{[n]}(\zeta_n),\\
  z_n^{\delta_{n - 1} \Phi_{n - 1}} F^{[0]} (x, \zeta_n) \zeta_n^{- \delta_{n
  - 1} \Phi_{n - 1}} T^{[n]} (\zeta_n) 
  & \xrightarrow{x \rightarrow 0}
  T^{[n]}(\zeta_n),
\end{align*}
By applying \eqref{PreCat} in Lemma \ref{Lem:PreCat},
we can complete the proof.
\end{prf}

\begin{prf}[Proof of Theorem \ref{Thm:Cat}]
Suppose that $\Phi_n(z)$ is 
the shrinking solution of 
the $n$-th isomonodromy equation \eqref{zisoeq}
with boundary value $\Phi_0$,
and $z_1,\ldots,z_n$ are large enough.
Denote $\hat{\Phi}_k =
(z_k\cdots z_1)^{\ad\delta_k\Phi_k}
\Phi_k$, 
note that $\delta_k\hat{\Phi}_k =
\delta_k\Phi_k$.
We inductively assume that
\begin{align}
\label{CatuRec}
\nu_{d}(u,\Phi_n(z)) & =
\tmop{Ad}\left(
C_{d}(u,\delta_k\hat{\Phi}_k)
\overrightarrow{\prod^n_{j=k+1}}
C_{d+\arg (u_j-u_{j-1})}\left(
E_j
,\delta_j\hat{\Phi}_k\right)
\right)
\mathe^{2\pi\mathi\hat{\Phi}_k}.
\end{align}
Let us rewrite \eqref{CdDeng} in Lemma \ref{Lem:PhikCat} as
\begin{align*}
\lim_{u_k\to\infty}
C_{d}(u,\delta_k\hat{\Phi}_k)
u_k^{-\delta_{k-1}\hat{\Phi}_{k-1}}
u_k^{\delta_{k}\hat{\Phi}_{k-1}}
& =
C_{d}(u,\delta_{k-1}\hat{\Phi}_{k-1})
C_{d+\arg (u_k-u_{k-1}) }\left(
E_k,
\delta_k\hat{\Phi}_{k-1}\right),
\end{align*}
from which we see that
\eqref{CatuRec} holds for $k+1$. 
Since
the induction starting point $k=n$
holds by definition, 
thus \eqref{CatuRec} inductively holds for $k=0$.
Note that $\hat{\Phi}_0=\Phi_0$, 
we thus complete the proof.
\end{prf}\\

\begin{rmk}
If we want to determine 
the value of shrinking solution $\Phi(u;\Phi_0)$,
beside the $n\times n$ matrix $\Phi_0$
and $u\in\mathbb{C}^n\setminus\Delta$,
we also need to know
the argument $\arg(u_k-u_{k-1})$. 
Therefore, the right side of the equation 
depends not only on $\Phi_0$, 
but also on $\arg(u_k-u_{k-1})$.
\end{rmk}

  \section{Almost Every Solution is a Shrinking Solution}
\label{Sect:GeneralSol}

Suppose that $\Phi_n(z)$ is 
any solution of 
the $n$-th isomonodromy equation \eqref{isoeq} with the monodromy data
\begin{align}
\label{Def:MD}
\tmop{MD}_{u,d} ( \Phi_n ) & \assign
(\nu_{u,d}(\Phi_n(z)),
\sigma(\Phi_n(z)),
\delta\Phi_n(z)),
\end{align} 
where $\nu_{u,d} (\Phi_n(z))\assign
\nu_{d} (u,\Phi_n(z))$.
In Section \ref{subSect:ECfromMD},
we will provide a criterion on the monodromy data, 
to determine whether $\Phi_n(z)$
(with non-resonant initial value)
is a shrinking solution.
In Section \ref{subSect:GP}, 
we will prove that 
the set of shrinking solutions $\Phi_n(z)$,
associated to all boundary values $\Phi_0$ 
contains almost all solutions.
Finally, we will provide a concrete condition 
on the monodromy data, 
such that $\Phi_n(z)$ is a shrinking solution.

\subsection{Effective Criterion from
Monodromy Data}
\label{subSect:ECfromMD}

In this section, 
from the monodromy data $\tmop{MD}_{u,d}(\Phi_n)$ of 
the solution $\Phi_n$,
we will use the explicit Riemann-Hilbert mapping \eqref{Cat}
provided by Theorem \ref{Thm:Cat}, 
to determine whether a solution $\Phi_n$
that satisfies the non-resonant condition is a shrinking solution.
In Section \ref{subSect:NGS}, 
we will see that this criterion will fail 
under the resonant condition. 
For convenience, we introduce the following notations
\begin{nota}
\label{Nota:Mat}
For $n\times n$ matrix $\Phi=(\varphi_{ij})_{n\times n}$,
\begin{itemize}
    \item denote $\Phi^{[k]}$ as the upper left $k\times k$
    submatrix, with eigenvalues
    $\lambda^{(k)}_1,\ldots,\lambda^{(k)}_k$;
    \item denote $\Phi^{[k]}_{\hat{k}k}$ as 
    the $k$-th column of $\Phi^{[k]}$
    but delete the $k$-th row;
    \item denote $\Phi^{[k]}_{k\hat{k}}$ as 
    the $k$-th row of $\Phi^{[k]}$
    but delete the $k$-th column.
\end{itemize}
For example, we have $\Phi^{[k]}=
\begin{pmatrix}
\Phi^{[k-1]} & \Phi^{[k]}_{\hat{k}k} \\
\Phi^{[k]}_{k\hat{k}} &  \varphi_{kk}
\end{pmatrix}$.
\end{nota}
\begin{defi}
\label{Def:vcat}
For $n\times n$ matrix $\Phi_0$ and
$\theta=(\theta_1,\ldots,\theta_n)\in\mathbb{R}^n$,
when $\Phi_0^{[2]},\ldots,\Phi_0^{[n-1]}$ is non-resonant,
we denote
\begin{align*}
\nu_{\tmop{cat}(\theta),d}(\Phi_0)
& \assign
\tmop{Ad}\left(
\overrightarrow{\prod^n_{k=1}}
C_{d+\theta_k}(E_k,\delta_k\Phi_0)
\right)
\mathe^{2\pi\mathi \Phi_0},\\
\tmop{MD}_{ \tmop{cat}(\theta),d } (\Phi_0)
& \assign
(\nu_{\tmop{cat}(\theta),d}(\Phi_0),
\sigma(\Phi_0),
\delta\Phi_0),\\
\mathscr{G}_{ + }
& \assign
\{ \tmop{MD}_{ \tmop{cat}(0), -\frac{\pi}{2} } (\Phi_0)
~|~ \Phi_0\in \mathfrak{c}_0 \}.
\end{align*}
\end{defi}
\begin{rmk}
Here the symbol ${\rm cat}(\theta)$ follows from the fact that the limit $\frac{u_{k+1}-u_{k}}{u_{k}-u_{k-1}}\rightarrow \infty$ for $k=2,\ldots,n$ is a caterpillar point on the De Concini-Procesi space of the space $\mathbb{C}^n\setminus \Delta$. See \cite{Xu} for more details.
\end{rmk}

\begin{lem}
\label{Lem:MtoExpM}
We have the following bijection 
\begin{align*}
(\mathe^{2\pi \mathi (\cdot)},\sigma):
\mathfrak{gl}_n^{\tmop{nr}}(\mathbb{C}) 
& \rightarrow
\{(V,\Sigma)
~|~
V\in{\rm GL}_n(\mathbb{C}),
\Sigma\subseteq \mathbb{C},
\sigma(V) = \mathe^{2\pi \mathi \Sigma},
\ \text{and} \
(\Sigma - \Sigma) \cap
\mathbb{Z} = \{0\}
\},\\
M & \mapsto
(\mathe^{2\pi \mathi M},
\sigma(M)),
\end{align*}
where we denote
$\mathfrak{gl}_n^{\tmop{nr}}(\mathbb{C}) \assign
\{A\in \mathfrak{gl}_n(\mathbb{C}) ~|~
A \text{ is non-resonant}\}$.
\end{lem}

\begin{lem}
\label{Lem:CoeftoMonoInj}
The monodromy mapping
$\tmop{MD}_{u,d}$ 
defined in \eqref{Def:MD}
is an injection on $\mathfrak{gl}_n^{\tmop{nr}}(\mathbb{C})$.
\end{lem}
\begin{prf}
Let us take
\begin{align*}
F_{d}(\xi;u,A)  = 
H_{d}(\xi;u,A) \xi^{\delta A}\mathe^{u \xi},\quad
F^{[\tmop{reg}]}(\xi;u,A)  =
H^{[\tmop{reg}]}(\xi;u,A) \xi^A,
\end{align*}
as the canonical fundamental solution of 
the confluent hypergeometric system
\begin{align*}
\frac{\mathd}{\mathd \xi} F(\xi) =
(u+A \xi^{-1})\cdot F(\xi),
\end{align*}
at $\xi = \infty$ and $\xi=0$, respectively. 
Let us introduce 
the monodromy matrix $M_d(u,A)$ 
associated to $F_d(\xi)$ by
\begin{align}
\label{Def:MdAd}
M_d(u,A) \assign
C_d (u, A) \cdot A
\cdot C_d (u, A)^{- 1}.
\end{align}
Now suppose that 
\[(\nu_{u,d}(A_1), \sigma(A_1), \delta A_1)
= (\nu_{u,d}(A_2), \sigma(A_2), \delta A_2).\]
From Lemma \ref{Lem:MtoExpM} 
and \eqref{Md=ConjPhi} in
Proposition \ref{AppPro:FormSC},
we have
$M_d(u,A_1)=M_d(u,A_2)$,
simply denoted by $M_d$.
According to \eqref{Def:MdAd} we have
\begin{align*}
F_{d}(\xi;u,A_1) F_{d}(\xi;u,A_2)^{-1} & =
H_{d}(\xi;u,A_1) H_{d}(\xi;u,A_2)^{-1}\\
& =
(H^{[\tmop{reg}]}(\xi;u,A_1) {C}_d(u,A_1)^{-1}
\xi^{M_d})
(H^{[\tmop{reg}]}(\xi;u,A_1) {C}_d(u,A_2)^{-1}
\xi^{M_d})^{-1}\\
& =
H^{[\tmop{reg}]}(\xi;u,A_1) {C}_d(u,A_1)^{-1}
{C}_d(u,A_2) H^{[\tmop{reg}]}(\xi;u,A_1)^{-1}.
\end{align*}
Therefore $F(\xi)\assign F_{d}(\xi;u,A_1) F_{d}(\xi;u,A_2)^{-1}$
is analytic on 
the entire Riemann sphere $\mathbb{C}\cup\{\infty\}$,
and we have
\begin{align*}
{I} = F(\infty) = F(0) = {C}_d(u,A_1)^{-1}{C}_d(u,A_2).
\end{align*}
It concludes the proof.
\end{prf}

\begin{lem}
\label{Lem:catComp} 
Denote $c_0 = 1$, for $c = (c_1, \ldots, c_n) \in (\mathbb{C}\backslash \{ 0
  \})^n$ we have
  \begin{align*}
    \tmop{Ad} \left( \overrightarrow{\prod_{k = 1}^n} C_{d +
    \theta_k} (c_k E_k, \delta_k \Phi_0) \right) \mathe^{2 \pi \mathi \Phi_0}
    & = 
    \nu_{\tmop{cat} (\theta + \arg c), d} \left( \left(
    \overrightarrow{\prod_{k = 1}^n} \left( \frac{c_k}{c_{k - 1}}
    \right)^{\tmop{ad} \delta_{k - 1} \Phi_0} \right) \Phi_0 \right).
  \end{align*}
\end{lem}

\begin{prf}
It can be verified directly from 
Proposition \ref{AppPro:FormSC}.
\end{prf}


\begin{pro}
\label{Pro:nrGenCri}
Suppose that $\Phi_n(z)$ is a solution of 
the $n$-th isomonodromy equation \eqref{zisoeq}
with non-resonant initial value,
and \begin{align*}
\tmop{Im}(u_1\mathe^{\mathi d})
    >\cdots>
\tmop{Im}(u_n\mathe^{\mathi d}),
\end{align*}
then $\Phi_n(z)$ is a shrinking solution
if and only if $\tmop{MD}_{u,d}(\Phi_n)\in
\mathscr{G}_{+}$.
\end{pro}

\begin{prf}
According to Lemma \ref{Lem:catComp} and
Theorem \ref{Thm:Cat}, 
we can assume that $d=-\frac{\pi}{2}$ and
\begin{align*}
\arg(u_k - u_{k-1} ) = 0,\quad
\text{for every }2\leqslant k \leqslant n,
\end{align*}
without loss of generality.
Therefore when $\Phi_0\in\mathfrak{c}_0$ we have
\begin{align*}
\tmop{MD}_{u,d}
(\mathcal{s}_n\cdots\mathcal{s}_1(\Phi_0)) =
\tmop{MD}_{\tmop{cat}(0), -\frac{\pi}{2} }(\Phi_0)
\in \mathscr{G}_{+} .
\end{align*}
Conversely, if there is a $\Phi_0\in\mathfrak{c}_0$ 
such that 
\begin{align*}
\tmop{MD}_{u,d} (\Phi_n)=
\tmop{MD}_{\tmop{cat}(0), -\frac{\pi}{2} }(\Phi_0) =
\tmop{MD}_{u,d} (\mathcal{s}_n\cdots\mathcal{s}_1(\Phi_0)),
\end{align*}
according to Lemma \ref{Lem:CoeftoMonoInj}, 
we have $\Phi_n =
\mathcal{s}_n\cdots\mathcal{s}_1(\Phi_0)$. It concludes the proof.
\end{prf}


\subsection{Inverse Monodromy Problem}
\label{subSect:GP}


This subsection will need to establish 
the inverse mapping of 
the Riemann-Hilbert mapping \eqref{Cat}
provided in Theorem \ref{Thm:Cat}. It will thus provide a sufficient condition 
for an element that belongs to $\mathscr{G}_{+}$. 
Although it only characterizes 
a proper subset of $\mathscr{G}_{+}$, 
it is enough
to prove the generic property of the shrinking solutions.

\begin{lem}
\label{Pro:SolveP0}
If there exists
\begin{align*}
V \in  
\tmop{GL}_n(\mathbb{C}),\quad
\sigma_k = 
(\lambda^{(k)}_j)_{j=1,\ldots,k},\quad
\Lambda = 
\tmop{diag}(\varphi_{11},\ldots,\varphi_{nn}),
\end{align*}
such that for $1\leqslant k\leqslant n-1$,
    \begin{subequations}
    \begin{eqnarray}
    \label{Cond:NonResk}
    \lambda^{(k)}_{j_1}-\lambda^{(k)}_{j_2}
    & \notin
    \mathbb{Z}\setminus\{0\},
    \quad
    \text{for every $j_1,j_2$},\\
    \label{Cond:NonResNear}
    \lambda^{(k+1)}_{j_1}-\lambda^{(k)}_{j_2}
    & \notin
    \mathbb{Z}\setminus\{0\},
    \quad
    \text{for every $j_1,j_2$},
    \end{eqnarray}
    \end{subequations}
and for $1\leqslant k\leqslant n$,
    \begin{subequations}
    \begin{align}
    \label{sigk:Exp}
        \mathe^{2\pi\mathi \sigma_k}
         & =
        \sigma(V^{[k]}),\\
    \label{sigk:Sum}
        \sum_{j=1}^k \lambda^{(k)}_j
         & =
        \sum_{j=1}^k \varphi_{j j},
    \end{align}
    \end{subequations}
then there is a unique
$\Phi_0\in\mathfrak{gl}_n(\mathbb{C})$
such that
\begin{subequations}
\begin{align}
\tmop{MD}_{\tmop{cat}(0), -\frac{\pi}{2} } (\Phi_0) 
& =
(V,\sigma_n,\Lambda)\in
\mathscr{G}_{+},\\
\label{Cond:sigPk}
\sigma(\Phi_0^{[k]}) 
& =
\sigma_k,\quad
\text{for every $1\leqslant k \leqslant n$}.
\end{align}
\end{subequations}
\end{lem}

\begin{prf}
  Denote $d = - \frac{\pi}{2}$. Suppose that we have recursively calculated
  $\Phi_0^{[k - 1]}$ such that it 
  satisfies condition \eqref{Cond:sigPk}, 
  and
  \begin{align}
  \label{Ck-1rec}
    V^{[k - 1]} = 
    \tmop{Ad} \left( \overrightarrow{\prod_{j = 1}^{k - 1}}
    C_d (E_j, \delta_j \Phi_0^{[k - 1]}) \right) \mathe^{2 \pi \mathi
    \Phi_0^{[k - 1]}}.
  \end{align}
  For simplicity, let us denote
  \begin{align*}
    c_{k - 1} =  
    \overrightarrow{\prod_{j = 1}^{k - 1}} C_d (E_j,
    \delta_j \Phi_0^{[k - 1]}), \quad
    M_{k - 1} =  
    \tmop{Ad} (c_{k - 1}) \Phi_0^{[k - 1]}.
  \end{align*}
  Without loss of generality,
  we additionally
  assume that $\Phi_0^{[k - 1]}$ 
  has distinct eigenvalues,
  with spectral decomposition of $M_{k - 1}$
  \begin{align*}
    M_{k - 1} = 
    P_1 \lambda^{(k - 1)}_1 + \cdots + P_{k - 1} \lambda^{(k
    - 1)}_{k - 1}.
  \end{align*}
Therefore we naturally have
  \begin{align}
  \label{CPCeigen}
    V^{[k]}_{k \hat{k}} P_i 
    V^{[k]}_{\hat{k} k}  =  - \prod_{p \neq i}
    \prod_j \frac{\mathe^{2 \pi \mathi \lambda^{(k - 1)}_i} - \mathe^{2 \pi
    \mathi \lambda^{(k)}_j}}{\mathe^{2 \pi \mathi \lambda^{(k - 1)}_i} -
    \mathe^{2 \pi \mathi \lambda^{(k - 1)}_p}} .
  \end{align}
Next, based on the \eqref{Ck-1rec} 
(provided by replacing $k-1$ as $k$)
which we want to verify, 
and the closed formula of 
the Stokes matrices $S_d^\pm (E_k,\delta_k \Phi_0^{[k]})$ 
given in Proposition \ref{Pro:S},
let us construct $\Phi_0^{[k]} =  
    \left(\begin{array}{cc}
      \Phi_0^{[k - 1]} & (\Phi_0)^{[k]}_{\hat{k} k}\\
      (\Phi_0)^{[k]}_{k \hat{k}} & \varphi_{k k}
    \end{array}\right)$
recursively for $k$,
under condition \eqref{Cond:NonResNear},
where
\begin{subequations}
  \begin{align}
  \label{Form:Vcol}
    (\Phi_0)^{[k]}_{\hat{k} k} & :=  
    \left( \frac{1}{2 \pi \mathi} \prod_{p, j}
    \frac{(\Phi_0^{[k - 1]} - \lambda^{(k)}_j I) !}{(\Phi_0^{[k - 1]} -
    \lambda^{(k - 1)}_p I) !} \right) \mathe^{- 2 \pi \mathi \Phi_0^{[k - 1]}}
    c_{k - 1}^{- 1} \cdot V^{[k]}_{\hat{k} k}, \\
  \label{Form:Vrow}
    (\Phi_0)^{[k]}_{k \hat{k}} & :=  
    V^{[k]}_{k \hat{k}} \cdot c_{k - 1} \mathe^{-
    \pi \mathi (\Phi_0^{[k - 1]} + \varphi_{k k} I)} \left( \frac{1}{2 \pi \mathi}
    \prod_{p, j} \frac{(\lambda^{(k)}_j I - \Phi_0^{[k - 1]}) !}{(\lambda^{(k -
    1)}_p I - \Phi_0^{[k - 1]}) !} \right) . 
  \end{align}
\end{subequations}
From the reflection formula 
of $x!\assign \Gamma(x+1)$
and \eqref{CPCeigen},
we have
  \begin{align*}
    &\quad (\Phi_0)^{[k]}_{k \hat{k}} 
    \frac{1}{z I - \Phi_0^{[k - 1]}}
    (\Phi_0)^{[k]}_{\hat{k} k} \nonumber\\
    & =  V^{[k]}_{\hat{k} k} \mathe^{\pi \mathi (- 3 M_{k - 1} - \varphi_{k k} I)}
    \left( \left( \frac{1}{2 \pi \mathi} \right)^2 \prod_{p, j}
    \frac{(\lambda^{(k)}_j I - M_{k - 1}) ! (M_{k - 1} - \lambda^{(k)}_j I)
    !}{(\lambda^{(k - 1)}_p I - M_{k - 1}) ! (M_{k - 1} - \lambda^{(k - 1)}_p
    I) !} \right) \frac{1}{z I - M_{k - 1}} V^{[k]}_{\hat{k} k} \nonumber\\
    & =  V^{[k]}_{\hat{k} k} \left( \mathe^{- 2 \pi \mathi M_{k - 1}}
    \prod_{p, j} \frac{\mathe^{2 \pi \mathi M_{k - 1}} - \mathe^{2 \pi \mathi
    \lambda^{(k - 1)}_p} I}{\mathe^{2 \pi \mathi M_{k - 1}} - \mathe^{2 \pi
    \mathi \lambda^{(k)}_j} I} \cdot \frac{M_{k - 1} - \lambda^{(k)}_j I}{M_{k
    - 1} - \lambda^{(k - 1)}_p I} \right) \frac{1}{z I - M_{k - 1}}
    V^{[k]}_{\hat{k} k} \nonumber\\
    & =  \sum_{i = 1}^{k - 1} V^{[k]}_{k \hat{k}} P_i V^{[k]}_{\hat{k} k}
    \cdot \left( \prod_{p \neq i} \prod_j \frac{\mathe^{2 \pi \mathi
    \lambda^{(k - 1)}_i} - \mathe^{2 \pi \mathi \lambda^{(k -
    1)}_p}}{\mathe^{2 \pi \mathi \lambda^{(k - 1)}_i} - \mathe^{2 \pi \mathi
    \lambda^{(k)}_j}} \cdot \frac{\lambda^{(k - 1)}_i -
    \lambda^{(k)}_j}{\lambda^{(k - 1)}_i - \lambda^{(k - 1)}_p} \right)
    \frac{1}{z - \lambda^{(k - 1)}_i} \nonumber\\
    & =  - \sum_{i = 1}^{k - 1} \left( \prod_{p \neq i} \prod_j
    \frac{\lambda^{(k - 1)}_i - \lambda^{(k)}_j}{\lambda^{(k - 1)}_i -
    \lambda^{(k - 1)}_p} \right) \frac{1}{z - \lambda^{(k - 1)}_i}. 
  \end{align*}
Therefore
  \begin{align*}
    \frac{\det (z I - \Phi_0^{[k]})}{\det (z I - \Phi_0^{[k - 1]})} = (z - \varphi_{k
    k}) - 
    (\Phi_0)^{[k]}_{k \hat{k}} 
    \frac{1}{z I - \Phi_0^{[k - 1]}}
    (\Phi_0)^{[k]}_{\hat{k} k} = \prod_{p, j} \frac{z - \lambda^{(k)}_j}{z -
    \lambda^{(k - 1)}_p},
  \end{align*}
and thus the set of eigenvalues $\sigma (\Phi_0^{[k]}) = \sigma_k$. 
From \eqref{Ck-1rec} and 
Proposition \ref{Pro:S} we have
  \begin{align}
  \nonumber
    V^{[k]} & =  \tmop{Ad} \left( \left(\begin{array}{cc}
      c_{k - 1} & 0\\
      0 & 1
    \end{array}\right) \right) \left( S_d^- (E_k, \Phi_0^{[k]})^{- 1}
    \left(\begin{array}{cc}
      \mathe^{2 \pi \mathi \Phi_0^{[k - 1]}} & 0\\
      0 & \mathe^{2 \pi \mathi \varphi_{k k}}
    \end{array}\right) S_d^+ (E_k, \Phi_0^{[k]}) \right) \nonumber\\
    & =  \tmop{Ad} \left( \overrightarrow{\prod_{j = 1}^{k - 1}}
    C_d (E_j, \delta_j \Phi_0^{[k]}) \right) \cdot \tmop{Ad}
    (C_d (E_k, \Phi_0^{[k]})) \mathe^{2 \pi \mathi \Phi_0^{[k]}}
    \nonumber\\
    & =  \tmop{Ad} \left( \overrightarrow{\prod_{j = 1}^k} C_d
    (E_j, \delta_j \Phi_0^{[k]}) \right) \mathe^{2 \pi \mathi \Phi_0^{[k]}} .
    \label{VReck}
  \end{align}
By induction, we get 
\[V=  \tmop{Ad} \left( \overrightarrow{\prod_{j = 1}^n} C_d
    (E_j, \delta_j \Phi_0) \right) \mathe^{2 \pi \mathi \Phi_0}. \]
Following Definition \ref{Def:vcat}, we see that $\tmop{MD}_{\tmop{cat}(0), -\frac{\pi}{2} } (\Phi_0) 
=
(V,\sigma_n,\Lambda)$. 
Since condition \eqref{Cond:NonResk} 
forces $\Phi_0$ to satisfy
\eqref{Form:Vcol}, \eqref{Form:Vrow},
thus guarantees the uniqueness of $\Phi_0$,
we complete the proof.
\end{prf}

\begin{lem}
\label{Lem:LogSq}
Given a sequence
$c=(c_i)_{i=1,\ldots,k}$ of complex numbers
and a constant $s$ 
such that
\begin{align*}
\mathe^{2\pi\mathi s} = c_1\cdots c_k,
\end{align*}
then there is a sequence
$\sigma=(\lambda_j)_{j=1,\ldots,k}$
such that
$\mathe^{2\pi\mathi \sigma} = c$,
$\sum_{j=1}^k \lambda_j = s$
and
\begin{eqnarray}
\label{Cond:ConsBoundy}
|\tmop{Re}(\lambda_{j_1}-\lambda_{j_2})|
\leqslant 1, \quad
\text{for every $j_1,j_2$}.
\end{eqnarray}
If we additionally have
\begin{eqnarray}
\label{Cond:ConsLem}
|\tmop{Re}(\lambda_{j_1}-\lambda_{j_2})|
 <  1, \quad
\text{for every $j_1,j_2$},
\end{eqnarray}
then the sequence $\sigma$ satisfying 
condition \eqref{Cond:ConsBoundy} is unique.
\end{lem}

\begin{prf}
Inductively assume that 
we have constructed
the sequence $\Sigma_m=(\lambda_{j,m})_{j=1,\ldots,k}$ 
satisfying
\begin{align}
\label{Recsm}
\mathe^{2\pi\mathi \Sigma_m} = c,
\quad
\sum_{j=1}^k \lambda_{j,m} = s.
\end{align}
Assume that
\begin{align*}
\tmop{Re} \lambda_{j_1,m}  = 
\underset{1\leqslant j\leqslant k}{\tmop{min}}
\tmop{Re} \lambda_{j,m},\quad
\tmop{Re} \lambda_{j_2,m}  = 
\underset{1\leqslant j\leqslant k}{\tmop{max}}
\tmop{Re} \lambda_{j,m}.
\end{align*}
If we have $\tmop{Re}(
\lambda_{j_2,m}-
\lambda_{j_1,m}) >1$,
then we can construct the sequence 
$\Sigma_{m+1}=(\lambda_{j,m+1})_{j=1,\ldots,k}$ as
\begin{align*}
  \lambda_{j, m + 1} = \left\{
  \begin{array}{ll}
    \lambda_{j, m} & ; j \neq j_1, j_2\\
    \lambda_{j, m} + 1 & ; j = j_1\\
    \lambda_{j, m} - 1 & ; j = j_2
  \end{array}\right..
\end{align*}
Therefore the sequence
$\Sigma_{m+1}$ also satisfy \eqref{Recsm},
and we either have
\begin{align*}
\underset{1\leqslant j_1,j_2\leqslant k}{\tmop{max}}
\tmop{Re}(
\lambda_{j_2,m+1}-
\lambda_{j_1,m+1}) <
\underset{1\leqslant j_1,j_2\leqslant k}{\tmop{max}}
\tmop{Re}(
\lambda_{j_2,m}-
\lambda_{j_1,m}),
\end{align*}
or
\begin{align*}
  \# \left\{ j_1 : \tmop{Re} \lambda_{j_1, m + 1} = \underset{1 \leqslant j
  \leqslant k}{\min} \tmop{Re} \lambda_{j, m + 1} \right\} & < \# \left\{
  j_1 : \tmop{Re} \lambda_{j_1, m} = \underset{1 \leqslant j \leqslant
  k}{\min} \tmop{Re} \lambda_{j, m} \right\}, \\
  \# \left\{ j_2 : \tmop{Re} \lambda_{j_2, m + 1} = \underset{1 \leqslant j
  \leqslant k}{\max} \tmop{Re} \lambda_{j, m + 1} \right\} & < \# \left\{
  j_2 : \tmop{Re} \lambda_{j_2, m} = \underset{1 \leqslant j \leqslant
  k}{\max} \tmop{Re} \lambda_{j, m} \right\}.
\end{align*}
Therefore, for a sufficiently large ${N}$,
we eventually have $\underset{1\leqslant j_1,j_2\leqslant k}{\tmop{max}}
\tmop{Re}(
\lambda_{j_2,{N}}-
\lambda_{j_1,{N}}) \leqslant 1$.
By taking $\sigma=\Sigma_N$,
we finish the proof.
\end{prf}

\begin{thm}
\label{Pro:GenSGenc}
Almost every solution of the equation \eqref{isoeq} is a shrinking solution.
\end{thm}

\begin{prf}
Denote $\mathscr{G}_{+}^{\tmop{cl}}$
as the set of triples
$(V,\sigma_n,\Lambda)$
that satisfies the following conditions, 
\begin{enumerate}
    \item $\sigma_n$ is a complex sequence 
    $(\lambda_j^{(n)})_{j=1,\ldots,n}$,
    such that
    $\mathe^{2\pi\mathi \sigma_n}=(\mathe^{2\pi\mathi\lambda_j^{(n)}})_{j=1,\ldots,n}$
    is the eigenvalues of $V$;
    
    \item $\Lambda$ is a diagonal matrix
    $\tmop{diag}(\varphi_{11},\ldots,\varphi_{nn})$, 
    such that
    $\mathe^{2\pi\mathi
    \sum_{j=1}^k \varphi_{jj}}=
    \tmop{det}V^{[k]}$.
\end{enumerate}
It is also the closure of
the space $\mathscr{G}_{+}$ given in Definition \ref{Def:vcat}.
We introduce the following notation
\begin{align*}
\mathfrak{c}_0^{\mathrm{nr}}
& \assign 
\{\Phi_0\in\mathfrak{c}_0:
\text{$\Phi_0$ is non-resonant}
\},\\
\mathfrak{c}_n^{\mathrm{nr}}
& \assign 
\{\Phi_n(z)\in\mathfrak{c}_n:
\text{$\Phi_n(z)$ is non-resonant}
\},\\
\mathscr{G}_{+}^{\mathrm{nr}}
& \assign 
\{(V,\sigma_n,\Lambda)\in
\mathscr{G}_{+}:
(\sigma_n-\sigma_n)\cap
(\mathbb{Z}\setminus\{0\})
= \varnothing
\}.
\end{align*}
According to Proposition \ref{Pro:nrGenCri}, 
we have the following diagram.
\begin{eqnarray*}
\begin{tikzcd}
& \mathfrak{c}_n^{\mathrm{nr}} 
\arrow[dd, hook, "\cong"',
"{ \tmop{MD}_{u,d} }"] 
\arrow[rr, hook, "\subseteq"] &  & \text{non-resonant solutions} 
\arrow[dd, hook, 
"{ \tmop{MD}_{u,d} }"] \\
\mathfrak{c}_0^{\mathrm{nr}} 
\arrow[ru, "\cong"',
"\mathcal{s}_n\cdots\mathcal{s}_1"] 
\arrow[rd, "\cong",
"{\tmop{MD}_{\tmop{cat}(0),- \frac{\pi}{2}}}"']
&  &  & \\ &
\mathscr{G}_{+}^{\mathrm{nr}} 
\arrow[rr, hook, "\subseteq"]   &  &
\mathscr{G}_{+}^{\tmop{cl}}              
\end{tikzcd}
\end{eqnarray*}
According to Lemma \ref{Lem:LogSq}, 
for a generic element 
$(V,\sigma_n,\Lambda)$
in $\mathscr{G}_{+}^{\tmop{cl}}$
satisfying the conditions 
\eqref{sigk:Exp}, 
\eqref{sigk:Sum}, 
we can find a unique sequence 
$\sigma_k=(\lambda^{(k)}_j)_{j=1,\ldots,k}$
for $1\leqslant k\leqslant n-1$ such that
\begin{eqnarray}
|\tmop{Re}(
\lambda^{(k)}_{j_1}-
\lambda^{(k)}_{j_2})|
 <  1, \quad
\text{for every $j_1,j_2$}.
\end{eqnarray}
At the same time, 
condition \eqref{Cond:NonResk} 
is automatically satisfied, 
and generically there is also 
condition \eqref{Cond:NonResNear}, 
thus Lemma \ref{Pro:SolveP0} ensures that 
$\mathscr{G}_{+}^{\mathrm{nr}}$ 
is open dense in
$\mathscr{G}_{+}^{\tmop{cl}}$.
Finally, note that 
the space of non-resonant solutions is open dense, 
and the monodromy data mapping
$\tmop{MD}_{u,d}$ 
on this space is an analytic homeomorphism
to its image, 
thus the solution space $\mathfrak{c}_n^{\mathrm{nr}}$
is also open dense, 
and we complete the proof.
\end{prf}
The following theorem characterizes those monodromy data that come from shrinking solutions.
\begin{thm}
\label{Thm:GenisGen}
Suppose that $\Phi_n=(\varphi_{i j})_{n\times n}$ 
is a solution of 
the $n$-th isomonodromy equation \eqref{zisoeq},
satisfying the following generic conditions
\begin{itemize}

\item $\Phi_n$ has non-resonant initial value;

\item Take $V=\nu_{u,d}(\Phi_n)$, where
\begin{eqnarray}
\label{Cond:uPos}
\tmop{Im}(u_1\mathe^{\mathi d})
    >\cdots>
\tmop{Im}(u_n\mathe^{\mathi d}).
\end{eqnarray}
For every $1\leqslant k \leqslant n-1$,
there is a sequence
$\sigma_k=(\lambda^{(k)}_j)_{j=1,\ldots,k}$
satisfy that
\begin{align*}
\mathe^{2\pi\mathi \sigma_k}
=
\sigma(V^{[k]}),\quad
\sum_{j=1}^k \lambda^{(k)}_j
=
\sum_{j=1}^k \varphi_{j j},\quad
|\tmop{Re}(\lambda^{(k)}_{j_1}-\lambda^{(k)}_{j_2})|
< 1, \quad
\text{for every $j_1,j_2$};
\end{align*}

\item Denote $\sigma(\Phi_n)=
(\lambda^{(n)}_j)_{j=1,\ldots,n}$,
for every $1\leqslant k \leqslant n-1$ we have
\begin{align*}
\lambda^{(k+1)}_{j_1}-\lambda^{(k)}_{j_2} \notin \mathbb{Z},
\quad
\text{for every $j_1,j_2$},
\end{align*}
\end{itemize}
then $\Phi_n$ is a shrinking solution.
\end{thm}

\begin{prf}
The proof can be completed by combining 
Proposition \ref{Pro:nrGenCri} and 
Lemma \ref{Pro:SolveP0}. 
Note that
Theorem \ref{Pro:GenSGenc} ensures that 
these conditions are generic.
\end{prf}\\

\begin{rmk}
It can be seen that, 
for any $u$ and anti-Stokes direction $d$,
we can always make
a suitable index permutation 
such the condition
\eqref{Cond:uPos} hold.
\end{rmk}

\subsection{Summary and the Proof of Theorem \ref{Thm:termwise}}
In the end, we summarize the structure of Section \ref{Sect:Series}--Section \ref{Sect:GeneralSol} 
and provide the proof of Theorem \ref{Thm:termwise}, 
based on the following diagram
\[\begin{tikzcd}
  \cdots 
& [-00pt]{\mathfrak{s}_{k-1}(\bm{z}_{k-2})} 
& [-43pt]
& [-43pt]{\mathfrak{b}_{k-1}(\rm{pt})}
& [-10pt]{\mathfrak{s}_{k}(\bm{z}_{k-1})} 
& [-35pt]
& [-35pt]{\mathfrak{b}_{k}(\rm{pt})} 
& \cdots \\
  \cdots 
& {\mathfrak{s}_{k-1}} 
&
& {\mathfrak{b}_{k-1}} 
& {\mathfrak{s}_{k}} 
&
& {\mathfrak{b}_{k}} 
& \cdots \\
  \cdots 
&
& {\mathfrak{c}_{k-1}} 
&
&
& {\mathfrak{c}_{k}} 
&
& \cdots
	\arrow[from=3-3, to=2-2]
	\arrow[from=3-3, to=2-4]
	\arrow[from=3-6, to=2-5]
	\arrow[from=3-6, to=2-7]
	\arrow["{\mathcal{s}_k}",from=1-4, to=1-5]
	\arrow["{\mathcal{s}_k}",from=2-4, to=2-5]
	\arrow["{\mathcal{s}_k}",from=3-3, to=3-6]
	\arrow[from=1-4, to=2-4]
	\arrow[from=1-5, to=2-5]
	\arrow[from=1-7, to=2-7]
	\arrow[from=1-2, to=2-2]
	\arrow["{\mathcal{s}_{k+1}}",from=1-7, to=1-8]
	\arrow["{\mathcal{s}_{k+1}}",from=2-7, to=2-8]
	\arrow["{\mathcal{s}_{k+1}}",from=3-6, to=3-8]
	\arrow["{\mathcal{s}_{k-1}}",from=3-1, to=3-3]
	\arrow["{\mathcal{s}_{k-1}}",from=2-1, to=2-2]
	\arrow["{\mathcal{s}_{k-1}}",from=1-1, to=1-2]
\end{tikzcd}\]
In Sections \ref{step1} and \ref{step2}, we construct
the mapping $\mathfrak{b}_{k-1}(\tmop{pt})
\xrightarrow{\mathcal{s}_k}
\mathfrak{s}_k(\bm{z}_{k-1})$. Here recall the spaces $\mathfrak{b}, \mathfrak{c}$ and $\mathfrak{s}$ are defined in Notation \ref{bcs}.
In Section \ref{step3}, 
we expand it into 
the mapping $\mathfrak{b}_{k-1}
\xrightarrow{\mathcal{s}_k}
\mathfrak{s}_k$, that is Theorem \ref{kmainthm}.
In Sections \ref{SubSect:SolAsy} and 
\ref{SubSect:EffCri}, 
we prove that 
the mapping $\mathfrak{b}_{k-1}(\tmop{pt})
\xrightarrow{\mathcal{s}_k}
\mathfrak{s}_k(\bm{z}_{k-1})$ is a bijection. 
Based on the existence and uniqueness of the solution with respect to the initial value, 
we see that 
the mapping $\mathfrak{b}_{k-1}
\xrightarrow{\mathcal{s}_k}
\mathfrak{s}_k$ is also a bijection, 
that is Theorem \ref{Thm:AsySol}.
According to Proposition \ref{Pro:TrivalFirst}, 
for each $0\leqslant k < n$, 
we have $\mathfrak{c}_k\subseteq \mathfrak{b}_k$.
Therefore,
for each $0< k\leqslant n$, 
we have $\mathfrak{c}_k\subseteq \mathfrak{s}_k$.
Thus we prove the uniqueness part 
of Theorem \ref{Thm:termwise}.

Furthermore, we can use Theorem \ref{Thm:Cat} to compute the monodromy data 
of the shrinking solutions 
with respect to the boundary values. See Appendix \ref{Sect:Diag} or the formula in Theorem \ref{thm: introcatformula}. The inverse formula is given in Theorem \ref{Thm:GenisGen}, from which we see that  the set of monodromy corresponding to shrinking solutions are dense in the space of all possible monodromy data.
Since the generic solutions and the monodromy data 
are one-to-one correspondence, 
we eventually prove
that almost every solution is a shrinking solution, 
thus complete the proof of Theorem \ref{Thm:termwise}.

  \section{Examples of Non-generic Solutions and Applications}
\label{Sect:ExAp}

\subsection{Non-Generic Solutions}
\label{subSect:NGS}

Although almost all the solutions of 
the isomonodromy equation \eqref{isoeq} are shrinking solutions, 
there are still many solutions that 
are not in the generic class. 
In this subsection, we will introduce some examples.

\subsubsection{Upper triangular matrix solutions}

The most simplest family of solutions is 
the upper triangular solution. 
It will contains many non-shrinking solutions, 
but the coefficients can still be given by 
Proposition \ref{Lem:TotallyFormalPhik}.

\begin{ex}[Constant Solution]
For $\Phi_{k-1}=\delta_{k-1}\Phi_{k-1}$ 
without the boundary condition \eqref{Cond:ConstraEigen}, 
we can verify that
the series $\Phi_k$ is trivial,
that is
\begin{eqnarray}
\Phi_k (z_1,\ldots,z_{k-1},z_k)=
\Phi_{k-1}(z_1,\ldots,z_{k-1}).
\end{eqnarray}
Furthermore,
for any $m\geqslant k$ we have $\Phi_m = \Phi_{k-1}$.
\end{ex}

Suppose that a solution of the $n$-th isomonodromy equation takes the form
$\Phi_n = 
\left(\begin{array}{cc}
    \Phi_n^{[n - 1]} & (\Phi_n)_{\hat{n} n}\\
    0 & a_{n n}
\end{array}\right)$,
then we have
\begin{align*}
  \frac{\partial}{\partial z_n} 
  \Phi_n^{[n-1]}
  & =  
  0, \\
  \frac{\partial}{\partial z_n} 
  \left(\begin{array}{c}
    (\Phi_n)_{\hat{n} n}\\
    a_{n n}
        \end{array}\right)
  & =  
  (a_{n n} I - \delta_{n
  - 1} \Phi_n) (z_1 \cdots z_{n - 1} D_n^{(n)}) \cdot   \left(\begin{array}{c}
    (\Phi_n)_{\hat{n} n}\\
    a_{n n}
        \end{array}\right). 
\end{align*}
Here we take the Notation \ref{Nota:Mat}. 
Therefore $(\Phi_n)_{\hat{n} n}$ satisfies the following linear system with
respect to $z_n$ at $z_n = \infty$
\begin{align}
  \frac{\partial}{\partial z_n} (\Phi_n)_{\hat{n} n} 
  & = 
  ((a_{n n} I -
  \Phi_n^{[n - 1]}) z_n^{- 1} + O (z_n^{- 2})) (\Phi_n)_{\hat{n} n} . 
  \label{Sys:Pncol}
\end{align}
All of the non-resonant solutions $\Phi_n$ in this form, can be obtained by
the analytic continuation of $\mathcal{s}_n$. 
More precisely, we have

\begin{pro}
  Suppose that $\Phi_{n - 1}^{[n - 1]}$ is non-resonant, denote
  \begin{align*}
    \Phi_{n - 1} & = 
    \left(\begin{array}{cc}
      \Phi_{n - 1}^{[n - 1]} & (\Phi_{n - 1})_{\hat{n} n}\\
      0 & a_{n n}
    \end{array}\right) . 
  \end{align*}
  The mapping $\mathcal{s}_n$
  defined on $\mathfrak{b}_{n-1}$ can be 
  analytically continued 
  to the $\Phi_{n - 1}$
  in the above form,
  and we have
  \begin{align*}
    \mathcal{s}_n (\Phi_{n - 1}) & = 
    \left(\begin{array}{cc}
      \Phi_{n - 1}^{[n - 1]} & (\mathcal{s}_n (\Phi_{n - 1}))_{\hat{n} n}\\
      0 & a_{n n}
    \end{array}\right), \\
    (\mathcal{s}_n (\Phi_{n - 1}))_{\hat{n} n} & = 
    \left( I + \sum_{p =
    1}^{\infty} z_n^{- p} F_p \right) z_n^{a_{n n} I - \Phi_{n - 1}^{[n - 1]}}
    (\Phi_{n - 1})_{\hat{n} n}, 
  \end{align*}
  where $I + \sum_{p = 1}^{\infty} z_n^{- p} F_p$ is the canonical fundamental
  solution of \eqref{Sys:Pncol} at $z_n = \infty$.
\end{pro}

\begin{prf}
  It is easy to see that the mapping $\mathcal{s}_n$ defined in this way is
  analytic, so we can assume that $\Phi_{n - 1}$ 
  satisfies condition \eqref{Cond:ConstraEigen}. 
  We can verify that 
  $\mathcal{s}_n (\Phi_{n - 1})$ 
  is exactly the series solution of 
  $n$-th isomonodromy equation, 
  and we have
  \begin{align*}
    z_n^{\tmop{ad} \delta_{n - 1} \Phi_{n - 1}} \phi_{n, m} (z_n) & =
    \left(\begin{array}{cc}
      0 & (z_n^{\tmop{ad} \Phi_{n - 1}^{[n - 1]}} F_m)
      \cdot
      (\Phi_{n - 1})_{\hat{n}
      n}\\
      0 & 0
    \end{array}\right), \quad m \geqslant 1, 
  \end{align*}
  thus we complete the proof.
\end{prf}\\
It can be seen that
$\mathcal{s}_n$
cannot be continued to the case that
$\Phi_{n - 1}^{[n - 1]}$ is resonant,
and it will be the singularity of $\mathcal{s}_n$. 
In order to parameterize the whole solution space of
the $n$-th isomonodromy equation, 
it is necessary to investigate 
the blow up at such singularity.

\subsubsection{Rational non-generic solutions}
Another interesting family of non-shrinking solutions are 
the following rational solutions. Recall that we denote $u={\rm diag}(u_1,...,u_n)$.
\begin{pro}
\label{Pro:Rational}
Suppose that $a$ and $b$ are $n \times m$ and $m \times n$ matrices such that $b a = 0$. If $b u a$ is invertible (thus $\tmop{rank}a=
    \tmop{rank}b=m
    \leqslant \frac{n}{2}$),
then the $n\times n$ matrix function
\begin{align*}
\Phi_n(u)
& =
\tmop{ad}_u \left(a \frac{1}{b u a} b \right)
\end{align*}
is a rational solution of
  the $n$-th isomonodromy equation \eqref{isoeq}.
The solution $\Phi_n(u)$ is diagonalizable and
has eigenvalues $1, -1, 0$ 
with multiplicities $m,m,n-2m$
and spectral decomposition
\begin{subequations}
\begin{align}
\label{Phin1part}
\mathcal{P}_1 \Phi_n (u)
& = 
\left( I - \frac{1}{2}
a \frac{1}{b u a} b u
\right) 
\cdot 
\left( u a \frac{1}{b u a} b \right) 
= u a \frac{1}{b u a} b -
  a \frac{1}{b u a} \frac{b u^2 a}{2} 
  \frac{1}{b u a} b,\\
\label{Phin-1part}
\mathcal{P}_{- 1} \Phi_n (u)
& = 
\left( a \frac{1}{b u a} b u \right) 
\cdot
\left( I - \frac{1}{2} 
u a \frac{1}{b u a} b \right) 
= a \frac{1}{b u a} b u -
  a \frac{1}{b u a} \frac{b u^2 a}{2} 
  \frac{1}{b u a} b,\\
\label{Phin0part}
\mathcal{P}_0 \Phi_n (u)
& = \left( I - a \frac{1}{b u a} b u \right) \cdot
  \left( I - u a \frac{1}{b u a} b \right).
\end{align}
\end{subequations}
\end{pro}

\begin{prf}
Using $\frac{\partial}{\partial u_k} \left( \frac{1}{b u a} \right) = -
  \frac{1}{b u a} b E_k a \frac{1}{b u a}$, we can verify that
\begin{align*}
    \frac{\partial}{\partial u_k} \tmop{ad}_u \left( a \frac{1}{b u a} b
    \right) 
& = \tmop{ad}_{E_k} \left( a \frac{1}{b u a} b \right) -
    \tmop{ad}_u \left( a \frac{1}{b u a} b E_k a \frac{1}{b u a} b \right), \\
    \left[ \tmop{ad}_{E_k} \left( a \frac{1}{b u a} b \right), \tmop{ad}_u
    \left( a \frac{1}{b u a} b \right) \right] 
& = \tmop{ad}_{E_k} \left( a
    \frac{1}{b u a} b u a \frac{1}{b u a} b \right) - \tmop{ad}_u \left( a
    \frac{1}{b u a} b E_k a \frac{1}{b u a} b \right),
\end{align*}
and thus $\Phi_n(u)$ is 
indeed a solution of 
the $n$-th isomonodromy equation.

Since it can be directly verified that 
$\Phi_n(u)^3=\Phi_n(u)$, 
therefore $\Phi_n(u)$ is diagonalizable and 
has eigenvalues $1, -1, 0$,
and we naturally have
\eqref{Phin1part}, \eqref{Phin-1part}
and \eqref{Phin0part}.
Note that we have the spectrum
\begin{align*}
\sigma\left( u a \frac{1}{b u a} b \right)=
\sigma\left( a \frac{1}{b u a} b u \right)=
\{0,1\},
\end{align*}
therefore the multiplicities of $1$ and $-1$
respectively are
\begin{align*}
\tmop{rank} \mathcal{P}_1 \Phi_n(u)
& =
\tmop{rank} 
\left( u a \frac{1}{b u a} b \right)=m,\\
\tmop{rank} \mathcal{P}_{-1} \Phi_n(u)
& =
\tmop{rank} 
\left( a \frac{1}{b u a} b u \right)=m.
\end{align*}
It concludes the proof.
\end{prf}

Although we can verify that
\begin{align*}
\lim_{z_n\to\infty}\Phi_n^{[n-1]}(u)=
\tmop{ad}_{ u^{[n-1]} }
\left(\Tilde{a}\frac{1}
{ \Tilde{b} u^{[n-1]} \Tilde{a} }
\Tilde{b}\right),
\end{align*}
but when $m>1$, the
$\lim_{z_n\to\infty}
z_n^{\ad\mathcal{l}_n(\delta_{n-1}\Phi_n)}\Phi_n$
may not exists,
therefore $\Phi_n(u)$ is not a shrinking solution.

The case $m=1$ and
\[ a = \left(\begin{array}{c}
     a_1\\
     \vdots\\
     a_{n - 1}\\
     1
   \end{array}\right), \quad b = \left(\begin{array}{cccc}
     b_1 & \cdots & b_{n - 1} & 1
   \end{array}\right), \]
reduce to the Example \ref{Ex:Rational}.

In the case $m=1$, 
by substituting $\Phi_0$ in Example \ref{Ex:Rational} into 
\eqref{Form:Vcol},
\eqref{Form:Vrow}
we find that the Stokes matrices are just identity matrix,
\begin{align*}
\nu_{u,d}( \Phi_n(z) ) = I, \quad
S_d^\pm ( \Phi_n(z) ) = I.
\end{align*}
Therefore, even though $\Phi_n$ is a shrinking solution,
due to the resonant property, 
we still cannot distinguish it through 
the monodromy data mapping $\tmop{MD}_{u,d}$.

\subsection{Birkhoff's Reduction Problem}
\label{Sect:BRP}

The concepts occur in this section 
can be referred to Section \ref{subSect:BasicODE}.
As an application of 
Theorem \ref{Thm:termwise} and
Theorem \ref{Thm:Cat},
this section will give a new proof of 
certain Birkhoff’s reduction problem, which was a special case proved 
by Turrittin \cite{Turrittin1963}.
Finally, we will outline 
a possible sketch for the case of Poincaré $1$
with diagonalizable leading term,
which still remains open.
For the background of this problem, 
we refer to the survey \cite{Balser2004,Bolibruch2006,Bolibruch1995}.

Consider the following linear differential system 
with Poincaré rank $r$,
\begin{align}
\label{SysBRP}
\frac{\mathd}{\mathd z} X(z) & = A(z) X(z),\\
A(z) & = z^{r-1} \sum_{p=0}^\infty A_p z^{-p},
\nonumber
\end{align}
where $A(z)$ is convergent at $z=\infty$
and $A_0\neq 0$.
If there is an invertible transformation $T(z)$ 
such that $X(z)=T(z)Y(z)$, 
then $Y(z)$ will satisfy the following linear differential system
\begin{align*}
\frac{\mathd}{\mathd z} Y(z) & = B(z) Y(z),\\
B(z) & = T(z)^{-1} A(z) T(z) - T(z)^{-1} 
\frac{\mathd}{\mathd z} T(z),
\end{align*}
which naturally leads to the following definition.
\begin{defi}
If there exists a matrix function $T(z)$ meromorphic at $z=\infty$, 
such that
\begin{align*}
\frac{\mathd}{\mathd z} T(z) = A_1(z) T(z) - T(z) A_2(z),
\end{align*}
then the systems $A_1(z)$ and $A_2(z)$ are called \textbf{meromorphically equivalent}. 
If we further have $T(\infty)=I$, 
then the systems $A_1(z)$ and $A_2(z)$ are called \textbf{analytically equivalent}.
\end{defi}

Birkhoff’s reduction problem is 
\begin{itemize}

\item is the system \eqref{SysBRP} always 
meromorphically equivalent to 
a system with coefficients 
\begin{align}
\label{BSF}
B(z)=z^{r-1}\sum_{p=0}^r B_p z^{-p};
\end{align}

\item under what conditions can the system \eqref{SysBRP} be 
analytically equivalent to 
a system with coefficients \eqref{BSF},

\end{itemize}
and we call \eqref{BSF} 
the \textbf{Birkhoff standard form} of \eqref{SysBRP},
if it exists.

For the positive answer of the meromorphically equivalent case, 
in 1963 \cite{Turrittin1963}, Turrittin proved that 
$A_0$ needs to be diagonalizable and the eigenvalues are distinct;
according to Balser’s result \cite{Werner1989}, 
we can obtain a more general condition \cite{Balser2004}, that 
the formal monodromy matrix of 
the formal fundamental solution of system \eqref{SysBRP} 
needs to be diagonalizable.

For the positive answer of the analytically equivalent case, 
in 1913 \cite{Birkhoff1913}, Birkhoff proved that 
the monodromy matrix of 
the fundamental solution of system \eqref{SysBRP} 
needs to be diagonalizable; 
in 1994 \cite{Bolibrukh1994}, Bolibrukh proved that 
the system \eqref{SysBRP} needs to be irreducible.

We will focus on Birkhoff's reduction problem 
in the case of Poincaré rank $r=1$,
and $A_0$ is diagonalizable.
Let us introduce the following lemma.
\begin{lem}
\label{Lem:MeAnEq}
Denote the canonical fundamental solution of 
system \eqref{SysBRP} at $z=\infty$ as $X_d(z;A)$. Assume that we have two systems with coefficients
\begin{align*}
A^{(1)}(z) = \sum_{p=0}^\infty A_{p}^{(1)} z^{-p},\quad
A^{(2)}(z) = \sum_{p=0}^\infty A_{p}^{(2)} z^{-p},
\end{align*}
and $u=A_{0}^{(1)}=A_{0}^{(2)}$ is a diagonal matrix. Then
\begin{itemize}

\item[(1).] they are meromorphically equivalent
if and only if 
$X_d(z;A^{(1)})$ and $X_d(z;A^{(2)})$ have the same monodromy factors;

\item[(2).] they are analytically equivalent
if and only if 
$X_d(z;A^{(1)})$ and $X_d(z;A^{(2)})$ have the same monodromy factors, and
$\delta_u A_{1}^{(1)}=\delta_u A_{1}^{(2)}$.
\end{itemize}
\end{lem}

\begin{prf}
Note that $X_d(z;A^{(i)})$ is asymptotic to
$F_d(z;A^{(i)}) \mathe^{u z} z^{ \delta_u A_{1}^{(i)} }$,
where \begin{align}
\label{lemmeroeq}
F_d(z;A^{(i)}) = I + O(z^{-1}),\quad
z\to\infty, \arg z = d.
\end{align}
Denote $X(z)=X_d(z;A^{(1)})^{-1} X_d(z;A^{(2)})$,
we have
\begin{align*}
\text{Systems with coefficients $A^{(1)},A^{(2)}$ are meromorphically equivalent }
& \Leftrightarrow
X(z) \text{ is meromorphic at $z=\infty$}\\
& \phantom{\Leftrightarrow}
\qquad\qquad \Updownarrow\\
\text{$X_d(z;A^{(1)}),X_d(z;A^{(2)})$ 
have the same monodromy factors}
& \Leftrightarrow
X(z) \text{ is single-valued at $z=\infty$}
\end{align*}
If $X(z)$ is meromorphic at $z=\infty$, 
and $\delta_u A_{1}^{(1)}=\delta_u A_{1}^{(2)}$, 
according to \eqref{lemmeroeq} 
we have $X(\infty)=I$, thus we complete the proof.
\end{prf}

If we require the eigenvalues of $u$ are distinct, 
then the meromorphic Birkhoff reduction problem can be reduced to 
whether we can find a $\Phi_0\in\mathfrak{c}_0$ 
such that $\nu_{\tmop{cat}(0),-\frac{\pi}{2}}(\Phi_0)=V$,
for any given matrix $V$.
If we additionally take 
the diagonal matrix $\Lambda = \tmop{diag}(a_{11},\ldots,a_{nn})$
such that $\mathe^{2\pi\mathi\sum_{j=1}^k a_{jj}}=
\tmop{det}V^{[k]}$, 
then the analytic Birkhoff reduction problem can be reduced to 
whether we can take find a $\Phi_0\in\mathfrak{c}_0$ 
such that $\nu_{\tmop{cat}(0),-\frac{\pi}{2}}(\Phi_0)=V$
and $\delta \Phi_0 = \Lambda$.

The following proposition is a
special case of \cite{Turrittin1963}.

\begin{pro}
Suppose that
$A(z)=\sum_{p=0}^\infty A_p z^{-p}$,
and $A_0$ is diagonal with distinct eigenvalues, then for fixed $d\in\mathbb{R}$, 
the monodromy factor $V$ of the
canonical fundamental solution $X_d(z;A)$
can take any invertible matrix. 
As a consequence, 
the system $A(z)$ is
meromorphically equivalent to its
Birkhoff standard form.
\end{pro}

\begin{prf}
Under the Definition \ref{Def:vcat},
according to Lemma \ref{Lem:MeAnEq} and Theorem \ref{Thm:Cat},
we only need to prove that 
for any given invertible matrix $V$,
there is a $\Phi_0\in\mathfrak{c}_0$ 
such that $\nu_{\tmop{cat}(0),-\frac{\pi}{2}}(\Phi_0)=V$.
For the eigenvalues set $\sigma_k$
of the submatrix $\Phi_0^{[k]}$,
from \eqref{Ck-1rec}
we have $\mathe^{2\pi\mathi \sigma_k } =
\sigma(V^{[k]})$.
Note that we can always choose $\sigma_1,\ldots,\sigma_n$
satisfy condition \eqref{Cond:NonResNear},
according to \eqref{Form:Vcol}, \eqref{Form:Vrow},
we can construct the $\Phi_0\in\mathfrak{c}_0$ that we need.
\end{prf}

For the case where the eigenvalues of 
the diagonal matrix $u=A_0$ have multiplicities, 
we can also establish the explicit Riemann-Hilbert mapping 
similar to Theorem \ref{Thm:Cat}. 
If we want to consider 
the meromorphic Birkhoff reduction problem in this case, 
we will need to investigate the properties of 
the corresponding Stokes matrix 
$S_d^\pm(E_k+\cdots+E_n,A)$,
and this case will be discussed in future work.


\begin{appendices}
   \section{Explicit Stokes Matrices via Boundary Values}
\label{Sect:Diag}

Suppose that $\Phi_n(z)$ is a shrinking solution
with boundary value $\Phi_0$.
In this appendix, 
we will introduce 
the explicit expression of $S^\pm_d(u,\Phi_n(z))$ 
given by the boundary value $\Phi_0$, 
based on Theorem \ref{Thm:Cat}.
We will specifically provide 
the explicit expressions for 
the entries.

We need to introduce the explicit expressions of 
the connection matrices $C_d(u,A)$
for confluent hypergeometric systems \eqref{Conflu}
with coefficients $(u,A)=(E_k,\delta_k A)$
under Notation \ref{Nota:Mat}.

\begin{pro}
\label{Pro:O}
  Suppose that 
  $A^{[k]}, A^{[k - 1]}$ 
  are non-resonant and ${q} \in
  \mathbb{Z}$. We have
  \begin{align*}
    C_d (E_k, A^{[k]}) 
    & = 
    \left(\begin{array}{c}
      \mathe^{(2{q} \pm 1) \pi \mathi (\tmmathbf{r}^{[k]} -\tmmathbf{l}^{[k -
      1]})} \omega_{\hat{k} \ast}\\
      \mathe^{2{q} \pi \mathi (\tmmathbf{r}^{[k]} - a_{k k})} \omega_{k \ast}
    \end{array}\right), \\
    C_d (E_k, A^{[k]})^{- 1} 
    & = 
    \left(\begin{array}{cc}
      \mathe^{2{q} \pi \mathi (\tmmathbf{r}^{[k - 1]} -\tmmathbf{l}^{[k]})}
      m_{\ast \hat{k}} & \mathe^{(2{q} \pm 1) \pi \mathi (a_{k k} I
      -\tmmathbf{l}^{[k]})} m_{\ast k}
    \end{array}\right) 
  \end{align*}
  where
  \begin{itemize}
    \item $| d - 2{q} \pi |, | d - (2{q} \pm 1) \pi | < \pi$;
    
    \item $\tmmathbf{l}^{[j]}$ is the operator of left multiplication by
    $A^{[j]}$;
    
    \item $\tmmathbf{r}^{[j]}$ is the operator of right multiplication by
    $A^{[j]}$.
  \end{itemize}
  \begin{align*}
    \omega_{\hat{k} \ast} 
    & = 
    \prod_{p, j} \left( (\tmmathbf{l}^{[k - 1]}
    -\tmmathbf{r}^{[k]}) ! \frac{(\tmmathbf{l}^{[k - 1]} - \lambda^{(k -
    1)}_p) !}{(\tmmathbf{l}^{[k - 1]} - \lambda^{(k)}_j) !}
    (\tmmathbf{r}^{[k]} -\tmmathbf{l}^{[k - 1]}) ! \frac{(\tmmathbf{r}^{[k]} -
    \lambda^{(k)}_j) !}{(\tmmathbf{r}^{[k]} - \lambda^{(k - 1)}_p) !} \right)
    \left(\begin{array}{cc}
      {I}_{k - 1} & 0
    \end{array}\right), \\
    \omega_{k \ast} 
    & = 
    \left( \prod_{p, j} \frac{(\tmmathbf{r}^{[k]} -
    \lambda^{(k)}_j) !}{(\tmmathbf{r}^{[k]} - \lambda^{(k - 1)}_p) !} \right)
    \left(\begin{array}{cc}
      O_{1 \times (k - 1)} & 1
    \end{array}\right), \\
    m_{\ast \hat{k}} 
    & = 
    \prod_{p, j} \left( (\tmmathbf{l}^{[k]}
    -\tmmathbf{r}^{[k - 1]}) ! \frac{(\lambda^{(k - 1)}_p -\tmmathbf{r}^{[k -
    1]}) !}{(\lambda^{(k)}_j -\tmmathbf{r}^{[k - 1]}) !} (\tmmathbf{r}^{[k -
    1]} -\tmmathbf{l}^{[k]}) ! \frac{(\lambda^{(k)}_j -\tmmathbf{l}^{[k]})
    !}{(\lambda^{(k - 1)}_p -\tmmathbf{l}^{[k]}) !} \right)
    \left(\begin{array}{c}
      {I}_{k - 1}\\
      0
    \end{array}\right), \\
    m_{\ast k} 
    & = 
    \left( \prod_{p, j} \frac{(\lambda^{(k)}_j
    -\tmmathbf{l}^{[k]}) !}{(\lambda^{(k - 1)}_p -\tmmathbf{l}^{[k]}) !}
    \right) \left(\begin{array}{c}
      O_{(k - 1) \times 1}\\
      1
    \end{array}\right) . 
  \end{align*}
We denote $x!\assign \Gamma(x+1)$.
\end{pro}

Given any direction $d$, 
such that $d\notin \tmop{aS}(u)$.
Let us take 
the following connected components $\mathcal{R}_d(J)$ 
labelled by 
the subsets $J\subseteq \{1,\ldots,n-1\}$: 
$u\in \mathcal{R}_d(J)$ if
\begin{align}
\tmop{Im}(u_{k+1} \mathe^{\mathi d}) &<
\min_{1\leqslant j \leqslant k} 
\tmop{Im}(u_j \mathe^{\mathi d}), \
\text{ for } k\in J,\\
\tmop{Im}(u_{k+1} \mathe^{\mathi d}) &>
\max_{1\leqslant j \leqslant k} 
\tmop{Im}(u_j \mathe^{\mathi d}), \ \text{ for } 
k\notin J.
\end{align}
\begin{defi}
For the $(k+1)\times(k+1)$ matrix 
$V=\begin{pmatrix}
V^{[k]} & \alpha \\
\beta & v_{k+1,k+1}
\end{pmatrix}$,
define
\begin{align*}
\mathcal{p}_{+}^{(k)} V \assign
V^{[k]},\quad
\mathcal{p}_{-}^{(k)} V \assign
((V^{-1})^{[k]})^{-1} =
V^{[k]} - \frac{\alpha \beta}{v_{k+1,k+1}}.
\end{align*}
For $u\in\mathcal{R}_d(J)$,
we recursively define
\begin{align*}
\mathcal{p}_{u,d}^{(n)}  \assign {\rm Id},\quad
\mathcal{p}_{u,d}^{(k)}  \assign 
\left\{
\begin{array}{ll}
     \mathcal{p}^{(k)}_+ \mathcal{p}^{(k + 1)}_{u, d} & ; k \in J\\
     \mathcal{p}^{(k)}_- \mathcal{p}^{(k + 1)}_{u, d} & ; k \nin J
\end{array}\right..
\end{align*}
\end{defi}

\begin{lem}
\label{Lem:AppCatReal}
Suppose that $u\in\mathcal{R}_d(J)$,
for $1\leqslant k \leqslant n$ we have
\begin{align}
\label{PreAppCatReal}
\mathcal{p}^{(k)}_{u, d}
\nu_{u,d}(\Phi_n(z))
& =
\tmop{Ad}\left(
\overrightarrow{\prod^{k}_{j=1}}
C_{d+\arg(u_j-u_{j-1})}(E_j,\delta_j\Phi_0^{[k]})
\right)
\mathe^{2\pi\mathi\Phi_0^{[k]}},
\end{align}
\end{lem}

\begin{prf}
Let us inductively assume that
\eqref{PreAppCatReal}
holds for $k+1$,
and denote $\theta_j=\arg(u_j-u_{j-1})$,
therefore
\begin{align*}
\mathcal{p}_{u,d}^{(k+1)}
\nu_{u,d}(\Phi_n(z))
& =
\tmop{Ad}\left(
\overrightarrow{\prod^{k+1}_{j=1}}
C_{d+\theta_j}(E_j,\delta_j\Phi_0^{[k+1]})
\right)
\mathe^{2\pi\mathi\Phi_0^{[k+1]}}\\
& = \tmop{Ad} 
\left(\begin{array}{cc}
\overrightarrow{\prod^{k}_{j=1}}
C_{d+\theta_j}(E_j,\delta_j\Phi_0^{[k]}) & 0\\
       0 & 1
\end{array}\right)\\
&\phantom{==}
\left(
S_{d+\theta_{k+1}}^-(E_{k+1},\Phi_0^{[k + 1]})^{- 1} 
\left(\begin{array}{cc}
\mathe^{2\pi\mathi\Phi_0^{[k]}} & 0\\
0 & \mathe^{2\pi\mathi a_{k+1,k+1}}
\end{array}\right)
S_{d+\theta_{k+1}}^+(E_{k+1},\Phi_0^{[k + 1]})
\right).
\end{align*}
Note that
\begin{align*}
\tmop{Im}(u_{k+1} \mathe^{\mathi d}) <
\tmop{Im}(u_k \mathe^{\mathi d})
&\Leftrightarrow
S_{d+\theta_{k+1}}^+(E_{k+1},\Phi_0^{[k + 1]})=
\left(\begin{array}{cc}
I & \ast\\
0 & 1
\end{array}\right),\
S_{d+\theta_{k+1}}^-(E_{k+1},\Phi_0^{[k + 1]})=
\left(\begin{array}{cc}
I & 0\\
\ast & 1
\end{array}\right),\\
\tmop{Im}(u_{k+1} \mathe^{\mathi d}) >
\tmop{Im}(u_k \mathe^{\mathi d})
&\Leftrightarrow
S_{d+\theta_{k+1}}^+(E_{k+1},\Phi_0^{[k + 1]})=
\left(\begin{array}{cc}
I & 0\\
\ast & 1
\end{array}\right),\
S_{d+\theta_{k+1}}^-(E_{k+1},\Phi_0^{[k + 1]})=
\left(\begin{array}{cc}
I & \ast\\
0 & 1
\end{array}\right),
\end{align*}
therefore 
it can be directly verified that
\eqref{PreAppCatReal}
holds for $k$, thus for any $1\leqslant k \leqslant n$. 
\end{prf}

Denote $(\lambda^{(k)}_j)_{j = 1, \ldots, k}$
as the eigenvalues of $\Phi_0^{[k]}$.
\begin{thm}\cite{Xu}
\label{Thm:AppCatReal}
Suppose that $\Phi_n(z)= (\varphi_{ij})_{n \times n}$ is 
the shrinking solution of 
the $n$-th isomonodromy equation \eqref{zisoeq} with 
the boundary value $\Phi_0$ and $u\in\mathcal{R}_d(J)$.
For every $1\leqslant k \leqslant n-1$,
we additionally require that
the upper left $k \times k$ submatrix $\Phi_0^{[k]}$ of $\Phi_0$ 
is non-resonant, and take 
$q^{\rm even}_i\in 2\mathbb{Z},
 q^{\rm odd}_i\in 2\mathbb{Z}-1$ such that
\begin{align*}
|d+\arg(u_{i+1}-u_i)-q^{\rm even}_i\pi|
& <\pi,\\
|d+\arg(u_{i+1}-u_i)-q^{\rm odd}_i\pi|
& <\pi,\quad
i=0,1,\ldots,n-1.
\end{align*}
Then we have
\begin{align*}
S_d (u, \Phi_n(z))_{k, k + 1} 
& = 
\mathe^{(q^{\rm even}_k\varphi_{k+1, k+1}-
(q^{\rm odd}_{k-1}-
 q^{\rm odd}_{k}+
 q^{\rm even}_{k})
\varphi_{k k})
\pi\mathi} 
\cdot s_{k, k+1},\\
S_d (u,\Phi_n(z))_{k+1,k} 
& = 
\mathe^{(q^{\rm odd}_{k-1}\varphi_{k k}
-q^{\rm odd}_{k}\varphi_{k+1, k+1})
\pi\mathi}
\cdot s_{k + 1, k},
\end{align*}
where $S_d=S_d^+-S_d^-$, and
\begin{align}
\label{Cat:skk+1}
s_{k, k+1} 
& = 
2 \pi \mathi \cdot 
\left( \prod_{p, j_1, j_2}
\frac{(\Phi_0^{[k]} - \lambda^{(k)}_p {I})! 
(\Phi_0^{[k]} - \lambda^{(k)}_p{I}) !}
{(\Phi_0^{[k]} - \lambda^{(k + 1)}_{j_1} {I})! 
(\Phi_0^{[k]} -\lambda^{(k - 1)}_{j_2} {I})!} \right)_{k \ast}
\mathe^{(q^{\rm odd}_{k-1}-q^{\rm odd}_k)\pi\mathi\Phi_0^{[k]}}
(\Phi_0)^{[k+1]}_{\widehat{k + 1}, k + 1},\\
\label{Cat:sk+1k}
s_{k+1, k} 
& = 
2 \pi \mathi \cdot (\Phi_0)^{[k + 1]}_{k + 1,\widehat{k + 1}}
\mathe^{(q^{\rm odd}_{k}-q^{\rm odd}_{k-1})
\pi\mathi\Phi_0^{[k]}}
\left( \prod_{p, j_1, j_2} \frac{(\lambda^{(k)}_p {I} -\Phi_0^{[k]})! 
(\lambda^{(k)}_p {I} - \Phi_0^{[k]}) !}
{(\lambda^{(k+1)}_{j_1} {I}- \Phi_0^{[k]})! 
 (\lambda^{(k-1)}_{j_2} {I}- \Phi_0^{[k]})!} 
\right)_{\ast k}.
\end{align}
\end{thm}

\begin{prf}
For $1 \leqslant k \leqslant n-1$,
denote $\theta_k=\arg(u_k-u_{k-1})$,
and $S_d^\pm\assign S_d^\pm(u,\Phi_n)$
for short.
Note that for $u\in\mathcal{R}_d(J)$ we have
$((S_d^\pm)^{-1})^{[k]}=
((S_d^\pm)^{[k]})^{-1}$, and
\begin{align*}
\mathcal{p}_{u,d}^{(k+1)}
\nu_{u,d}(\Phi_n(z))
& = 
((S_d^-)^{-1} )^{[k+1]}
\cdot 
\mathe^{2\pi\mathi\delta\Phi_0^{[k+1]}}
\cdot 
(S_d^+)^{[k+1]},\\
\tmop{Ad} \left(
\overrightarrow{\prod^{k}_{j=1}}
C_{d+\theta_j}(E_j,\delta_j\Phi_0^{[k]})
\right)
\mathe^{2 \pi \mathi \Phi_0^{[k]}} 
& = 
((S_d^-)^{-1} )^{[k]} \cdot 
\mathe^{2\pi\mathi\delta\Phi_0^{[k]}}
\cdot 
(S_d^+)^{[k]}.
\end{align*}
By the uniqueness of the LU decomposition, 
Lemma \ref{Lem:AppCatReal} and
$(S_d^{\pm})^{-1}=S_{d\pm\pi}^{\mp}$, 
we have
\begin{align*}
(S_d^{\pm})^{[k+1]}
& = 
\left(\begin{array}{cc}
(S_d^{\pm})^{[k]} & 0\\
0 & 1
\end{array}\right) 
\left( 
\tmop{Ad} \left(
\overrightarrow{\prod^{k}_{j=1}}
C_{d+\theta_j}(E_j,\delta_j\Phi_0^{[k+1]})
\right)
S_{d+\theta_{k+1}}^{\pm} (E_{k + 1},\Phi_0^{[k+1]})\right)\\
& =
\left( \tmop{Ad} \left(
\overrightarrow{\prod^{k}_{j=1}}
C_{d\pm\pi+\theta_j}(E_j,\delta_j\Phi_0^{[k+1]})
\right)
S_{d+\theta_{k+1}}^{\pm} (E_{k + 1},\Phi_0^{[k+1]})\right)
\left(\begin{array}{cc}
(S_d^{\pm})^{[k]} & 0\\
0 & 1
\end{array}\right) .\\
\end{align*}
For $1\leqslant j \leqslant k$,
if we denote $j\times j$ matrix
\begin{align*}
C_{d\pm\pi+\theta_j} (E_j, \Phi^{[j]})
= \left(
\begin{array}{c}
\Tilde{\omega}_{\widehat{j}} \\
\Tilde{\omega}_{j}
\end{array}\right), \quad
C_{d+\theta_j} (E_j, \Phi^{[j]})^{- 1} 
= \left(
\begin{array}{cc}
\Tilde{m}_{\widehat{j}} & 
\Tilde{m}_{j}
\end{array}\right),
\end{align*}
then we have
\begin{align}
\label{Form:jkentry}
(S_d^{\pm})_{j+1,k+1}
& =\Tilde{\omega}_{j+1}
\Tilde{\omega}_{\widehat{j+2}}\cdots
\Tilde{\omega}_{\widehat{k}}\cdot
S_{d+\theta_{k+1}}^{\pm} 
(E_{k + 1},\Phi_0^{[k+1]})_{\widehat{k+1},k+1},\\
\label{Form:kjentry}
(S_d^{\pm})_{k+1,j+1}
& =
S_{d+\theta_{k+1}}^{\pm}(E_{k+1},\Phi_0^{[k+1]})_{k+1,\widehat{k+1}}\cdot
\Tilde{m}_{\widehat{k}}\cdots
\Tilde{m}_{\widehat{j+2}}
\Tilde{m}_{j+1}.
\end{align}
It concludes the proof.
\end{prf}
\begin{rmk}
Note that we have
\begin{align*}
\tmop{Im}(u_{k} \mathe^{\mathi d})<
\tmop{Im}(u_{j} \mathe^{\mathi d})
& \Rightarrow
(S_d^+)_{k j}=(S_d^-)_{j k}=0,\\
\tmop{Im}(u_{k} \mathe^{\mathi d})>
\tmop{Im}(u_{j} \mathe^{\mathi d})
& \Rightarrow
(S_d^+)_{j k}=(S_d^-)_{k j}=0,
\end{align*}
therefore we can directly recover the Stokes matrix $S_d^\pm$ from $S_d\assign S_d^+-S_d^-$.
Denote
\begin{align*}
\tmop{sgn}_J (k) \assign 
\left\{\begin{array}{ll}
1 & ; k \in J\\
- 1 & ; k \nin J
\end{array}\right.,
\end{align*}
we can verify that
$q^{\rm odd}_k = q^{\rm even}_k - \sgn_J(k)$.
\end{rmk}

When $\Phi_0^{[k]}$ can be diagonalized, 
\eqref{Cat:skk+1} and \eqref{Cat:sk+1k} 
can be written in the more explicit form given in Theorem \ref{thm: introcatformula}. It relies on the following diagonalization formula, see e.g., \cite{Denton2022}.

\begin{pro}
\label{For:Diag}
Suppose that $\lambda^{(k)}_1, \ldots, \lambda^{(k)}_k$ are distinct, and $\{ \lambda^{(k - 1)}_1, \ldots, \lambda^{(k - 1)}_{k - 1} \} \cap
    \{ \lambda^{(k)}_1, \ldots, \lambda^{(k)}_k \} = \varnothing$.
For arbitrary $d_1, \ldots, d_k \neq 0$, 
we have
$P=(P_{r c})_{n\times n}$ and
$P^{-1}=((P^{-1})_{r c})_{n\times n}$
with entries
\begin{align*}
P_{r c}  = 
(- 1)^{r + k}
\frac{\det
(\lambda^{(k)}_c {I} - A)
^{1, \ldots, k\backslash k}
_{1, \ldots, k\backslash r}}
{\prod_{j = 1}^{k - 1}
(\lambda^{(k)}_c - \lambda^{(k - 1)}_j)} d_c,\quad
(P^{- 1})_{r c}  =
(- 1)^{k + c} \frac{
\det
(\lambda^{(k)}_r {I} - A)
^{1, \ldots, k\backslash c}
_{1, \ldots, k\backslash k}
}{\prod_{p \neq r}
(\lambda_r^{(k)} - \lambda^{(k)}_p)} d_r^{- 1},
\end{align*}
  such that
\begin{align*}
P^{- 1} A^{[k]} P & = 
\tmop{diag} (\lambda^{(k)}_1, \ldots,\lambda^{(k)}_k).
\end{align*}
We denote $X
^{1, \ldots, k\backslash a}
_{1, \ldots, k\backslash b}$
as the $(k-1)\times(k-1)$ submatrix of $X$, 
which is formed by the first $k$ rows 
excluding the $a$-th row, 
and the first $k$ columns 
excluding the $b$-th column.
\end{pro}

\begin{pro}
\label{For:DiagEnt}
Suppose that the 
$k\times k$ matrix $P$ is defined in 
the Proposition \ref{For:Diag}. Denote
\begin{align*}
\left(\begin{array}{cc}
      P^{- 1} & 0\\
      0 & {I}_{n - k}
    \end{array}\right) A 
\left(\begin{array}{cc}
      P & 0\\
      0 & {I}_{n - k}
    \end{array}\right) 
    & =
\left(\begin{array}{cccccc}
      \lambda^{(k)}_1 &  & 0 & \alpha^{(k)}_{1, k + 1} & \cdots &
      \alpha^{(k)}_{1, n}\\
      & \ddots &  & \vdots & \ddots & \vdots\\
      0 &  & \lambda^{(k)}_k & \alpha^{(k)}_{k, k + 1} & \cdots &
      \alpha^{(k)}_{k, n}\\
      \beta^{(k)}_{k + 1, 1} & \cdots & \beta^{(k)}_{k + 1, k} & a_{k + 1, k +
      1} & \cdots & a_{k + 1, n}\\
      \vdots & \ddots & \vdots & \vdots & \ddots & \vdots\\
      \beta^{(k)}_{n, 1} & \cdots & \beta^{(k)}_{n, k} & a_{n, k + 1} & \cdots
      & a_{n n}
\end{array}\right),
\end{align*}
  we have
\begin{align*}
\alpha^{(k)}_{r, c} 
 = 
- \frac{\det
(\lambda^{(k)}_r {I} - A)
^{1,\ldots,k-1, k}
_{1, \ldots, k - 1, c}}
{\prod_{p \neq r} 
(\lambda_r^{(k)} - \lambda^{(k)}_p)} 
d_r^{- 1},\quad 
\beta^{(k)}_{r, c} = 
- \frac{\det 
(\lambda^{(k)}_c {I} - A)
^{1,\ldots,k-1,r}
_{1,\ldots,k-1,k}}
{\prod_{j = 1}^{k - 1} (\lambda^{(k)}_c -
    \lambda^{(k - 1)}_j)} d_c .
\end{align*}
We denote $X
^{a_1, \ldots, a_k}
_{b_1, \ldots, b_k}$
as the $k\times k$ submatrix of $X$, 
which is formed by the $(a_1, \ldots, a_k)$ rows 
and the $(b_1, \ldots, b_k)$ columns.
\end{pro}
Now we can give the explicit expression of the Stokes matrices via the boundary values. 
\begin{pro}
Suppose that $u\in\mathcal{R}_d(J)$,
take 
$q^{\rm even}_i\in 2\mathbb{Z},
 q^{\rm odd}_i\in 2\mathbb{Z}-1$ such that
\begin{align*}
|d+\arg(u_{i+1}-u_i)-q^{\rm even}_i\pi|
& <\pi,\\
|d+\arg(u_{i+1}-u_i)-q^{\rm odd}_i\pi|
& <\pi,\quad
i=0,1,\ldots,n-1.
\end{align*}
and $j<k,\Phi_0^{[j]},\ldots,\Phi_0^{[k-1]}$ 
is diagonalizable.
For $S_d\assign S_d^+-S_d^-$
we have
\begin{align*}
S_d(u,\Phi_n(z))_{j k} & =
\mathe^{(q_{k - 1}^{\tmop{even}} \varphi_{k k} - 
(q_{j - 1}^{\tmop{odd}} -
q_{k - 1}^{\tmop{odd}} + 
q_{k - 1}^{\tmop{even}}) 
\varphi_{j j}) \pi\mathi}
\cdot s_{j k},\\
S_d(u,\Phi_n(z))_{k j} & =
\mathe^{(q_{j - 1}^{\tmop{odd}} \varphi_{j j} - 
q_{k - 1}^{\tmop{odd}}
\varphi_{k k}) \pi\mathi}
\cdot s_{k j},
\end{align*}
\begin{align*}
s_{j k} & =
- 2 \pi \mathi 
\sum_{i_j, \ldots, i_{k - 1}} 
\mathe^{\left(
(q_{j-1}^{\tmop{odd}} - q_{j - 1}^{\tmop{even}}) \lambda^{(j)}_{i_j} + \left(
   \sum_{s = j}^{k - 1} (q_{s - 1}^{\tmop{even}} - q_s^{\tmop{even}})
   \lambda^{(s)}_{i_s} \right) + (q_{k - 1}^{\tmop{even}} - q_{k -
   1}^{\tmop{odd}}) \lambda^{(k - 1)}_{i_{k - 1}} 
\right) \pi \mathi}\\
&\quad\
\frac{\prod_{s = j}^{k - 1} ((\lambda^{(s)}_{i_s}
   -\tmmathbf{\lambda}^{(s)}_{\ast}) ! (\lambda^{(s)}_{i_s}
   -\tmmathbf{\lambda}^{(s)}_{\hat{i}_s} - 1) ! \cdot \det
   (\lambda^{(s)}_{i_s} I - \Phi_0)^{1, \ldots, s - 1, s}_{1, \ldots, s - 1, s
   + 1})}{(\lambda^{(j)}_{i_j} -\tmmathbf{\lambda}^{(j - 1)}_{\ast}) ! \left(
   \prod_{s = j}^{k - 2} (\lambda^{(s)}_{i_s} -\tmmathbf{\lambda}^{(s +
   1)}_{\hat{i}_{s + 1}}) ! (\lambda^{(s + 1)}_{i_{s + 1}}
   -\tmmathbf{\lambda}^{(s)}_{\hat{i}_s}) ! \cdot (\lambda^{(s)}_{i_s} -
   \lambda^{(s + 1)}_{i_{s + 1}}) \right) (\lambda^{(k - 1)}_{i_{k - 1}}
   -\tmmathbf{\lambda}^{(k)}_{\ast}) !}
\\
s_{k j} & = 
2 \pi \mathi \sum_{i_j, \ldots, i_{k - 1}} 
\mathe^{\left( 
(q_{j-1}^{\tmop{even}} - q_{j - 1}^{\tmop{odd}}) \lambda^{(j)}_{i_j} + 
\left(\sum_{s = j}^{k - 1} (q_s^{\tmop{even}} - q_{s - 1}^{\tmop{even}})
   \lambda^{(s)}_{i_s} \right) + 
(q_{k - 1}^{\tmop{odd}} - 
 q_{k - 1}^{\tmop{even}}) 
\lambda^{(k - 1)}_{i_{k - 1}} \right) \pi \mathi}\\
&\quad\
\frac{\prod_{s = j}^{k - 1} ((\tmmathbf{\lambda}^{(s)}_{\ast} -
   \lambda^{(s)}_{i_s}) ! (\tmmathbf{\lambda}^{(s)}_{\hat{i}_s} -
   \lambda^{(s)}_{i_s} - 1) ! \cdot \det (\Phi_0 - \lambda^{(s)}_{i_s} I)^{1,
   \ldots, s - 1, s + 1}_{1, \ldots, s - 1, s})}{(\tmmathbf{\lambda}^{(j -
   1)}_{\ast} - \lambda^{(j)}_{i_j}) ! \left( \prod_{s = j}^{k - 2}
   (\tmmathbf{\lambda}^{(s + 1)}_{\hat{i}_{s + 1}} - \lambda^{(s)}_{i_s}) !
   (\tmmathbf{\lambda}^{(s)}_{\hat{i}_s} - \lambda^{(s + 1)}_{i_{s + 1}}) !
   \cdot (\lambda^{(s + 1)}_{i_{s + 1}} - \lambda^{(s)}_{i_s}) \right)
   (\tmmathbf{\lambda}^{(k)}_{\ast} - \lambda^{(k - 1)}_{i_{k - 1}}) ! }.
\end{align*}
Here we denote
$f(\tmmathbf{\lambda}^{(s)}_{\ast})\assign
\prod_{i=1}^{s}f(\lambda^{(s)}_{i})$,
$f(\tmmathbf{\lambda}^{(s)}_{\hat{i}_s})\assign
\prod_{i\neq i_s} f(\lambda^{(s)}_{i})$, for example, $(\tmmathbf{\lambda}^{(s)}_{\ast} -
\lambda^{(s)}_{i_s}) ! =
\prod_{i=1}^s ({\lambda}^{(s)}_i -
   \lambda^{(s)}_{i_s}) ! $.
\end{pro}

\begin{prf}
For simplicity, 
we only provide the case of $(k, k+1)$-entry.
The remaining cases can be directly obtained from 
\eqref{Form:jkentry} and \eqref{Form:kjentry}
in a similar way.
Suppose that
$P^{-1} \Phi_0^{[k]} P =
\diag(\lambda^{k}_1,\ldots,\lambda^{k}_k)$,
denote
\begin{align*}
D_k =
\diag(\lambda^{k}_1,\ldots,\lambda^{k}_k),\quad
\left(\begin{array}{cc}
      P^{- 1} & 0\\
      0 & {I}_{n - k}
      \end{array}\right)
\Phi_0
\left(\begin{array}{cc}
      P & 0\\
      0 & {I}_{n - k}
      \end{array}\right)
= A_k.
\end{align*}
We can rewrite
\eqref{Cat:skk+1} as
\begin{align*}
s_{k, k + 1} 
& = 
2 \pi \mathi \cdot 
P_{k\ast}
\left( \prod_{p, j_1, j_2}
    \frac{(D_k - \lambda^{(k)}_p I)! 
    (D_k - \lambda^{(k)}_p I)!}
    {(D_k -\lambda^{(k + 1)}_{j_1} I)! 
     (D_k -\lambda^{(k - 1)}_{j_2} I)!} \right)
    (A_k)^{[k +
    1]}_{\widehat{k + 1}, k + 1}.
\end{align*}
From Proposition \ref{For:Diag} and 
Proposition \ref{For:DiagEnt}
we have
\begin{align*}
P_{k i}\cdot(A_k)_{i,k+1} &=
- \frac{\det
(\lambda^{(k)}_i {I}-\Phi_0)
^{1, \ldots,k-1, k}
_{1, \ldots,k-1, k+1}}
{\prod_{p \neq i} (\lambda_i^{(k)} -
    \lambda^{(k)}_p)}.
\end{align*}
Thus we finish the proof.
\end{prf}

\begin{rmk}\label{diffconv}
Our convention of Stokes matrices
are different from \cite{Xu}. In particular, according to the convention in \cite{Xu}, the Stokes matrices of the equation \ref{Conflu} with respect 
to the direction $d=-\frac{\pi}{2}$ are
\[S_+( u,A)=\mathe^{-\pi\mathi\delta A}
S_{-\frac{\pi}{2}}^{-}(\mathi u,A), \ \ S_-(u,A)=(\mathe^{\pi\mathi\delta A}
S_{-\frac{\pi}{2}}^{+}(\mathi u,A))^{-1},\]
where $S_{-\frac{\pi}{2}}^{\pm}(\mathi u,A)$ are the Stokes matrices as in our convention, i.e., Definition \ref{Def:MonoData} with $d=-\frac{\pi}{2}$ 
and $u$ replacing by $\mathi u$. Taking the difference of notations into account, 
and further requiring
\[
{\rm arg}(u_{2}-u_1)=\cdots
={\rm arg}(u_{n}-u_{n-1})=\frac{\pi}{2},
\]
then the case $d=-\frac{\pi}{2}$ of the above theorem becomes Theorem \ref{thm: introcatformula}.    
\end{rmk}

\end{appendices}

\normalem
\bibliography{paper.bib}
\bibliographystyle{alpha}

\Addresses
\end{document}